\documentclass{amsart}
\usepackage{amssymb}
\usepackage{amsmath}
\usepackage{longtable}

\usepackage[pdftex,colorlinks,hypertexnames=false,
            linkcolor=blue,urlcolor=blue,citecolor=blue]{hyperref}

\allowdisplaybreaks

\newtheorem{thm}{Theorem}[section]
\newtheorem{lem}[thm]{Lemma}
\newtheorem{prop}[thm]{Proposition}
\newtheorem{cor}[thm]{Corollary}
\newtheorem{rem}[thm]{Remark}
\newtheorem{defn}[thm]{Definition}
\newtheorem{hypo}[thm]{Hypothesis}

\newcommand{\Z}{{\mathbb{Z}}}

\newcommand{\R}{{\mathbb{R}}}
\newcommand{\C}{{\mathbb{C}}}
\newcommand{\F}{{\mathbb{F}}}

\newcommand{\Irr}{\operatorname{Irr}}

\newcommand{\la}{\langle}
\newcommand{\ra}{\rangle}

\def\pr {{\bf Proof :}~}
\def\epr{\hfill {\bf QED}}

\def\Irr{{\rm Irr}}
\def\Ker{{\rm Ker}}

\def\GL{{\rm GL}}
\def\SL{{\rm SL}}

\def\rk{{\rm rk}}
\def\rz{{\rm rz}}
\def\rs{{\rm rs}}
\def\rKer{{\rm rKer}}
\def\rZ{{\rm rZ}}
\def\rS{{\rm rS}}
\def\height{{\rm ht}}
\def\ob{{\rm ob}}
\def\kob{{\rm kob}}
\def\lin{{\rm lin}}
\def\mida{{\rm mida}}
\def\norml{{\rm normal}}
\def\infl{{\rm Infl}}

\def\a{{\alpha}}
\def\b{{\beta}}

\begin{document}

\title[On the characters of the Sylow $p$-subgroups of $Y_n(p^a)$]{On the
       characters of the Sylow $\mathbf{p}$-subgroups of untwisted Chevalley groups $\mathbf{Y_n(p^a)}$} 

\author{Frank Himstedt, Tung Le and Kay Magaard}

\address{F.H.: Technische Universit\"at M\"unchen, Zentrum Mathematik --
         M11, Boltzmannstr. 3, 85748 Garching, Germany}
\address{T.L.: North-West University, Mafikeng, South Africa}
\address{K.M.: School of Mathematics, University of Birmingham,
  Edgbaston, Birmingham B15 2TT, U.K.}

\email{F.H.: himstedt@ma.tum.de}
\email{T.L.: lttung96@yahoo.com}
\email{K.M.: k.magaard@bham.ac.uk}


\begin{abstract}
Let $UY_n(q)$ be a Sylow $p$-subgroup of an untwisted Chevalley group
$Y_n(q)$ of rank $n$ defined over $\F_q$ where $q$ is a power of a
prime $p.$ We partition the set $\Irr(UY_n(q))$ of irreducible characters of 
$UY_n(q)$ into families indexed by antichains of positive roots of the root system
of type $Y_n$. 
We focus our attention on the families of characters of $UY_n(q)$ which are indexed
by antichains of length $1$.
Then for each positive root $\alpha$ we establish a one to one correspondence
between the minimal degree
members of the family indexed by $\alpha$ and the linear characters of a certain
subquotient $\overline{T}_\alpha$ of $UY_n(q)$.
For $Y_n = A_n$ our single root character construction recovers amongst other things
the
elementary supercharacters of these groups. Most importantly though this paper lays
the groundwork for our classification of the
elements of $\Irr(UE_i(q))$, $ 6 \leq i \leq 8$ and  $\Irr(UF_4(q))$. 
\end{abstract}

\maketitle


\section{Introduction}  \label{sec:intro}

Let $p$ be a prime, $q = p^a$ and $Y_n(q)$ be a finite quasisimple
group of untwisted rank $n$ defined over the field $\F_q$. By
$UY_n(q)$ we denote a Sylow $p$-subgroup of $Y_n(q)$ and by $\Irr(X)$ we denote 
the set of ordinary irreducible characters of the group $X.$

This paper lays the groundwork for our study of $\Irr(UE_i(q))$ where $ 6
\leq i \leq 8$, and $\Irr(UF_4(q))$. Our approach is to construct the 
characters explicitly using as primary parameters the underlying root system 
and the field. Our focus here is on the families of characters which we parameterize by 
a single root. For the classical groups, these families can be described recursively 
via character correspondences which can be achieved using Lemma \ref{nred}.
Establishing similar character correspondences for families parameterized by more than one root 
requires iterated applications of Lemma \ref{nred}. Given the length of the current paper
we treat the recursive method in a sequel. 

The solution of the dual problem,
the determination of the conjugacy classes $UY_n(q)$, has been achieved for rank up to 
$6$ in \cite{GR} and \cite{GMR} by Goodwin, Mosch, and R\"ohrle.
Combining these results with the results of this paper and its planned sequel 
opens the way for the construction of the generic character tables of these groups.  
For the groups $UD_4(q)$ this is presently being carried out by Goodwin, Le, and Magaard,
see \cite{GHLM}.

One motivation for constructing the generic character tables of 
$UY_n(q)$ is to aid in the construction of the cross characteristic
representations and in the determination of the decomposition numbers
of the exceptional groups of Lie type. Following Okuyama and Waki
\cite{OkuyamaWaki}, \cite{Waki_G2} and Himstedt, Huang and Noeske
\cite{H}, \cite{HH}, \cite{HN}, we see that characters of parabolic
subgroups are a useful tool in the computation of the decomposition
numbers of finite groups of Lie type in the cross characteristic case.

A second motivation is to explain exactly why the
primes $3$ and $5$ are bad for the exceptional groups of Lie type from
the point of view of the representation theory of $UY_n(q)$, where
$Y_n(q)$ is exceptional. A partial explanation is supplied in Le,
Magaard \cite{LeMag} where families of characters are exhibited whose
degree is not a power of $q= 3^a$ or $q = 5^a$. This generalizes a
construction for $UD_4(q)$ in \cite{HML_D4}, where it is shown that
there exists exactly one family of characters whose behavior at
the prime $2$ is different than for odd primes. At present we do not
know how many families of $UE_i(q)$ characters behave differently at
bad primes than at others. 

A third motivation is the conjecture of Higman from 1962 that the number of elements
of $\Irr(UA_r)$ is a polynomial in $\Z[q]$ and the generalization to other classical 
groups. This problem has led to the development of supercharacter theories.
By grouping conjugacy classes into so called superclasses and characters into supercharacters
these theories allow one to construct supercharacter tables which may be viewed 
as summarized versions of ordinary character tables. These theories were introduced by
Diaconis and Isaacs \cite{DI} for algebra groups such as $UA_n(q)$. Subsequently Andr\'e and Neto \cite{AndreNeto2}
developed supercharacter theories for $UB_n(q)$, $UC_n(q)$, and $UD_n(q)$. A common key feature in these theories 
is that the supercharacters are constructed as tensor products of elementary characters. 
An open problem in this area is the question of how to split non elementary supercharacters into ordinary irreducible characters.
In his thesis Le \cite{TungPhD} shows that the splitting of supercharacters into irreducibles
is governed by certain pattern subgroups of $UY_n$ and thus it would suffice to know 
that Higman's conjecture holds for pattern subgroups. This however is not true in general 
as was shown by Halasi in \cite{Hal}.

Our approach for classifying the irreducible characters of $UY_n(q)$ 
is based on an analysis of the supports of the centers of the characters
and character correspondences. 
We proceed as follows.
Let $\Phi$ be a root system of type $Y_n$ and let $\Phi^+$
denote the set of positive roots with respect to some choice of
simple roots. The group $UY_n$ is generated by the root subgroups 
$X_\alpha$ where $\alpha \in \Phi^+$. For $\chi \in \Irr(UY_n)$ we define a set
$\rs(\chi)$, see Definition \ref{def:rootker}, which consists of those elements
$\beta \in \Phi^+$ whose root subgroup $X_\beta$ lies in $Z(\chi)$ but not
in $\Ker(\chi)$. Subsets of
$\Phi^+$ which arise in this way are called \emph{representable}. 
We show that the number $N(\Phi)$ of representable sets in $\Phi^+$ is
a sum of (generalized) Narayana numbers and that $N(\Phi)$ is equal to
the number of antichains in the poset of positive roots of $\Phi$ as
well as the number of clusters in a cluster algebra of type $\Phi$,
see Proposition \ref{prop:numrepsets}. 

We say that a character $\chi$ is a \emph{single root character} if its
representable set is non-empty and as small as possible, that is $|\rs(\chi)| =1.$ Each
family of single root characters contains a collection of characters
of minimal degree which we call {\it midafi} characters and each
collection contains a special element that we call {\it standard
  midafi}. In case the root system is of type $A_n$ our standard
midafi characters are called basic characters by Andr\'e \cite{Andre2002} and are 
called elementary supercharacters by Diaconis and Isaacs \cite{DI}. 
For the root systems of types $B_n$, $C_n$ and $D_n$ our midafi characters 
differ from those defined in Andr\'e and Neto \cite{AndreNeto2}. To see this 
we note that all of our standard midafi characters are irreducible, whereas 
not all the elementary characters defined by Andr\'e and Neto are. 

For $\alpha \in \Phi^+$ we define 

\vspace{-0.3cm}

\[
\Irr(UY_n(q))_\alpha := \{ \chi \in \Irr(UY_n(q)) \mid \rs(\chi) = \{ \alpha \} \},
\]
the set of {\it single root characters lying over  $\alpha$} and the set

\vspace{-0.3cm}

\[ 
 \Irr^\mida(UY_i)_\alpha:= \{ \mu \in \Irr(UY_n(q))_\alpha \mid \mu \ \mbox{is midafi} \}.
\]

The observation that $X_\beta $ must act faithfully on any module affording $\chi$ 
for all $\beta \in \Phi^+$ such that $\alpha -\beta \in \Phi^+$  leads 
to the definition of the hook $h(\alpha)$ of $\alpha$, see \ref{def:hook}.
We show that the largest pattern subgroup contained in the kernel of $\chi \in \Irr(UY_n(q))_\alpha$ 
is the group generated by the set of root subgroups $X_\kappa$ with $\kappa \in k(\alpha)$.
We will see that $h(\alpha) \cap k(\alpha) = \emptyset$ and that in
general $\beta + \gamma \in h(\alpha) \cup k(\alpha)$ for all $\beta, \gamma \in h(\alpha)$. Also we will see that 
typically $H_\alpha:= \langle X_\gamma \ | \ \gamma \in h(\alpha) \rangle$
acts as special group on any module affording $\chi$. As a result we obtain that 
$\chi(1) = cq^d$ with $d = (|h(\alpha)|-1)/2$. When $c=1$, then $\chi$ is a midafi.

More generally, using our Reduction Lemma \ref{nred}, we can interpret $c$ as the degree of 
an irreducible character of a suitable quotient $\overline{T}_\alpha$ of a certain subgroup
$S_\alpha < UY_n(q)$. 
Our main theorems can now be stated under mild hypotheses on the prime
$p$, see Hypothesis \ref{hyp:p}. 

\begin{thm} \label{thm:main_classical}
Let $\Phi_n$ be an irreducible root system of type $A_n$, $B_n$, $n \geq 2$, $C_n$, $n \geq 3$
or $D_n$, $n \geq 4$ and $\F_q$ a finite field of characteristic $p$ such that Hypothesis \ref{hyp:p} holds. 
For every positive root $\alpha \in \Phi_n^+$ the map
\[
\Psi: \Irr(\overline{T}_\alpha) \times \Irr(X_\alpha)^* \to
\Irr(UY_n)_\alpha, \quad (\mu, \lambda) \mapsto 
 (\infl_{\overline{T}_\alpha}^{S_\alpha} \mu \cdot
   \infl_{X_\alpha}^{S_\alpha} \lambda)^{UY_n}
\]
is a one to one correspondence. 
\end{thm}

\begin{thm} \label{thm:main_exceptional}
Let $i \in \{2,4,6,7,8\}$ and let $\Phi_i$ be a root system of type 
$G_2$, $F_4$ or $E_i$ respectively and $\F_q$ a finite field of characteristic $p$ such that Hypothesis \ref{hyp:p} holds.
For every positive root $\alpha \in \Phi_i^+$ the map
\[
\Psi: \Irr^\lin(\overline{T}_\alpha) \times \Irr(X_\alpha)^* \to
\Irr^\mida(UY_i)_\alpha, \quad (\mu, \lambda) \mapsto 
 (\infl_{\overline{T}_\alpha}^{S_\alpha} \mu \cdot
   \infl_{X_\alpha}^{S_\alpha} \lambda)^{UY_i}
\]
is a one to one correspondence.
\end{thm}

In case $\Phi$ is classical, we determine for every $\alpha \in \Phi^+$ the structure 
of $\overline{T}_\alpha$, whereas for $\Phi$ exceptional we determine for every 
$\alpha \in \Phi_i^+$ the number of midafis for $\alpha$ and their degrees.

To explain why it is that our results differ for classical and exceptional groups we need to 
address the issue of how $\overline{T}_\alpha$ is constructed. We partition the set 
$h(\alpha) \setminus \{ \alpha \}$ into two subsets of equal size $a(\alpha)$ and $\ell(\alpha)$ 
which we call the arm, respectively the leg of $h(\alpha)$. Then we define the source 
$s(\alpha) := \Phi^+ \setminus a(\alpha)$. A priory we have $2^{|a(\alpha)|}$ 
choices for the sets $a(\alpha)$ and $\ell(\alpha)$. However, if our choice for the arm and 
leg satisfy the following conditions 
\begin{enumerate}
\item $s(\alpha)$ is closed under addition of roots, and that
\item $\lambda + \sigma \in \ell(\alpha) \cup k(\alpha)$ for all 
$\lambda \in \ell(\alpha)$ and all $\sigma \in s(\alpha)$ such that $\lambda + \sigma \in \Phi^+$, 
\end{enumerate}

then we can achieve the hypotheses of Lemma \ref{nred} to establish the 
correspondence in Theorem \ref{thm:main_classical}. In this case the group 
$S_\alpha := \langle X_\sigma \ | \ \sigma \in s(\alpha) \rangle $ contains a normal 
subgroup  $K_\alpha:= \langle X_\lambda \ | \ \lambda \in k(\alpha) \rangle$ such that 
$\overline{S}_\alpha$, the image of $S_\alpha$ in $UY_n(q)/K_\alpha$, contains a normal 
subgroup $\overline{L}_\alpha:= \langle \overline{X}_\lambda \ | \ \lambda \in \ell(\alpha) \rangle$ 
such that $\overline{S}_\alpha/\overline{L}_\alpha = \overline{T}_\alpha \times \overline{X}_\alpha$.
 
Conditions (1) and (2) can always be achieved when $\Phi$ is classical or of type $G_2$. 
When $\Phi$ is exceptional condition (1) can always be achieved. However when 
$\Phi$ is of type $E_8$, then condition (2)
can not be achieved for $46$ of the $120$ roots of $\Phi^+$. 
The numbers for types $E_7$, $E_6$ and $F_4$ are  $11$ out of $63$, $2$ out of $36$, and $2$ out $24$, respectively. 
Nevertheless we have 

\begin{thm} \label{thm:main_exceptionalnormal}
Let $i \in \{2,4,6,7,8\}$ and let $\Phi_i$ be a root system of type 
$G_2$, $F_4$ or $E_i$ respectively and $\F_q$ a finite field of characteristic $p$ such that Hypothesis \ref{hyp:p} holds.
For every positive root $\alpha \in \Phi_i^+$ for which conditions (1) and (2) above can be achieved the map
\[
\Psi: \Irr(\overline{T}_\alpha) \times \Irr(X_\alpha)^* \to
\Irr(UY_i)_\alpha, \quad (\mu, \lambda) \mapsto 
 (\infl_{\overline{T}_\alpha}^{S_\alpha} \mu \cdot
   \infl_{X_\alpha}^{S_\alpha} \lambda)^{UY_i}
\]
is a one to one correspondence.
\end{thm}

We remark that if $p>3$ then  Hypothesis \ref{hyp:p} is satisfied for
the groups $UG_2(q)$ and in this case every irreducible character of $UG_2(q)$
has degree $1$, $q$ or $q^2$ and all irreducible characters of degree $>1$ are midafis.
We remark further that the number of possible choices for $a(\alpha)$ so that (1) and (2) above are satisfied 
is $(|h(\alpha)|-1)/2$ in type $A_n$ and much smaller in all other cases.

In case no choice of $a(\alpha)$ achieves condition (2) we pick from those choices 
which satisfy condition (1) the one that minimizes the index of $\overline{L}_\alpha$ in its normal 
closure in $\overline{S}_\alpha$. This amounts to minimizing the size of 

\vspace{-0.2cm}

$$\overline{\ell}(\alpha) := \{\tau \in s(\alpha) \ | \ \mbox{There exist} \  \lambda \in \ell(\alpha) \ \mbox{and} \
 \sigma_i \in s(\alpha) \ \mbox{such that} \ \tau = \lambda + \sum_i \sigma_i \}.$$

\vspace{-0.2cm}

Let $\overline{\tilde{L}}_\alpha := \langle \overline{X}_\mu \ | \  \mu \in \overline{\ell}(\alpha)  \rangle $.
Then $\overline{\tilde{L}}_\alpha$ is normal in $\overline{S}_\alpha$ and finally we can define 
the group $\overline{T}_\alpha$ in the statement of Theorem \ref{thm:main_exceptional} via 
$\overline{S}_\alpha / \overline{\tilde{L}}_\alpha \cong \overline{T}_\alpha \times \overline{X}_\alpha$.

With the machinery set up in this paper we are able to give full
descriptions of the character correspondences for the exceptional
groups. It should be noted that if $X_\alpha$ projects faithfully into a 
classical quotient of $UY_n(q)$, then condition (2) can always be achieved. Of the 
$120$ positive roots of $E_8$ exactly $63$ have the property that 
$X_\alpha$ projects faithfully into a classical quotient. In $46$ of the $57$ 
remaining cases condition (2) can not be achieved. (For $F_4$, $E_6$ and $E_7$
 the numbers are $2$ out of $10$, $2$ out of $7$, and $11$ out of $23$, respectively.) 
Compounding this is the fact that in those cases where condition (2) can not be achieved, 
the descriptions of the single root characters involve up to possibly $16$  (generally 
as many as there are subhooks listed in Table \ref{tab:armse8}) recursive 
applications of our Reduction Lemma \ref{nred} and are thus beyond the scope of this article. 

A full treatment of the single root characters  will appear in a 
forthcoming article. The case $\Phi$ of type $F_4$ will be considered in a forthcoming article by 
Goodwin, Le and Paolini \cite{GLP}, where we see that our machinery also generalizes to the case of multiple 
root characters which we encounter in root systems of exceptional Lie type.

It is worth remarking that the structure of the group $\overline{T}_\alpha$ is of the form $V \rtimes K$ where $K$ is classical 
and $V$ is a non-faithful $K$-module. This observation, which Thompson had already made for the unitriangluar 
groups, begins to reveal why inductive approaches to the solution of Higman's conjecture have eluded us so far. 

The paper is organized as follows. In Section \ref{sec:nota} we fix
notation and prove Lemma \ref{nred} which is fundamental to our construction of characters.
In Section \ref{sec:pattern} we define closed patterns
(additively closed subsets of roots) and the corresponding pattern
subgroups, and establish some of their basic properties. This is used
in Sections \ref{sec:rkeretc} and \ref{sec:representable} to define
the key terms of this paper, such as hooks, root kernels, representable
sets and to establish their basic properties. We show that the number of representable 
sets is equal to the number of antichains in the poset of positive
roots. In Section \ref{sec:single} we study single root characters to
lay the foundations for establishing our correspondences. For the
classical root systems our main theorem is established in Section
\ref{sec:single_classical} and for the exceptional root systems in
Section \ref{sec:exc_midafis}.


\section{Notation and a reduction lemma}
\label{sec:nota}
In this preliminary section we fix our notation and prove the key
lemma which is needed to establish the character correspondences in
Sections \ref{sec:single_classical} and~\ref{sec:exc_midafis}.  

\subsection{Character theoretic setting}
\label{charset}

For any finite group $U$ let $\Irr(U)$ be the set of complex
irreducible characters~of~$U$ and $\Irr^\lin(U) := \{\chi \in \Irr(U)
\mid \chi(1) = 1\}$ the set of linear characters. Let $(\cdot,\cdot)_U$  
or $(\cdot,\cdot)$ be the usual scalar product on the space of
$\C$-valued class functions of $U$. We write $\mathbf{1}_U$ or
$\mathbf{1}$ for the trivial character of $U$ and set 
$\Irr(U)^* := \Irr(U) \setminus \{\mathbf{1}_U\}$. Suppose that 
$N$ is a normal subgroup of $U$ and that $H$ is a (not necessarily
normal) subgroup of $U$. If $\chi$ is a character of $U$ and $\lambda$
is a character of $H$ and $\psi$ is a character of the factor group
$U/N$, we write $\lambda^U$ for the character of $U$ induced by
$\lambda$ and $\chi|_H$ for the restriction of $\chi$ to $H$ and
$\infl_{U/N}^U \psi$ for the inflation of $\psi$ to~$U$. For 
$\lambda \in \Irr(H)$ we set 
\begin{eqnarray*}
\Irr(U, \lambda) & := & \{\chi \in \Irr(U) \, | \, (\chi,\lambda^U)>0\}
\ \text{and}\\
\Irr(U/N,\lambda) & := & \{ \chi \in \Irr(U , \lambda) \, | \, N
\subseteq \Ker(\chi)\},
\end{eqnarray*}
where $\Ker(\chi)$ denotes the kernel of $\chi$. The following lemma
provides a character correspondence between finite groups and certain
subgroups. 

\begin{lem}
\label{nred}
Let $H$ be a subgroup of a finite group $U$ and $X$ a set of 
representatives for $U/H$. Furthermore, let $Y$, $Z$ be subgroups of
$H$ and $\lambda \in \Irr(Z)$ such that
\begin{enumerate}
\item[(a)] $Z \subseteq Z(U)$,

\item[(b)] $Y \unlhd H$,

\item[(c)] $Z \cap Y = \{1\}$,

\item[(d)] $ZY \unlhd U$,

\item[(e)] for the extension 
$\tilde{\lambda} \in \Irr(ZY) = \Irr(Z \times Y)$ of $\lambda$
with $Y \subseteq \Ker(\tilde{\lambda})$ we have
 ${}^x\tilde{\lambda}\neq \tilde{\lambda}$ for
  all $x \in X \setminus H$.
\end{enumerate}
Then the map
$\Phi: \Irr(H/Y, \lambda) \rightarrow \Irr(U, \lambda) \cap \Irr(U,
\mathbf{1}_Y), \, \chi \mapsto \chi^U$ 
is bijective.
If additionally
\begin{enumerate}
\item[(f)] $|X|=|Y|$
\end{enumerate}
holds, then $\Irr(U, \lambda) \cap \Irr(U, \mathbf{1}_Y) = \Irr(U, \lambda)$.
\end{lem}

\pr
Suppose that (a)-(e) are true. By (c),(d),(e) the character
$\tilde{\lambda}$ is an irreducible character of the normal subgroup 
$ZY$ of $U$ so that we can apply Clifford theory. For all elements 
$h \in H$, $y \in Y$, $z \in Z$ we have 
${}^h\tilde{\lambda}(zy) = \tilde{\lambda}(z^h y^h) =
\tilde{\lambda}(zy)$ by (a) and (b). So $H$ is contained in the
inertia subgroup $I_U(\tilde{\lambda})$, and from (e) we get 
$H = I_U(\tilde{\lambda})$. By Clifford theory
\cite[Theorem~(6.11)]{Isaacs:76} the map 
\[
\Phi: \Irr(H, \tilde{\lambda}) \rightarrow \Irr(U,
\tilde{\lambda}), \chi \mapsto \chi^U
\]
is a bijection. Since 
$\Irr(H, \tilde{\lambda}) = \Irr(H/Y, \lambda)$ and
$\Irr(U, \tilde{\lambda}) = \Irr(U, \lambda) \cap \Irr(U,
\mathbf{1}_Y)$ the first claim follows.

Assume additionally that (f) holds. By (a), (e) and (f) the
character~$\lambda$ has at least $|X|=|Y|$ distinct extensions to 
$ZY = Z\times Y$. It follows that $\lambda$ has exactly $|Y|$ distinct
extensions to $ZY = Z\times Y$ and that these are permuted transitively 
by the conjugation action of $U$. One of these extensions is
$\tilde{\lambda}$. Thus each $\mu \in \Irr(ZY)$ with 
$\mu|_Z = \lambda$ is conjugate in $U$ to $\tilde{\lambda}$. Hence 
\[
\Irr(U, \lambda) \subseteq \Irr(U, \tilde{\lambda}) =
\Irr(U, \lambda) \cap \Irr(U, \mathbf{1}_Y) \subseteq \Irr(U, \lambda).
\]
\epr

\subsection{Lie theoretic setting}
\label{lie}

We fix a power $q=p^a$ of a prime $p$ and write $\F_q$ for a field
with $q$ elements. Let $Y_n(q)$ be an untwisted Chevalley group
defined over~$\F_q$, constructed from a simple Lie algebra with the
irreducible root system $\Phi$ of Dynkin type $Y$ and rank $n$ as
described in \cite[Section~4.4]{Carter1}. So $Y_n(q)$ is generated 
by elements~$x_\alpha(t)$ for $\alpha \in \Phi$ and $t \in \F_q$. Let 
$X_\alpha := \la x_\alpha(t) \mid t \in \F_q \ra$ be the root
subgroup corresponding to a root $\alpha \in \Phi$.

We fix a set $\Delta = \{\delta_1, \dots, \delta_n\}$ of simple roots 
and write $\Phi^+$ for the corresponding set of positive roots. So
each $\alpha \in \Phi^+$ can be written as 
$\alpha = \sum_{i=1}^n m_i \delta_i$ where the coefficients $m_i \ge 0$ 
are integers. We write $\height(\alpha) := \sum_{i=1}^n m_i$ for the
height of $\alpha$. 

Let $UY_n(q)$ or $UY_n$ be the subgroup generated by 
$\{ x_\alpha(t) \, | \, \alpha \in \Phi^+, t\in\F_q\}$. So $UY_n$ is a
maximal unipotent subgroup and a Sylow $p$-subgroup of $Y_n(q)$. For
example, it is well-known that $UA_5(q)$ is isomorphic to the subgroup
of $\SL_6(q)$ consisting of all upper unitriangular matrices.
Let $\le$ be a total ordering on $\Phi^+$. Then each element
$u \in UY_n$ can be written uniquely as
\begin{equation} \label{eq:elemu}
u = \prod_{\alpha \in \Phi^+} x_\alpha(t_\alpha),
\end{equation}
where the product is taken over all positive roots in increasing
order. The multiplication of the elements of $UY_n$ is determined
by commutator relations after fixing the signs of certain structure
constants corresponding to the so-called extraspecial pairs of roots;
see \cite[Sections~4.2 and 5.2]{Carter1} for details.

We say that a non-empty subset $\Psi \subseteq \Phi$ is a root
subsystem if $\sigma_\alpha(\Psi)=\Psi$ for all reflections
$\sigma_\alpha$ corresponding to roots $\alpha \in \Psi$. Let 
$\Psi$ be a root subsystem of $\Phi$ of Dynkin type $Y'$
and rank $n'$ and let $\Z\Psi$ be the $\Z$-span of $\Psi$. We 
define $\Psi^+ := \Psi \cap \Phi^+$ and 
$U_\Psi := \prod_{\alpha \in \Psi^+} X_\alpha$ where the product
is taken over all $\alpha \in \Psi^+$ in increasing order. If
$\Z\Psi \cap \Phi = \Psi$ then the commutator relations and the 
properties of the structure constants in~\cite[p.58-59]{Carter1} imply 
that $U_\Psi$ is a subgroup of $UY_n$ isomorphic to $UY'_{n'}$. 

Recall that we can define a partial order $\preceq$ on $\Phi^+$
as follows (see~\cite[10.1]{Humphreys}): For roots $\alpha,\beta \in \Phi^+$
we write $\alpha \prec \beta$ if $\beta-\alpha$ is a non-zero sum 
of positive roots and we write $\alpha \preceq \beta$ 
if $\alpha \prec \beta$ or $\alpha=\beta$. The following lemma is a
special case of \cite[Lemma~3.2]{Sommers}. 

\begin{lem} \label{la:maxchain}
For all $\alpha, \beta \in \Phi^+$ with $\alpha \prec \beta$
there are roots $\gamma_i \in \Phi^+$ such that $\alpha=\gamma_0 \prec
\gamma_1 \prec \gamma_2 \prec \dots \prec \gamma_{s-1} \prec
\gamma_s=\beta$ and $\gamma_i-\gamma_{i-1} \in \Delta$ for all $i=1,2,\dots,s$.
\end{lem}
\pr
The lemma follows from \cite[Lemma~3.2]{Sommers}.
\epr

\medskip

A chain in $\Phi^+$ is a subset $C = \{\gamma_1,\gamma_2, \dots ,
\gamma_s\} \subseteq \Phi^+$ such that $\gamma_i \preceq \gamma_{i+1}$ 
for all $1 \leq i < s-1$. We say that $C$ is \emph{unrefinable} if
there is no root $\gamma \in \Phi^+ \setminus C$ such that 
$\gamma_i \prec \gamma \prec \gamma_{i+1}$ for some $i$. By
Lemma~\ref{la:maxchain} unrefinable chains have the property 
that the difference of consecutive elements is a simple root. By
\cite[Lemma~10.4A]{Humphreys} the highest positive root of $\Phi$ is
the unique maximal element of the poset $(\Phi^+, \preceq)$.


\section{Pattern subgroups}
\label{sec:pattern}

In \cite[Section~2]{IsaacsTriang}, Isaacs defines pattern
subgroups of the subgroup $U_m(q)$ of $\GL_m(\F_q)$ consisting of the
upper unitriangular matrices. Since $U_m(q)$ is isomorphic to
$UA_{m-1}(q)$ these pattern subgroups can be identified with
subgroups of $UA_{m-1}(q)$ in a natural way. The following
definition generalizes the notion of pattern subgroups to other Dynkin 
types. 
We assume the setting described in Section~\ref{sec:nota}.
In particular, $q$ is a power of a prime $p$ and $Y_n(q)$ is an
untwisted Chevalley group defined over~$\F_q$ with the irreducible
root system~$\Phi$. The set of positive roots is denoted by~$\Phi^+$ 
and $UY_n(q)$ or $UY_n$ is the subgroup of $Y_n(q)$ generated by the
root subgroups $X_\alpha$ for $\alpha \in \Phi^+$. 

\begin{defn} \label{def:pattern}
Let $S$ be a subset of $\Phi^+$.
\begin{enumerate}
\item[(a)] The set $S$ is called a \emph{closed pattern} if for all
  roots $\alpha, \beta \in S$ we have $\alpha+\beta \in S$ or
  $\alpha+\beta \not\in \Phi^+$. 

\item[(b)] For a closed pattern $S$ let $P(S)$ be the subgroup of $UY_n$ 
  generated by the root subgroups $X_\alpha$ for $\alpha \in S$. We
  call $P(S)$ the \emph{pattern subgroup} corresponding to $S$. 
\end{enumerate}
\end{defn}

Note that when we have a total ordering $\le$ on $\Phi^+$ and a
closed pattern $S$ then~(\ref{eq:elemu}), \cite[Lemma~3.6.3]{Carter1}
and the commutator relations imply that 
$P(S) = \prod_{\alpha \in S} X_\alpha$ where the product is taken
over the roots in $S$ in increasing order. Obviously, each root
subgroup $X_\alpha$ for $\alpha \in \Phi^+$ and the trivial subgroups 
$UY_n = P(\Phi^+)$ and $\{1\} = P(\emptyset)$ are pattern subgroups.

\begin{defn} \label{def:normal}
Let $M, N$ be subsets of a closed pattern $S \subseteq \Phi^+$.
\begin{enumerate}
\item[(a)] We say that $M$ is a \emph{closed subpattern} of $S$ if $M$
  is a closed pattern.

\item[(b)] The closed pattern \emph{generated} by $M$ is the
  intersection of all closed patterns containing~$M$.

\item[(c)] We say that $M$ \emph{normalizes} $N$ if for
  all $\alpha \in M$, $\beta \in N$ we have $\alpha + \beta \in N$ 
  or $\alpha + \beta \not \in \Phi^+$. We say that $N$ is
  \emph{normal} in $S$ if $S$ normalizes $N$. In this case we write $N
  \unlhd S$ and call $P(S)/P(N)$ the \emph{quotient pattern group} of
  $P(S)$ corresponding to~$N$. We say that $N$ is \emph{normal} if $N$
  is normal in $\Phi^+$.

\item[(d)] The \emph{normal closure} of $M$ in $S$ is the intersection
  of all normal closed subpatterns of $S$ containing~$M$. 
\end{enumerate}
\end{defn}

Obviously, the closed pattern generated by $M \subseteq \Phi^+$ is the
smallest closed pattern containing $M$ and the normal closure of $N$
in $S$ is the smallest normal closed subpattern of $S$ containing $M$.
Note that the normal closed patterns are exactly the (upper) order
ideals of the root poset $(\Phi^+, \preceq)$ in the sense of~\cite{Sommers}.

\begin{rem} \label{rem:normpat}
If $N \unlhd S$ then~(\ref{eq:elemu}), \cite[Lemma~3.6.3]{Carter1} and
the commutator relations imply that~$N$ is a closed subpattern of $S$
and $P(N) \unlhd P(S)$. In particular, $P(S)/P(N)$ is a well-defined
factor group. Let $\pi: P(S) \rightarrow P(S)/P(N)$ be the canonical
projection and $\gamma \in S \setminus N$. Since $\pi$ maps the root
subgroup $X_\gamma$ injectively into $P(S)/P(N)$ we often identify
$X_\gamma$ with $\pi(X_\gamma)$.
\end{rem}

To avoid degeneracies in the commutator relations we will often assume
the following hypothesis.

\begin{hypo} \label{hyp:p}
If the Dynkin diagram of $\Phi$ has a double or triple edge assume
that
\begin{itemize}
\item $p>2$ if $\Phi$ is of type $B_m$ or $C_m$ $(m \ge 2)$ or $F_4$,
\item $p>3$ if $\Phi$ is of type type $G_2$.
\end{itemize}
\end{hypo}

Next, we consider the connection between normal closed patterns 
and normal subgroups. For $x,y \in UY_n$ we set $[x,y] := x^{-1}y^{-1}xy$. 
Part (c) of the following lemma is stated in~\cite{BT} for connected
reductive groups; see also~\cite[1.12, 1.13]{DM}. Since we need it for
finite groups we sketch a proof.

\begin{lem} \label{la:patnorm}
Let $\alpha, \beta \in \Phi^+$ such that $\alpha+\beta \in \Phi^+$. We
assume that Hypothesis~\ref{hyp:p} is satisfied and set
$\Phi_{\alpha,\beta}^{>0} := \{i\alpha+j\beta \in \Phi^+ \mid i,j \in \Z_{>0} \}$.
\begin{enumerate}
\item[(a)]  For all $s,t \in \F_q$ there are constants 
  $c_{ij\alpha\beta} \in \F_q$ such that $c_{11\alpha\beta} \neq 0$ and 
  \begin{equation} \label{eq:rootprod}
  x_\alpha(s)^{-1}x_\beta(t)^{-1}x_\alpha(s)x_\beta(t) =
  \prod_{i,j > 0} x_{i\alpha+j\beta}(c_{ij\alpha\beta} \cdot (-t)^is^j),
  \end{equation}
  where the product is taken over all $i,j \in \Z_{>0}$, 
  such that $i\alpha+j\beta \in \Phi^+$ and the terms in the product
  are ordered from left to right so that $i+j$ is increasing. 

\smallskip

\item[(b)] The set $\Phi_{\alpha,\beta}^{>0}$ is a closed pattern
  and a subset of the normal closure of $\{\alpha+\beta\}$ in
  the closed pattern generated by $\alpha$ and $\beta$. 

\smallskip

\item[(c)] (\cite[Remarque 2.5]{BT}) \, $\displaystyle [X_\alpha,
  X_\beta] = \prod_{\gamma \in \Phi_{\alpha,\beta}^{>0}} \hspace{-0.1cm} X_\gamma$.
\end{enumerate}
\end{lem}

Note that because $\Phi_{\alpha,\beta}^{>0}$ is a closed pattern the
product in (c) does not depend on the order of the factors.

\medskip

\pr 
(a) By assumption, $\Phi$ is irreducible.
Let $\Phi_{\alpha,\beta} := \Phi \cap (\Z \alpha + \Z \beta)$. 
Since $\alpha+\beta \in \Phi$ we know from \cite[Lemma~3.6.3]{Carter1}
that $\Phi_{\alpha,\beta}$ is a root system of type $A_2$, $B_2$ or $G_2$. 
Let $-m\alpha+\beta, \dots, \beta, \dots, m'\alpha+\beta$
be the $\alpha$-string through $\beta$. Then (\ref{eq:rootprod}) holds
with $c_{11\alpha\beta}=\pm(m+1)$ by 
\cite[Corollary~5.2.3 and Sections~4.1, 4.2]{Carter1}. We have to
show that $c_{11\alpha\beta} \neq 0$. 

Suppose that $\Phi$ has type $A_n$, $D_n$, $E_6$, $E_7$ or $E_8$. 
Since all roots of $\Phi$ have the same length, $\Phi_{\alpha,\beta}$
has type $A_2$ and hence $m=0$; see \cite[9.3]{Humphreys}. Thus
$c_{11\alpha\beta} = \pm1$ and $c_{11\alpha\beta} \neq 0$. Suppose
that $\Phi$ has type $B_n$, $C_n$ or $F_4$. Again the lengths 
of the roots imply that $\Phi_{\alpha,\beta}$ has type $A_2$ or $B_2$
and so $m \in \{0,1\}$ by \cite[9.3]{Humphreys}. Thus  
$c_{11\alpha\beta} = \pm 1$ or $c_{11\alpha\beta} = \pm2$ and so 
$c_{11\alpha\beta} \neq 0$ by Hypothesis~\ref{hyp:p}. Finally suppose
that $\Phi$ has type $G_2$. Then $\Phi_{\alpha,\beta}$ has type $A_2$
or $G_2$ and so $m \in \{0,1,2\}$ by \cite[9.3]{Humphreys}. Thus
$c_{11\alpha\beta} = \pm 1, \pm 2$ or $\pm 3$. Thus 
$c_{11\alpha\beta} \neq 0$ because $p>3$.

\medskip

(b), (c) Obviously, $\Phi_{\alpha,\beta}^{>0}$ is a closed
pattern and thus the product $\prod_{\gamma \in \Phi_{\alpha,\beta}^{>0}}
X_\gamma$ is a subgroup of $UY_n$. The commutator relations imply that
$[X_\alpha, X_\beta] \subseteq \prod_{\gamma \in \Phi_{\alpha,\beta}^{>0}}
X_\gamma$. Since $\Phi_{\alpha,\beta}$ has type $A_2$, $B_2$ or $G_2$
we see that $\Phi_{\alpha,\beta}^{>0}$ is one of the following sets:
$\{\alpha+\beta\}$, $\{\alpha+\beta, 2\alpha+\beta\}$,
$\{\alpha+\beta, \alpha+2\beta\}$,
$\{\alpha+\beta, 2\alpha+\beta, \alpha+2\beta\}$,
$\{\alpha+\beta, 2\alpha+\beta, 3\alpha+\beta, 3\alpha+2\beta\}$,
$\{\alpha+\beta, \alpha+2\beta, \alpha+3\beta, 2\alpha+3\beta\}$.

We only treat the most complicated case $\Phi_{\alpha,\beta}^{>0} =
\{\alpha+\beta, 2\alpha+\beta, 3\alpha+\beta, 3\alpha+2\beta\}$. 
Because $2\alpha+\beta = \alpha + (\alpha+\beta)$, $3\alpha+\beta =
\alpha + (2\alpha+\beta)$ and $3\alpha+2\beta = \beta+(3\alpha+\beta)$
the statement in (b) follows. 

To prove (c) note that $[X_\alpha, X_\beta] \unlhd \langle X_\alpha, X_\beta\rangle$; 
see~\cite[Hilfssatz III.1.6 (b)]{HuppertI}.
By (a) we have $[X_\alpha, X_\beta] = X_{\alpha+\beta}$ mod
$X_{2\alpha+\beta} X_{3\alpha+\beta} X_{3\alpha+2\beta}$. Thus,
there are $d, d', d'' \in \F_q$ such that 
$u := x_{\alpha+\beta}(1) x_{2\alpha+\beta}(d) x_{3\alpha+\beta}(d')
x_{3\alpha+2\beta}(d'') \in [X_\alpha, X_\beta]$. Again it follows
from (a) that $[X_\alpha, u] = X_{2\alpha+\beta}$ mod  
$X_{3\alpha+\beta} X_{3\alpha+2\beta}$ so there are 
$f, f' \in \F_q$ such that $u' := x_{2\alpha+\beta}(1)
x_{3\alpha+\beta}(f) x_{3\alpha+2\beta}(f') \in [X_\alpha, X_\beta]$.
Again from (a) we get that $[X_\alpha, u'] = X_{3\alpha+\beta}$ mod
$X_{3\alpha+2\beta}$ and then
$u'' := x_{3\alpha+\beta}(1) x_{3\alpha+2\beta}(g) \in [X_\alpha, X_\beta]$ 
for some $g \in \F_q$. From (a) we get $[X_\beta, u''] =
X_{3\alpha+2\beta} \subseteq [X_\alpha, X_\beta]$. Now we can work
backwards and get $X_{3\alpha+\beta} \subseteq [X_\alpha, X_\beta]$ 
and then $X_{2\alpha+\beta} \subseteq [X_\alpha, X_\beta]$ and finally
$X_{\alpha+\beta} \subseteq [X_\alpha, X_\beta]$ proving (c). 
\epr

\begin{cor} \label{cor:patnorm}
Let $N \subseteq S \subseteq \Phi^+$ be closed patterns and suppose
that Hypothesis~\ref{hyp:p} holds. Then $N \unlhd S$ if and only if
$P(N) \unlhd P(S)$.
\end{cor}

\pr
By the remarks after Definition~\ref{def:normal} we already know that
$N \unlhd S$ implies $P(N) \unlhd P(S)$ even without the condition on
$p$. Now suppose that $P(N) \unlhd P(S)$ and let $\alpha \in S$, 
$\beta \in N$ such that $\alpha+\beta \in \Phi^+$. By
Lemma~\ref{la:patnorm}~(c) we have 
$X_{\alpha+\beta} \subseteq [X_\alpha, X_\beta] \subseteq [X_\alpha, P(N)] 
\subseteq P(N)$. The remark after Definition~\ref{def:pattern} and
the uniqueness in (\ref{eq:elemu}) imply that $\alpha+\beta \in N$. 
\epr

\medskip

Without Hypothesis~\ref{hyp:p}, the converse in
Corollary~\ref{cor:patnorm} is not true in general: Suppose that
$\Phi$ has type $B_2$ and that $\{ \alpha, \beta \}$ is a set of
simple roots. If $p=2$ then $X_{\alpha+\beta} \unlhd UB_2(q)$ but the
closed pattern $\{ \alpha + \beta \}$ is not normal. 

The following lemmas are of theoretical and computational use. They
show that derived subgroups and centers of quotient pattern groups are
compatible with the root structure.

\begin{lem} \label{la:derquopat}
Let $N \subseteq S \subseteq \Phi^+$ be closed patterns such that 
$N \unlhd S$ and assume that Hypothesis~\ref{hyp:p} holds. Then
$D := (\{\alpha+\beta \mid \alpha, \beta \in S\} \cap \Phi^+) \cup N$
is a normal closed subpattern of $S$ and 
$[P(S)/P(N), P(S)/P(N)] = P(D)/P(N)$. In particular, the derived
subgroup of a pattern subgroup is also a pattern subgroup.
\end{lem}

\pr
Since $S$ is a closed pattern and $N \subseteq S$ we have 
$N \subseteq D \subseteq S$. For all $\alpha \in S$, 
$\beta \in D \subseteq S$ with $\alpha+\beta \in \Phi^+$ we have
$\alpha+\beta \in \{\alpha+\beta \mid \alpha, \beta \in S\} \cap
\Phi^+ \subseteq D$. Hence~$D$ is a normal closed subpattern of $S$
and $P(N) \subseteq P(D) \unlhd P(S)$. It~follows from
Lemma~\ref{la:patnorm} (c) that 
$P(D)/P(N) \subseteq [P(S)/P(N), P(S)/P(N)]$. The commutator relations
and Lemma~\ref{la:patnorm} (b) show that $(P(S)/P(N)) / (P(D)/P(N))$ 
is abelian. Thus $P(D)/P(N)$ is the derived subgroup of $P(S)/P(N)$.
\epr

\begin{lem} \label{la:cenquopat}
Let $N \subseteq S \subseteq \Phi^+$ be closed patterns such that 
$N \unlhd S$ and assume that Hypothesis~\ref{hyp:p} holds. Then
\[
Z := \{\alpha \in S \mid \text{for all} \ \gamma \in S:
\alpha+\gamma \not\in \Phi^+ \ \text{or} \ \alpha+\gamma \in N\} \cup N
\]
is a normal closed subpattern of $S$ and $Z(P(S)/P(N)) = P(Z)/P(N)$. 
In particular, the center $Z(P(S)/P(N))$ is isomorphic to a direct product of
root subgroups.
\end{lem}

\pr
Suppose that $\alpha \in S$ and $\beta \in Z$ such that 
$\alpha + \beta \in \Phi^+$. By the definition of~$Z$ and the
normality of $N$ in $S$ we have 
$\alpha + \beta \in N \subseteq Z$ and hence $Z$ is a normal closed
subpattern of $S$ containing $N$. The commutator relations and
Lemma~\ref{la:patnorm}~(b) imply that $P(Z)/P(N) \subseteq Z(P(S)/P(N))$.

Now suppose that $u \, P(N) \in Z(P(S)/P(N)) \setminus \{1\}$.
According to (\ref{eq:elemu}) we write 
\[
u = x_{\alpha_1}(t_1) \cdots x_{\alpha_s}(t_s)
x_{\alpha_{s+1}}(t_{s+1}) \cdots x_{\alpha_{s'}}(t_{s'})
\]
where $\alpha_i \in S \setminus N$ and 
$t_i \neq 0$ for $i = 1, \dots, s'$ and  
$m := \height(\alpha_1) = \dots = \height(\alpha_s) <
\height(\alpha_{s+1}) \le \height(\alpha_{s+2}) \le \dots \le
\height(\alpha_{s'})$. We show that $\alpha_i \in Z$ for 
$i=1, \dots, s'$ by downwards induction on $m$. 

If $m = \max\{\height(\alpha) \mid \alpha \in \Phi^+\}$ 
then $s=s'=1$ and $\alpha_1+\gamma \not\in \Phi^+$ for all 
$\gamma \in S$ and thus $\alpha_1 \in Z$. Assume that 
$m < \max\{\height(\alpha) \mid \alpha \in \Phi^+\}$. Let 
$\gamma \in S$ and $M$ be the normal subgroup of $P(S)$ generated 
by all $X_\alpha$ to roots $\alpha \in S$ such that 
$\alpha \in N$ or $\height(\alpha) > m + \height(\gamma)$. 
By Lemma~\ref{la:patnorm} (a) there exist 
$\tilde{t}_1, \dots, \tilde{t}_s \in \F_q \setminus \{0\}$,
$\tilde{u} \in M$ such that
\[
x_\gamma(1)^{-1} u x_\gamma(1) = x_{\alpha_1}(t_1)
x_{\alpha_1+\gamma}(\tilde{t}_1) \cdots x_{\alpha_s}(t_s)
x_{\alpha_s+\gamma}(\tilde{t}_s) x_{\alpha_{s+1}}(t_{s+1}) \cdots
x_{\alpha_s}(t_{s'}) \tilde{u},
\]
where we set $x_{\alpha_i+\gamma}(\tilde{t}_i) := 1$ if
$\alpha_i+\gamma \not\in \Phi^+$. It follows from the uniqueness in
(\ref{eq:elemu}) that $\alpha_i+\gamma \not\in \Phi^+$ or 
$\alpha_i+\gamma \in N$ for $i=1, \dots, s$ and therefore 
$\alpha_1, \dots, \alpha_s \in Z$. By induction we get
$\alpha_{s+1}, \dots, \alpha_{s'} \in Z$ completing the proof.
\epr

\begin{defn} \label{def:rzpatquot}
Let $N \subseteq S \subseteq \Phi^+$ be closed patterns such that
$N \unlhd S$. Assume that Hypothesis~\ref{hyp:p} holds and let $Z$ be 
the closed subpattern of $S$ defined in Lemma~\ref{la:cenquopat}. We call 
the set $\rz(P(S)/P(N)) := Z \setminus N$ the \emph{root center} of
$P(S)/P(N)$. 
\end{defn}


\section{Root kernels, root centers and hooks}
\label{sec:rkeretc}

We keep the setting from the previous sections. In particular, $q$ is
a power of a prime $p$ and $Y_n(q)$ is an untwisted Chevalley group
defined over~$\F_q$ with root system~$\Phi$. In this section we
associate with each irreducible character $\chi$ of a pattern
subgroup~$P(S)$ of $UY_n(q)$ certain sets of roots and pattern
subgroups of $P(S)$.

\medskip

\noindent\emph{We assume throughout this section that Hypothesis~\ref{hyp:p}
  holds.} 

\smallskip

\begin{defn} \label{def:rootker}
Let $S \subseteq \Phi^+$ be a closed pattern. For $\chi \in \Irr(P(S))$ we
set
\begin{eqnarray*}
\rk(\chi) & := & \{ \alpha \in S \mid X_\alpha \subseteq \Ker(\chi) \},\\
\rz(\chi) & := & \{ \alpha \in S \mid X_\alpha \subseteq Z(\chi) \},\\
\rs(\chi) & := & \rz(\chi) \setminus \rk(\chi),
\end{eqnarray*}
and call $\rk(\chi)$ the \emph{root kernel}, $\rz(\chi)$ the
\emph{root center} and $\rs(\chi)$ the \emph{central root support} of
$\chi$. Associated with these sets of roots are the following groups:
\begin{eqnarray*}
\rKer(\chi) & := & \la X_\alpha \mid \alpha \in \rk(\chi) \ra
\subseteq P(S) ,\\
\rZ(\chi) & := & \la X_\alpha \mid \alpha \in \rz(\chi) \ra \subseteq
P(S) ,\\
\rS(\chi) & := & \rZ(\chi) / \rKer(\chi) \subseteq P(S) / \rKer(\chi).
\end{eqnarray*}
\end{defn}

The next lemma shows that root kernels and root centers behave in much
the same way as the usual kernels and centers of irreducible
characters. It also implies that $\rS(\chi)$ is indeed a factor group.

\begin{lem} \label{la:rootker}
Let $S \subseteq \Phi^+$ be a closed pattern and 
$\chi \in \Irr(P(S))$. Then: 
\begin{enumerate}
\item[(a)] $\rk(\chi)$ and $\rz(\chi)$ are closed patterns which are
  normal in $S$.

\item[(b)] $\rKer(\chi)$ and $\rZ(\chi)$ are normal subgroups of
  $P(S)$. 

\item[(c)] The factor groups $P(S) / \rKer(\chi)$,
  $P(S) / \rZ(\chi)$ are quotient pattern groups.

\item[(d)] $\rZ(\chi)/\rKer(\chi) = Z(P(S)/\rKer(\chi))$.
\end{enumerate}
\end{lem}

\pr
(a) Let $\alpha \in S$ and $\beta \in \rk(\chi)$ such that 
$\alpha+\beta \in \Phi^+$. By definition, we have 
$X_\beta \subseteq \Ker(\chi) \unlhd P(S)$. Thus
Lemma~\ref{la:patnorm}~(c) implies $X_{\alpha+\beta} \subseteq
[X_\alpha, X_\beta] \subseteq \Ker(\chi)$ and so
$\alpha+\beta \in \rk(\chi)$. It follows that $\rk(\chi)$ is a closed
pattern which is normal in~$S$. The proof for $\rz(\chi)$ is analogous.

\medskip

(b), (c) follow from (a) and Remark~\ref{rem:normpat}.

\medskip

(d) By (a) and Lemma~\ref{la:cenquopat} there is a closed pattern $Z$ 
such that $\rk(\chi) \subseteq Z \subseteq S$ and
$P(Z)/\rKer(\chi) = Z(P(S)/\rKer(\chi))$. Obviously, 
$Z \subseteq \rz(\chi)$. Let $\alpha \in S$ and $\beta \in \rz(\chi)$. 
By Lemma~\ref{la:patnorm} (c) we have $X_{\alpha+\beta} \subseteq
[X_\alpha, X_\beta] \subseteq [X_\alpha, Z(\chi)] \subseteq \Ker(\chi)$.
Hence $\alpha+\beta \in \rk(\chi) \unlhd S$ and
Lemma~\ref{la:patnorm}~(b) implies that 
$i\alpha+j\beta \in \rk(\chi)$ for all positive integers $i,j$.
Hence, $[X_\alpha, X_\beta] = 1$ modulo $\rKer(\chi)$ by
Lemma~\ref{la:patnorm}~(c) and we can conclude $\rz(\chi) \subseteq Z$.
\epr

\medskip

Without Hypothesis~\ref{hyp:p} the statements in
Lemma~\ref{la:rootker} are not always true. Suppose that
$\Phi$ has type $B_2$ and that $\{ \alpha, \beta \}$ is a set of
simple roots, where $\alpha$ is short. If $p=q=2$ then $N :=
\{x_\alpha(d_1) x_\beta(d_2) x_{\alpha+\beta}(d_3)
x_{2\alpha+\beta}(d_3) \mid d_1, d_2, d_3 \in \F_q\}$ is a normal
subgroup of $UB_2(2)$ of index 2. Let $\chi \in \Irr(UB_2(2))$
be the nontrivial linear character with $\Ker(\chi) = N$. Then
$\rk(\chi) = \{\alpha, \beta\}$ and this is no closed pattern.

\begin{lem} \label{la:rs_almostfaith}
Let $N \unlhd \Phi^+$ and 
$\Sigma = \{\alpha_1, \dots, \alpha_s\} = \rz(UY_n/P(N))$. 
\begin{enumerate}
\item[(a)] The normal closed pattern $N$ is the unique maximal element
  of the set $\{M \unlhd \Phi^+ \mid M \cap \Sigma = \emptyset\}$.

\item[(b)] For all irreducible characters $\chi \in \Irr(UY_n)$ the
  following are equivalent:
  \begin{enumerate}
  \item[(i)] $\rs(\chi) = \Sigma$.

  \item[(ii)] $\chi$ is a constituent of $\lambda^{UY_n}$ for some
    linear character $\lambda$ of $P(\Sigma \cup N)$ satisfying 
    $P(N) \subseteq \Ker(\lambda)$ and 
    $\lambda|_{X_{\alpha_i} } \neq \mathbf{1}_{X_{\alpha_i}}$ for $i=1,2,\dots,s$.
  \end{enumerate}
  If the conditions (i) and (ii) hold then $N = \rk(\chi)$.
\end{enumerate}
\end{lem}

\pr
(b) By definition we have $N \unlhd \Phi^+$ and $\Sigma \cap N = \emptyset$.
We identify $X_{\alpha_i}$ with its image in $UY_n/P(N)$, so 
$Z(UY_n/P(N)) = X_{\alpha_1} \times \cdots \times X_{\alpha_s}$.

(i) $\Rightarrow$ (ii) Suppose that $\rs(\chi) = \Sigma$ and let 
$M$ be a normal closed pattern with $M \cap \Sigma = \emptyset$. If 
$M \not\subseteq \rk(\chi)$ then we  choose 
$\beta \in M \setminus \rk(\chi)$ such that $\height(\beta)$ is
maximal. Let $\alpha \in \Phi^+$. Since $M \unlhd \Phi^+$ 
the maximality of $\height(\beta)$ implies that 
$\alpha+\beta \not\in \Phi^+$ or $\alpha+\beta \in \rk(\chi)$ 
and therefore $\beta \in \rz(\chi)$ by Lemma~\ref{la:cenquopat} 
and Lemma~\ref{la:rootker}~(d). Thus 
$\beta \in \rz(\chi) \setminus \rk(\chi) = \rs(\chi) = \Sigma$, 
a contradiction. Hence $M \subseteq \rk(\chi)$. In particular, we 
have $N \subseteq \rk(\chi)$ and we can identify $\chi$ with some
$\tilde{\chi} \in \Irr(UY_n/P(N))$. Because $\rs(\chi) = \Sigma$ we
have $X_{\alpha_i} \not\subseteq \Ker(\tilde{\chi})$ for all $i$. 
Thus there is a linear character $\tilde{\lambda} \in
\Irr(X_{\alpha_1} \times \cdots \times X_{\alpha_s})$ with
$\tilde{\lambda}|_{X_{\alpha_i} } \neq \mathbf{1}_{X_{\alpha_i}}$ 
for all $i$ such that $\tilde{\chi}$ is a constituent of 
$\tilde{\lambda}^{UY_n/P(N)}$. Hence 
$\lambda := \infl_{P(\Sigma \cup N)/P(N)}^{P(\Sigma \cup N)} \tilde{\lambda}$ 
has properties described in (ii).

(ii) $\Rightarrow$ (i) Let $\chi$ be a constituent of $\lambda^{UY_n}$
for some linear character $\lambda$ of $P(\Sigma \cup N)$ with
$P(N) \subseteq \Ker(\lambda)$ and 
$\lambda|_{X_{\alpha_i} } \neq \mathbf{1}_{X_{\alpha_i}}$ for $i=1,2,\dots,s$. 
By construction, we have $N \subseteq \rk(\chi)$ and 
$\Sigma \cap \rk(\chi) = \emptyset$. Assume that $N \neq \rk(\chi)$. 
Then there is $\beta \in \rk(\chi) \setminus N$ such that
$\height(\beta)$ is maximal. Let $\alpha \in \Phi^+$. Since 
$\rk(\chi) \unlhd \Phi^+$ the maximality of $\height(\beta)$ implies
that $\alpha+\beta \not\in \Phi^+$ or $\alpha+\beta \in N$ and
therefore $\beta \in \rz(UY_n/P(N))$ by Lemma~\ref{la:cenquopat}. 
Hence $\beta \in \Sigma$, a contradiction. Thus $\rk(\chi) = N$.
Lemma~\ref{la:rootker} (d) implies that $\rs(\chi) = \Sigma$ so that
(i) holds.

\medskip

(a) Let $\chi$ be an irreducible constituent of $\lambda^{UY_n}$ where
$\lambda$ is a linear character of $P(\Sigma \cup N)$ as in (b) (ii).
We have already seen in the proof of (i) $\Rightarrow$ (ii) that $N$
is an element of the set 
$\{M \unlhd \Phi^+ \mid M \cap \Sigma = \emptyset\}$ and that each
element $M$ of this set is a subset of $\rk(\chi)$. In the proof of 
(ii) $\Rightarrow$ (i) we showed that $\rk(\chi)=N$. This completes the
proof of the lemma.
\epr


\section{Representable sets}
\label{sec:representable}

We keep the setting from the previous sections. In particular, $q$ is
a power of a prime $p$ and $Y_n(q)$ is an untwisted Chevalley group
defined over~$\F_q$ with root system~$\Phi$. In this section we
use the central root support $\rs(\chi)$ of characters $\chi$ of
a pattern subgroup~$P(S)$ of $UY_n(q)$ to obtain a partition of
$\Irr(P(S))$ which is well-adapted to the Lie theoretic setting. 

\medskip

\noindent\emph{We assume throughout this section that Hypothesis~\ref{hyp:p}
  holds.} 

\smallskip

\begin{defn} \label{def:repres}
Let $\Sigma \subseteq \Phi^+$. We say that the set $\Sigma$ is
\emph{representable} if there exists $\chi \in \Irr(UY_n)$ such
that $\rs(\chi) = \Sigma$. In this case we define
\[
\Irr(UY_n)_\Sigma := \{ \chi \in \Irr(UY_n) \mid \rs(\chi) = \Sigma\}.  
\]
\end{defn}

Obviously, $\Irr(UY_n)$ is partitioned by the sets
$\Irr(UY_n)_\Sigma$ where $\Sigma$ ranges over the representable
subsets of $\Phi^+$. The next lemma gives a characterization of
the representable sets.

\begin{lem} \label{la:representable}
For a subset $\Sigma \subseteq \Phi^+$ the following are equivalent:
\begin{enumerate}
\item[(a)] The set $\Sigma$ is representable. 
\item[(b)] There is a closed pattern $N \unlhd \Phi^+$ such that
  $\Sigma = \rz(UY_n/P(N))$.
\item[(c)] There is a unique closed pattern $N \unlhd \Phi^+$ such that
  $\Sigma = \rz(UY_n/P(N))$.
\end{enumerate}
\end{lem}
\pr
(a) $\Rightarrow$ (b) Suppose that $\Sigma$ is
representable and let $\chi \in \Irr(UY_n)$ such that 
$\rs(\chi) = \Sigma$. We set  $N := \rk(\chi)$. Then 
$\Sigma = \rz(UY_n/P(N))$ by Lemma~\ref{la:rootker} (a), (d) with
$S=\Phi^+$.

\medskip

(b) $\Rightarrow$ (c) follows from Lemma~\ref{la:rs_almostfaith} (a).

\medskip

(c) $\Rightarrow$ (a) Choose a linear character
$\lambda \in P(\Sigma \cup N)$ as in Lemma~\ref{la:rs_almostfaith}~(b)
(ii) and let $\chi \in \Irr(UY_n)$ be a constituent of $\lambda^{UY_n}$. 
By Lemma~\ref{la:rs_almostfaith}~(b) we have $\rs(\chi) = \Sigma$ and
(a) follows.
\epr

\medskip

\begin{defn} \label{def:kSigma_almfaith}
Let $\Sigma \subseteq \Phi^+$ be representable.
\begin{enumerate}
\item[(a)] We write $k(\Sigma)$ for the unique normal closed pattern with
  \[
  \Sigma = \rz(UY_n/P(k(\Sigma))).
  \]

\item[(b)] We say that a character $\chi \in \Irr(UY_n/P(k(\Sigma)))$ 
  is \emph{almost faithful} with respect to $\Sigma$ if 
  $X_\alpha \not\subseteq \Ker(\chi)$ for all $\alpha \in \Sigma$.  
\end{enumerate}
\end{defn}

The following remark is an immediate consequence of
Lemmas~\ref{la:rootker} (d) and \ref{la:representable}.

\begin{rem} \label{rem:kSigma}
Let $\Sigma \subseteq \Phi^+$ be representable. Then:
\begin{enumerate}
\item[(a)] $k(\Sigma)$ is the largest element of the set 
  $\{M \unlhd \Phi^+ \mid M \cap \Sigma = \emptyset\}$.

\item[(b)] Inflation induces a one to one correspondence between the
  set of almost faithful irreducible characters of $UY_n/P(k(\Sigma))$
  and $\Irr(UY_n)_\Sigma$.

\item[(c)] $\Irr(UY_n)_\Sigma = \{\chi\in\Irr(UY_n) \mid \rk(\chi) = k(\Sigma)\}$.
\end{enumerate}
\end{rem}

\medskip

Next, we investigate some elementary properties of $k(\Sigma)$.

\begin{lem} \label{la:reprker}
Let $\Sigma \subseteq \Phi^+$ be representable. For all roots
$\alpha \in \Sigma$ and $\gamma \in \Phi^+$ we have 
$\alpha+\gamma \not\in \Phi^+$ or $\alpha+\gamma \in k(\Sigma)$.
\end{lem}
\pr
This follows from Lemmas~\ref{la:cenquopat} and \ref{la:rootker} 
(a), (d). 
\epr

\medskip

\begin{defn} \label{def:n}
For $\alpha \in \Phi^+$ let $n_0(\alpha) := \{ \alpha \}$ and for
$i \ge 1$ 
\[
n_i(\alpha) := \{\beta+\gamma \, | \, \beta \in n_{i-1}(\alpha),
\gamma \in \Phi^+\} \cap \Phi^+.
\]
Finally, we set $n(\alpha) := \bigcup_{i \ge 1} n_i(\alpha)$.
\end{defn}

Note that $\alpha \not\in n(\alpha)$ since the union is over all
$i \ge 1$. In fact, $n(\alpha)$ is the normal closure of $\{\alpha\}$
in $\Phi^+$ with $\alpha$ removed. Also note that by construction,
for all $\beta \in n(\alpha)$ and $\gamma \in \Phi^+$ we have 
$\beta + \gamma \not\in \Phi^+$ or $\beta + \gamma \in n(\alpha)$;
hence $n(\alpha) \unlhd \Phi^+$.

\begin{lem} \label{la:nker}
If $\Sigma \subseteq \Phi^+$ is representable then
$\bigcup_{\alpha \in \Sigma}n(\alpha) \subseteq k(\Sigma)$.
\end{lem}
\pr
Let $\alpha \in \Sigma$. By Lemma~\ref{la:reprker} we have 
$n_1(\alpha) \subseteq k(\Sigma)$. Suppose that $i>1$. Let
$\gamma\in\Phi^+$ and $\beta \in n_{i-1}(\alpha)$ such that
$\beta+\gamma\in\Phi^+$. By induction, $\beta \in k(\Sigma)$. Since
$k(\Sigma) \unlhd \Phi^+$ we get $\beta+\gamma \in k(\Sigma)$ and the
claim follows. 
\epr

\medskip

In the following we show that for each $\alpha \in \Phi^+$ the set
$\{\alpha\}$ is representable and we obtain a recursive description of 
$k(\{\alpha\})$. 

\begin{defn} \label{def:k}
Let $\alpha \in \Phi^+$. We define $k_0(\alpha) := n(\alpha)$ and for 
$i \ge 1$
\[
k_i(\alpha) := (\rz(UY_n/P(k_{i-1}(\alpha)) \setminus \{\alpha\}) \cup
k_{i-1}(\alpha). 
\]
Finally, we set $k(\alpha) := \bigcup_{i \ge 0} k_i(\alpha)$. 
\end{defn}

Note that $k_i(\alpha) \unlhd \Phi^+$ for all $i \ge 0$. This is clear
for $i=0$ since $n(\alpha) \unlhd \Phi^+$. Suppose that $i>0$. By induction 
$k_{i-1}(\alpha) \unlhd \Phi^+$. Let $\beta \in k_i(\alpha)$ and
$\gamma \in \Phi^+$ such that $\beta+\gamma \in \Phi^+$. Thus 
$\beta \in \rz(UY_n/P(k_{i-1}(\alpha))$ or $\beta \in k_{i-1}(\alpha)$. 
In both cases we have 
$\beta+\gamma \in k_{i-1}(\alpha) \subseteq k_i(\alpha)$. Hence
$k_i(\alpha) \unlhd \Phi^+$. It follows that $k(\alpha)$ is
well-defined and that $k(\alpha) \unlhd \Phi^+$.

\begin{lem} \label{la:alpha_rep}
For all $\alpha \in \Phi^+$ the set $\{\alpha\}$ is representable and
$k(\{\alpha\}) = k(\alpha)$.
\end{lem}
\pr
We have just seen that $k(\alpha) \unlhd \Phi^+$ and by
definition of $k(\alpha)$ we have $\alpha \not\in k(\alpha)$. 
Since $n_1(\alpha) \subseteq n(\alpha) \subseteq k(\alpha)$
we have $\{\alpha\} \subseteq \rz(UY_n/P(k(\alpha)))$. Suppose that
$\beta \in \rz(UY_n/P(k(\alpha)))$. Then $\beta+\gamma \not\in \Phi^+$
or $\beta+\gamma \in k(\alpha)$ for all positive roots $\gamma$.
By the definition of $k(\alpha)$ we have $\beta \in k_i(\alpha)$ for
some $i$ or $\beta=\alpha$. Because $\rz(UY_n/P(k(\alpha))) \cap
k(\alpha) = \emptyset$ we have $\beta=\alpha$. Hence 
$\{\alpha\} = \rz(UY_n/P(k(\alpha)))$ so that condition (b) of
Lemma~\ref{la:representable} is satisfied with $N = k(\alpha)$ 
and the claim follows.

\epr

\medskip

Next, we define a set $w(\alpha) \subseteq \Phi^+ \setminus k(\alpha)$.
The definition of $w(\alpha)$ requires the concept of hooks. 

\begin{defn} \label{def:hook}
Let $\alpha \in \Phi^+$. We call the set
\[
h(\alpha) := \{\gamma \in \Phi^+ \, | \, \alpha - \gamma \in \Phi^+ \cup \{ 0 \} \}
\]
the \emph{hook} corresponding to $\alpha$. A subset
$h'(\alpha) \subseteq h(\alpha)$ is called a \emph{subhook} of
$h(\alpha)$ if $\alpha \in h'(\alpha)$ and 
$\alpha - \gamma \in h'(\alpha) \cup \{ 0 \}$ for all $\gamma \in h'(\alpha)$.
\end{defn}

The terminology \emph{hooks} is motivated by the case that $\Phi$ is 
of type~$A_n$; see~\cite[Section 3]{HML_D4}. The following
definition is in some sense \emph{dual} to that of $k(\alpha)$. 

\begin{defn} \label{def:w}
Let $\alpha \in \Phi^+$. We define $w_0(\alpha) := \{\alpha\}$ and for 
$i \ge 1$
\[
w_i(\alpha) := \bigcup_{\beta \in w_{i-1}(\alpha)} h(\beta).
\]
Finally, we set $w(\alpha) := \bigcup_{i \ge 0} w_i(\alpha)$. 
\end{defn}

The set $w(\alpha)$ has the following interpretation in terms of the
root poset $(\Phi^+, \preceq)$:

\begin{lem} \label{la:w}
For all $\alpha \in \Phi^+$ we have $w(\alpha) = \{\beta \in \Phi^+
\mid \beta \preceq \alpha\}$.
\end{lem}
\pr
$\subseteq$ Let $\beta \in w(\alpha)$. Then there is some $i$ such
that $\beta \in w_i(\alpha)$. We use induction on $i$ to show that
$\beta \preceq \alpha$. If $i=0$ we have $\beta=\alpha \preceq \alpha$
so that we can assume $i>0$. Hence there is a root 
$\gamma \in w_{i-1}(\alpha)$ such that $\beta \in h(\gamma)$ and by
induction we have $\beta \preceq \gamma \preceq \alpha$ and therefore
$\beta \preceq \alpha$.

\medskip

$\supseteq$ Let $\beta \in \Phi^+$ with $\beta \preceq \alpha$. We can
assume that $\beta \neq \alpha$. By Lemma~\ref{la:maxchain} there
are $\gamma_i \in \Phi^+$ such that $\alpha=\gamma_0 \succ \gamma_1
\succ \gamma_2 \succ \dots \succ \gamma_{s-1} \succ \gamma_s=\beta$
and $\gamma_i-\gamma_{i+1}$ is a simple root for all~$i$. In
particular, we have $\gamma_{i+1} \in h(\gamma_i)$ for all $i$. It
follows that $\gamma_i \in w_i(\alpha)$ and hence 
$\beta \in w_s(\alpha) \subseteq w(\alpha)$.
\epr

\medskip

If $\alpha_h \in \Phi^+$ is the highest root then Lemma~\ref{la:w} and
the remarks after Lemma~\ref{la:maxchain} imply that $w(\alpha_k)=\Phi^+$.

\begin{lem} \label{la:k_w}
For each $\alpha \in \Phi^+$ we have: $w(\alpha) \cap k(\alpha)=\emptyset$ 
and $w(\alpha) \cup k(\alpha)=\Phi^+$. 
\end{lem}
\pr 
We show $w_i(\alpha) \cap k(\alpha)=\emptyset$ for all $i \ge 0$ 
by induction on $i$. By definition of $n(\alpha)$ and $k(\alpha)$ we 
have $\alpha \not \in k(\alpha)$ so that we can assume $i>0$. 
Suppose that $\beta \in w_i(\alpha) \cap k(\alpha)$. By definition 
of $w_i(\alpha)$ there is $\gamma \in w_{i-1}(\alpha)$ such that 
$\beta \in h(\gamma)$. Since $\beta \in k(\alpha)$ and 
$k(\alpha) \unlhd \Phi^+$ we get $\gamma \in k(\alpha)$ which is
impossible by induction. This proves $w(\alpha) \cap k(\alpha) = \emptyset$.

\medskip

We now show that $m(\alpha) := \Phi^+ \setminus w(\alpha)$ is a normal
closed pattern. Let $\beta \in m(\alpha)$ and $\gamma \in \Phi^+$ with
$\beta+\gamma \in \Phi^+$. If $\beta+\gamma \in w(\alpha)$ then 
$\beta \prec \beta+\gamma \preceq \alpha$. Hence $\beta \preceq \alpha$
and $\beta \in w(\alpha)$, a contradiction. Thus $m(\alpha)$ is a
normal closed pattern not containing~$\alpha$. Since $\Sigma:=\{\alpha\}$ 
is representable the maximality property of $k(\alpha)=k(\{\alpha\})$
in Remark~\ref{rem:kSigma} (a) implies that $m(\alpha) \subseteq k(\alpha)$ 
and therefore $\Phi^+ \subseteq w(\alpha) \cup k(\alpha)$.
\epr

\medskip

In analogy with $k(\alpha)$ and $k(\Sigma)$ we want to replace
$\alpha$ in $w(\alpha)$ by arbitrary representable sets $\Sigma$.

\begin{defn} \label{def:wSigma}
For subsets $\Sigma \subseteq \Phi^+$ we define 
$w(\Sigma) := \bigcup_{\alpha \in \Sigma} w(\alpha)$.
\end{defn}

In particular, we have $w(\{\alpha\}) = w(\alpha)$. Statement (b) in
the following proposition generalizes Lemma~\ref{la:k_w} to arbitrary
representable sets $\Sigma$.

\begin{prop} \label{prop:kker}
If $\Sigma \subseteq \Phi^+$ is representable then 
\begin{enumerate}
\item[(a)] $k(\Sigma) = \bigcap_{\alpha \in \Sigma} k(\alpha)$,
\item[(b)] $w(\Sigma) \cap k(\Sigma)=\emptyset$ and 
$w(\Sigma) \cup k(\Sigma)=\Phi^+$.
\end{enumerate}
\end{prop}
\pr
Suppose that $\beta \in w(\Sigma) \cap k(\Sigma)$. Then there is
$\alpha \in \Sigma$ such that $\beta \in w(\alpha)$ and hence 
$\beta \prec \alpha$. By Lemma~\ref{la:maxchain} there are 
roots $\gamma_i \in \Phi^+$ such that $\beta=\gamma_0 \prec \gamma_1
\prec \gamma_2 \prec \dots \prec \gamma_{s-1} \prec \gamma_s=\alpha$
and $\gamma_i-\gamma_{i-1}$ is a simple root for all~$i$. It follows
that the roots $\gamma_i$ belong to the normal closure $M$ of 
$\{\beta\}$ in $\Phi^+$. Since $\beta \in k(\Sigma)$ and 
$k(\Sigma) \unlhd \Phi^+$ we get $\alpha=\gamma_s \in M \subseteq
k(\Sigma)$. Thus $\alpha \in \Sigma \cap k(\Sigma)$, a contradiction
to Remark~\ref{rem:kSigma} (a). Hence $w(\Sigma) \cap k(\Sigma)=\emptyset$.

From Definition~\ref{def:wSigma} and Lemma~\ref{la:k_w} we get
\begin{equation} \label{eq:kwSigma}
k(\Sigma) \subseteq \Phi^+ \setminus w(\Sigma) = \Phi^+ \setminus
(\bigcup_{\alpha \in \Sigma} w(\alpha)) = \bigcap_{\alpha \in \Sigma}
(\Phi^+ \setminus  w(\alpha)) = \bigcap_{\alpha \in \Sigma} k(\alpha). 
\end{equation}
For all $\alpha \in \Sigma$ we have $k(\alpha) \unlhd \Phi^+$ and
$\alpha \not\in k(\alpha)$. Thus we have
$\bigcap_{\alpha \in \Sigma} k(\alpha) \unlhd \Phi^+$ and
$\left(\bigcap_{\alpha \in \Sigma} k(\alpha)\right) \cap \Sigma = \emptyset$.
From Remark \ref{rem:kSigma} (a) we can conclude 
$\bigcap_{\alpha \in \Sigma} k(\alpha) \subseteq k(\Sigma)$. 
Hence we have equality in \eqref{eq:kwSigma} and the claims (a) and
(b) follow. 
\epr

\medskip

In the following we give an interpretation of representable sets in
terms of the root poset $(\Phi^+,\preceq)$. This allows us to deduce
the number of representable sets of size~$k$ for all $k$ from results
in~\cite{FR}. Recall that $\Phi$ is an irreducible root system of
rank~$n$. An antichain of the poset $(\Phi^+, \preceq)$ is a subset  
$A \subseteq \Phi^+$ such that the elements of $A$ are 
pairwise incomparable. A variation of the next proposition for
pattern subgroups if $\Phi$ is of type~$A_n$ was already obtained by
Isaacs in~\cite[Theorem~3.1]{IsaacsTriang}.

\begin{prop} \label{prop:antichain}
For a subset $\Sigma \subseteq \Phi^+$ the following are equivalent:
\begin{enumerate}
\item[(a)] $\Sigma$ is representable.
\item[(b)] $\Sigma$ is an antichain of the root poset $(\Phi^+,\preceq)$.  
\end{enumerate}
\end{prop}
\pr 
(a) $\Rightarrow$ (b) Suppose that $\Sigma$ is a representable set and
that $\alpha,\beta \in \Sigma$ with $\alpha \neq \beta$. We need to
show that $\alpha$ and $\beta$ are not comparable, so suppose otherwise.
Without loss we may assume $\alpha \prec \beta$. The remarks after
Definition~\ref{def:n} and Lemma~\ref{la:nker} imply that 
$\beta \in n(\alpha) \subseteq k(\Sigma)$, hence 
$k(\Sigma) \cap \Sigma \neq \emptyset$ contradicting
Remark~\ref{rem:kSigma} (a). Thus $\alpha$ and $\beta$ are
incomparable showing that every representable set is an antichain. 

\medskip

(b) $\Rightarrow$ (a) Let $\Sigma$ be an antichain. We claim that
$N := \Phi^+ \setminus w(\Sigma)$ is a normal closed pattern. Let
$\alpha \in \Phi^+$ and $\beta \in N$ with $\alpha+\beta \in \Phi^+$.
If $\alpha+\beta \not\in N$ then $\alpha+\beta \in w(\Sigma)$. Hence
there is $\gamma \in \Sigma$ such that $\alpha+\beta \in w(\gamma)$.
By Lemma~\ref{la:w} we have $\beta \preceq \alpha+\beta \preceq \gamma$.
Thus $\beta \preceq \gamma$ and then $\beta \in w(\gamma) \subseteq w(\Sigma)$
which is impossible since $\beta \in N = \Phi^+ \setminus w(\alpha)$. 
This shows that $N \unlhd \Phi^+$.

Let $\alpha$ be a maximal element of $w(\Sigma)$ with respect to
$\preceq$. The definition of $w(\Sigma)$ and Lemma~\ref{la:w} 
imply that $\alpha \in \Sigma$. Conversely, let $\alpha \in \Sigma$
and suppose that there is $\gamma \in w(\Sigma)$ such that 
$\alpha \prec \gamma$. By the definition of $w(\Sigma)$ there is
$\beta \in \Sigma$ such that $\gamma \in w(\beta)$ and hence 
$\gamma \preceq \beta$ by Lemma~\ref{la:w}. It follows that 
$\alpha \preceq \beta$ contradicting the fact that $\Sigma$ is an
antichain. Thus $\alpha$ is a maximal element of $w(\Sigma)$. This
shows that the elements of $\Sigma$ are exactly the maximal elements
of $w(\Sigma)$. We get from Lemma~\ref{la:cenquopat} that 
$\Sigma = \rz(UY_n/P(N))$ so that condition (b) of
Lemma~\ref{la:representable} is satisfied. Hence $\Sigma$ is
representable. 
\epr

\medskip

We can now apply the results in~\cite[Sections~5.1-5.2]{FR} on the
number of antichains in root posets.

\begin{prop} \label{prop:numrepsets}
Let $N_k(\Phi) := |\{\Sigma \subseteq \Phi^+ \mid \Sigma \ \text{is
  representable and } |\Sigma|=k\}|$ be the number of representable
  subsets of size $k$ of $\Phi^+$ and $N(\Phi) := \sum_k N_k(\Phi)$
  the number of all representable subsets of $\Phi$. Then:
\begin{enumerate}
\item[(a)] $N_k(\Phi)=0$ for all $k>n$, i.e., $|\Sigma| \le n$ for all
  representable subsets $\Sigma$ of~$\Phi^+$. 

\item[(b)] For all $k$ the number $N_k(\Phi)$ is the coefficient of
  $t^k$ in the polynomial $N(\Phi,t)$, where $N(\Phi,t)$ is given by
  Table~\ref{tab:Ng}. 

\begin{table}[!ht] 
\caption{Generating functions $N(\Phi,t)$ for the numbers $N_k(\Phi)$.} 
\label{tab:Ng}

\begin{center}
\begin{tabular}{lcl}
$N(A_n,t)$ & $=$ & $\sum_{k = 0}^n \frac{1}{n+1}\binom{n+1}{k} \binom{n+1}{k+1} t^k$ \\
& & \\
$N(B_n,t)$ & $=$ & $\sum_{k = 0}^n \binom{n}{k}^2  t^k$ \\
& & \\
$N(C_n,t)$ & $=$ & $\sum_{k = 0}^n \binom{n}{k}^2  t^k$ \\
& & \\
$N(D_n,t)$ & $=$ & $1 + t^n + \sum_{k = 1}^{n-1}  \left[\binom{n}{k}^2  - \frac{n}{n-1} \binom{n-1}{k-1}\binom{n-1}{k} \right] t^k$ \\
& & \\
$N(E_6,t)$ & $=$ &  $1 + 36t + 204t^2 + 351t^3 + 204t^4 + 36t^5 + t^6$ \\
& & \\
$N(E_7,t)$ &$=$  & $1 + 63t + 546t^2 + 1470t^3 + 1470t^4 + 546t^5 + 63t^6 + t^7$ \\
& & \\
$N(E_8,t)$ & $=$ & $1 + 120t + 1540t^2 + 6120t^3 + 9518t^4 + 6120t^5$ \\
& &  \ \ $+ 1540t^6 + 120t^7 + t^8$\\
& & \\
$N(F_4,t)$ &$=$  & $1 + 24t + 55 t^2 + 24t^3 + t^4$ \\
& & \\
$N(G_2,t)$ & $=$ &  $1 + 6t + t^2 $\\
\end{tabular}
\end{center}
\end{table}

\item[(c)] $N(\Phi)=\prod_{i = 1}^n \frac{e_i+h+1}{e_i+1}$, 
  where $e_1, e_2, \dots, e_n$ are the exponents of the Weyl group $W$ of
  $\Phi$ and $h$ is the Coxeter number of $W$. The numbers $N(\Phi)$
  for the various root systems are given by Table~\ref{tab:N}. 

\begin{table}[!ht] 
\caption{The numbers $N(\Phi)$.} 
\label{tab:N}

\begin{center}
\begin{tabular}{|l|l|l|l|l|l|l|l|}
\hline
$A_n$ & $B_n$, $C_n$ & $D_n$ & $E_6$ & $E_7$ & $E_8$ & $F_4$ & $G_2$\\
\hline 
$\frac{1}{n+2} \binom{2n+2}{n+1} $ & $\binom{2n}{n}$ & $\frac{3n-2}{n}
\binom{2n-2}{n-1}$ & $833$ & $4160$ & $25080$ & $105$ & $8$ \\ 
\hline
\end{tabular}

\end{center}
\end{table}
\end{enumerate}
\end{prop}

\pr
This follows from Proposition~\ref{prop:antichain} and 
\cite[Theorems 5.1 and 5.9]{FR}.

\epr

\medskip

The numbers $N_k(\Phi)$ are called \emph{(generalized) Narayana numbers} 
in \cite{FR}. According to \cite{FR} the numbers $N(\Phi)$ appeared
first in a paper by Djorkovi\'{c} \cite{D} on the enumeration of conjugacy
classes of elements of finite order in compact and complex semisimple
Lie groups. The numbers $N(\Phi)$ also count clusters in cluster
algebras and related objects; see Fomin and Zelevinsky \cite{FZ1},
\cite{FZ2}, \cite{FZ}, and Fomin and Reading \cite{FR} and the 
references therein. For other interesting connections between classical
combinatorial objects and the characters of $UA_n(q)$ see Marberg~\cite{Mar}.


\section{Single root characters}
\label{sec:single}

We keep the setting from the previous sections. In particular, $q$ is
a power of a prime $p$ and $Y_n(q)$ is an untwisted Chevalley group 
defined over~$\F_q$ with an irreducible root system~$\Phi$. The set of
positive roots is denoted by $\Phi^+$.

\medskip

\emph{We assume throughout this section that Hypothesis~\ref{hyp:p} holds.}

\medskip

For $N=\Phi^+$ we have $\rz(UY_n/P(N))=\emptyset$. It follows from 
Lemma~\ref{la:representable} that $\emptyset$ is representable and
that $\Irr(UY_n)_\emptyset = \{\mathbf{1}_{UY_n}\}$. 

In what follows we investigate $\Irr(UY_n)_\Sigma$ when $|\Sigma|=1$. 
We have already seen in Lemma~\ref{la:alpha_rep} that $\{\alpha\}$ is
representable for all $\alpha \in \Phi^+$ and we gave a recursive
description of $k(\{\alpha\})=k(\alpha)$ in Definition~\ref{def:k}.

\begin{defn} \label{def:single}
For $\alpha \in \Phi^+$ we set $\Irr(UY_n)_\alpha := \Irr(UY_n)_{\{\alpha\}}$. 
We say that a character $\chi \in \Irr(UY_n)_\alpha$ is a \emph{single
root character lying over $\alpha$} and that it is \emph{almost
faithful at $\alpha$}. A \emph{midafi character} or just \emph{midafi}
for $\alpha$ is a  character $\chi \in \Irr(UY_n)_\alpha$ such that 
$\chi(1) = \min\{\psi(1) \mid \psi \in \Irr(UY_n)_\alpha\}$. We set
\[
\Irr^\mida(UY_n)_\alpha := \{\chi \in \Irr(UY_n)_\alpha \mid \chi
\text{ is a midafi for } \alpha\}.
\]
\end{defn}

\medskip

Note that $q$ is fixed in Definition~\ref{def:single}. The term midafi
is an abbreviation for {\it minimal degree almost faithful
  irreducible}. The origin of this terminology is the
thesis~\cite{TungPhD} of the second author where the term was first
used in the case that $\Phi$ has Dynkin type $A_n$. To study single
root characters we will use the concept of arms and legs. 

\begin{defn} \label{def:armleg}
Let $\alpha \in \Phi^+$ and $a(\alpha), \ell(\alpha) \subseteq h(\alpha)$. 
\begin{enumerate}
\item[(a)] We say that $a(\alpha)$ is an \emph{arm} and
$\ell(\alpha)$ is the corresponding \emph{leg} of~$h(\alpha)$
if the hook $h(\alpha)$ is the disjoint union
$h(\alpha) = \{\alpha\} \cup a(\alpha) \cup \ell(\alpha)$ and
\begin{itemize}
\item[1)] $|a(\alpha)|=|\ell(\alpha)|$ and 
\item[2)] for each $\beta \in a(\alpha)$ there is a unique
  $\gamma \in \ell(\alpha)$ with $\beta+\gamma = \alpha$.
\end{itemize}

\item[(b)] If $a(\alpha)$ is an arm of $h(\alpha)$ we call 
the set $s(\alpha) := \Phi^+ \setminus a(\alpha)$ of roots the
\emph{source} corresponding to $\alpha$. If additionally
$s(\alpha)$ is a closed pattern then we call $S_\alpha:=
P(s(þ\alpha))$, the \emph{source group} corresponding to $\alpha$. 
\end{enumerate}
\end{defn}

The terminology of hooks, arms and legs is motivated by the case
that $\Phi$ is an irreducible root system of type $A_n$;
see~\cite[Section~3.3]{HML_D4} and \cite[Section~2.3]{TungPhD}. 

Note that Lemma~\ref{la:k_w} implies that we always have 
$k(\alpha) \subseteq s(\alpha)$.
We point out that the definition above allows for $2^m$ choices for
the set $a(\alpha)$, where $m = (|h(\alpha)| -1)/2$. Relative to a
choice of $a(\alpha)$ the sets $\ell(\alpha)$ and $s(\alpha)$ are
unique. We discuss this in more detail at the end of this section. The
next lemma is needed for the classification of the midafi characters when
$UY_n$ is of exceptional type and $n \ge 4$. 

\begin{defn} \label{def:armlegsub}
Let $\alpha \in \Phi^+$. We call a subset $a'(\alpha)$ of a subhook
$h'(\alpha)$ of $h(\alpha)$ an \emph{arm} of $h'(\alpha)$ if
there is an arm $a(\alpha)$ of $h(\alpha)$ with $a'(\alpha) =
a(\alpha) \cap h'(\alpha)$. In this case we call 
$\ell'(\alpha) := h'(\alpha) \setminus (a'(\alpha) \cup \{\alpha\})$
the corresponding \emph{leg} of $h'(\alpha)$.
\end{defn}

The following technical lemma will be used in
Section~\ref{sec:exc_midafis}. A stronger version of this lemma for
type $A_n$ is implicitly contained in~\cite{TungPhD}. Recall that a 
$p$-group $H$ is \emph{special} if $Z(H) = [H,H] = \Phi(H)$. If
$H$ is special with $|H| = q^{2m+1}$ and $Z(H)$ is elementary abelian of
order $q = p^a$, then we say that $H$ is \emph{special of type $q^{1+2m}$}. 

\begin{lem} \label{la:twohook}
Let $\alpha, \beta \in \Phi^+$. Let $h(\alpha)$ be the hook
corresponding to $\alpha$ with arm~$a(\alpha)$ and leg~$\ell(\alpha)$,
let $h'(\beta)$ be a subhook of $h(\beta)$ with arm $a'(\beta)$ and
leg $\ell'(\beta)$ and let~$h'_{\alpha\beta}$ be the closed pattern
generated by $h(\alpha) \cup h'(\beta) \cup k(\alpha)$. We define 
$\ell'_{\alpha\beta} := \{\alpha, \beta\} \cup \ell(\alpha) \cup k(\alpha)$. 
Assume that the following holds:
\begin{enumerate}
\item[(a)] $h(\alpha) \cup k(\alpha)$, $h'(\beta) \cup k(\alpha)$,
  $\Phi^+ \setminus a(\alpha)$ are closed patterns.

\item[(b)] $h'_{\alpha\beta} \setminus a(\alpha)$ normalizes $\ell(\alpha) \cup k(\alpha)$.

\item[(c)] $(h(\alpha) \cup k(\alpha)) \cap h'(\beta) = \emptyset$.

\item[(d)] The group $H_\alpha := P(h(\alpha) \cup
  k(\alpha))/P(k(\alpha))$ is special of type $q^{1+2|a(\alpha)|}$
  with $[x, H_\alpha] = Z(H_\alpha)$ for all $x \in H_\alpha \setminus Z(H_\alpha)$.

\item[(e)] The group $H'_\beta := P(h'(\beta) \cup
  k(\alpha))/P(k(\alpha))$ is special of type $q^{1+2|a'(\beta)|}$
  with $[x, H'_\beta] = Z(H'_\beta)$ for all $x \in H'_\beta \setminus Z(H'_\beta)$.
\end{enumerate}
Set $H'_{\alpha\beta} := P(h'_{\alpha\beta})/P(k(\alpha))$ and 
$L'_{\alpha\beta} := P(\ell'_{\alpha\beta})/P(k(\alpha))$ and let
$\mu \in \Irr(L'_{\alpha\beta})$ with 
$X_\alpha, X_\beta \not\subseteq \Ker(\mu)$ and 
$\prod_{\gamma \in \ell(\alpha)} X_\gamma \subseteq \Ker(\mu)$. 
Then each character $\psi \in \Irr(H'_{\alpha\beta}, \mu)$ has degree 
$\psi(1) \ge  q^{|a(\alpha)|+|a'(\beta)|}$.
\end{lem}

\pr
To apply the Reduction Lemma~\ref{nred} we introduce the following
notation (only for this proof):
\begin{itemize}
\item $U := H'_{\alpha\beta}$,

\item $H := P(h'_{\alpha\beta} \setminus a(\alpha))/P(k(\alpha))$,

\item $X := \prod_{\gamma \in a(\alpha)} X_\gamma$,

\item $Y := P(\ell(\alpha) \cup k(\alpha))/P(k(\alpha))$,

\item $Z := X_\alpha$,

\item $\lambda := \mu_Z$ and $\tilde{\lambda} := \mu_{ZY}$.
\end{itemize}
Since $k(\alpha)$ is a normal closed pattern $U$ is a quotient pattern
group and it follows from (a) that $H$ is a subgroup of~$U$. We know
from Lemma~\ref{la:alpha_rep} and Definition~\ref{def:kSigma_almfaith}~(a) 
that $Z$ is a subgroup of $U$ with $Z \subseteq Z(U)$. It follows 
from assumption (b) that $Y \unlhd H$ and we have $Z \cap Y = \{1\}$. 
By (\ref{eq:elemu}) the set $X$ is a set of representatives
for $U/H$. From assumption (b) and Lemma~\ref{la:derquopat} we get 
$ZY \unlhd U$. Since $\prod_{\gamma \in \ell(\alpha)} X_\gamma \subseteq \Ker(\mu)$ 
we have $Y \subseteq \Ker(\tilde{\lambda})$ and by (a), (d) we have 
${}^u\tilde{\lambda}\neq \tilde{\lambda}$ for all $u \in X \setminus \{1\}$. 
We also have $|X|=|Y|$. Thus we can apply Lemma~\ref{nred}. 

Let $\psi \in \Irr(H'_{\alpha\beta}, \mu)$. We consider the
restriction $\psi_H$ of $\psi$ to $H$. Since $Y \unlhd H$ the
character $\psi_H$ has a constituent $\chi \in \Irr(H/Y, \lambda)$
with $X_\beta \not\subseteq \Ker(\chi)$. We define
$H'_\beta := P(h'(\beta) \cup k(\alpha))/P(k(\alpha))$. 
Note that assumption (c) implies that $H'_\beta$ is a subgroup of $H$.
It follows from (e) that $\chi(1) = \chi_{H'_\beta}(1) \ge q^{|a'(\beta)|}$. 
By Frobenius reciprocity $\psi$ is a constituent of $\chi^U$. Since
$\chi^U$ is irreducible by Lemma~\ref{nred} we have $\psi = \chi^U$
and hence $\psi(1) = \chi^U(1) = q^{|a(\alpha)|} \cdot \chi(1) \ge
q^{|a(\alpha)|+|a'(\beta)|}$.  
\epr

\bigskip

Using the theory developed in Sections~\ref{sec:nota}-\ref{sec:single}
we can now study the single root characters of~$UY_n$ for all
irreducible root systems $\Phi$. In the following two sections we
give a proof of Theorems~\ref{thm:main_classical}-\ref{thm:main_exceptionalnormal} 
on a case-by-case basis and construct the midafi characters in
$\Irr(UY_n)_\alpha$ for all Dynkin types~$Y$, all ranks~$n$ and all
positive roots~$\alpha$. We now outline our approach, details
will be given in Sections~\ref{sec:single_classical} and
\ref{sec:exc_midafis}.

Let $\Phi$ be an irreducible root system and $\alpha \in \Phi^+$. 
By construction, the size of the corresponding hook is of the 
form $|h(\alpha)| = 1+2k$ for some nonnegative integer~$k$. 
Hence, there are~$2^k$ possible choices for the arm $a(\alpha)$ 
and the leg~$\ell(\alpha)$, and in general, we could not find a 
canonical choice. An important fact, which will be proved in 
Sections~\ref{sec:single_classical} and \ref{sec:exc_midafis}, is 
that one can always choose an arm $a(\alpha)$ of $h(\alpha)$ with 
the following property: 
\begin{enumerate}
\item The source $s(\alpha) = \Phi^+ \setminus a(\alpha)$ is a 
  closed pattern. 
\end{enumerate}
In this case we can work with the source group 
$S_\alpha = P(s(\alpha))$ and it follows from Lemma~\ref{la:k_w} that 
$\ell(\alpha) \cup k(\alpha) \subseteq s(\alpha)$. For us a good
choice for $a(\alpha)$ is one in which condition (1) and the 
condition
\begin{enumerate}
\item[(2)] The set $\ell(\alpha) \cup k(\alpha)$ is a normal closed
  subpattern of $s(\alpha)$.
\end{enumerate}
are met. 
If for a root $\alpha \in \Phi^+$ and for some choice of $a(\alpha)$ we 
achieve conditions (1) and (2), then we define 
\[
\overline{T}_\alpha := S_\alpha/P( \{\alpha \} \cup \ell(\alpha) \cup k(\alpha)).
\]
We now show that conditions (1), (2) together with conditions (3) and (4) below 
suffice to establish the
description of $\Irr(UY_n)_\alpha$ in terms of character
correspondences as stated in Theorems \ref{thm:main_classical} and
\ref{thm:main_exceptionalnormal}. 

\begin{prop} \label{prop:red_to_thm} 
Let $\Phi$ be a root system of type $Y_n$. Suppose that for some root 
$\alpha \in \Phi^+$ of $\height(\alpha) > 1$ and some choice of $a(\alpha)$
the following are true: 
\begin{enumerate}
\item[(1)] The source $s(\alpha)$ is a closed pattern.

\item[(2)] $\ell(\alpha) \cup k(\alpha) \unlhd s(\alpha)$.

\item[(3)] $h(\alpha) \cup k(\alpha)$ is a closed pattern and
  $a(\alpha)$ normalizes $\{\alpha\} \cup \ell(\alpha) \cup k(\alpha)$.

\item[(4)] For each element 
  $y \in (\prod_{\gamma \in a(\alpha)} X_\gamma) \setminus \{1\}$
  there is some root $\beta \in \ell(\alpha)$ such that   
  $\{[y,x_\beta(t)] \mid t \in \F_q\}=X_\alpha$.
\end{enumerate}
Then the map
$\Psi_\alpha: \Irr(\overline{T}_\alpha) \times \Irr(X_\alpha)^*
\to \Irr(UY_n)_\alpha$, with 
\[
(\mu, \lambda) \mapsto (\infl_{\overline{T}_\alpha}^{S_\alpha} \mu \cdot
 \infl_{X_\alpha}^{S_\alpha} \lambda)^{UY_n}
\]
is a one to one correspondence with the property
$\Psi_\alpha(\mu, \lambda)(1) = q^{|a(\alpha)|} \cdot \mu(1)$.
Moreover
\[
S_\alpha/P(\ell(\alpha) \cup k(\alpha)) \cong \overline{T}_\alpha \times X_\alpha.
\] 
\end{prop}

\pr
We apply the Reduction Lemma~\ref{nred}. The role of $U$ is played by
the quotient pattern group $UY_n/P(k(\alpha))$ and the role of $H$ is
played by $S_\alpha/P(k(\alpha))$. The role of $X$ is played by
$\prod_{\gamma \in a(\alpha)} X_\gamma$, the role of $Y$ is played by
$\prod_{\gamma \in \ell(\alpha)} X_\gamma$ and the role of $Z$ is
played by $X_\alpha$, where we identify the root subgroups $X_\gamma$
with their images in $UY_n/P(k(\alpha))$. 

Condition (a) of Lemma~\ref{nred} is satisfied by 
Definition~\ref{def:kSigma_almfaith} (a) and Lemma~\ref{la:alpha_rep}.
Hypothesis (2) implies that condition (b) of Lemma~\ref{nred} is
satisfied. Combining hypothesis (2) and hypothesis (3) implies
that $a(\alpha)$ normalizes $\{\alpha\} \cup \ell(\alpha) \cup k(\alpha)$  
and $\{\alpha\} \cup \ell(\alpha) \cup k(\alpha) \unlhd s(\alpha)$.
As $\Phi^+ = s(\alpha) \cup a(\alpha)$ we now see that
condition (d) of Lemma~\ref{nred} is satisfied. Conditions (c) and (f)
in Lemma~\ref{nred} hold by the definition of $a(\alpha)$ and 
$\ell(\alpha)$.

We still have to check condition (e) of Lemma~\ref{nred}. 
Let $\lambda \in \Irr(X_\alpha)^*$ and let $\tilde{\lambda}$ 
be the inflation of $\lambda$ to 
$ZY = X_\alpha \times P(\ell(\alpha) \cup k(\alpha))/P(k(\alpha))$.  
Suppose that there is $1 \neq x \in (\prod_{\gamma \in a(\alpha)} X_\gamma) \setminus \{1\}$
such that ${}^x\widetilde{\lambda} = \widetilde{\lambda}$. Because the
linear character $\lambda$ is nontrivial on $X_\alpha$ there is some $t \in \F_q$ 
such that $\widetilde{\lambda}(x_\alpha(t)) = \lambda(x_\alpha(t)) \neq 1$.
By hypothesis (4) there is $\beta \in \ell(\alpha)$ 
and $t' \in \F_q$ such that $[x, x_\beta(t')^{-1}] = x_\alpha(t)$. Hence
\[
{}^x\widetilde{\lambda}(x_\beta(t')) = 
\widetilde{\lambda}(x^{-1} x_\beta(t') x x_\beta(t')^{-1} x_\beta(t')) =
\widetilde{\lambda}( x_\alpha(t)) \widetilde{\lambda}( x_\beta(t'))
\neq \widetilde{\lambda}(x_\beta(t')),
\]
contradicting ${}^x\widetilde{\lambda} = \widetilde{\lambda}$. Thus we
have ${}^x\widetilde{\lambda} \neq \widetilde{\lambda}$ for all
$x \in (\prod_{\gamma \in a(\alpha)} X_\gamma) \setminus \{1\}$. Now
the correspondence follows from Lemma~\ref{nred}. 

Finally we observe that hypotheses (1), (2) combined with 
Lemmas~\ref{la:derquopat} and \ref{la:representable} imply
that
\[
S_\alpha/P(\ell(\alpha) \cup k(\alpha)) \cong \overline{T}_\alpha \times X_\alpha,
\] 
completing the proof.
\epr

\medskip

We remark that generally the hypotheses (3) and (4) are easily verified and do not present a 
serious obstacle. 

If $\Phi$ is a classical root system, i.e., $\Phi$ is of type
$A_n$, $B_n$, $C_n$ or $D_n$, then we show in Section \ref{sec:single_classical}
that it is always possible to make a
choice for the arm of a root such that 
the hypotheses of Proposition \ref{prop:red_to_thm} 
are satisfied and for these Dynkin types we are able to derive 
some information on the structure of $\overline{T}_\alpha$. 

However when $\Phi$ is of exceptional type and rank $\ge 4$, 
then there exist roots in $\Phi$ for which no choice of arm 
achieves the hypotheses of Proposition \ref{prop:red_to_thm}. 
For the exceptional types of root systems one can always choose
$a(\alpha)$ such that hypothesis (1) holds, but then, in general,
$\ell(\alpha) \cup k(\alpha)$ is no longer normal in $s(\alpha)$. 
In this case we replace $\ell(\alpha)$ by a subset
$\bar{\ell}(\alpha) \supsetneq \ell(\alpha)$ such that 
$\bar{\ell}(\alpha) \cup k(\alpha) \unlhd s(\alpha)$, 
$\{\alpha\} \cup \bar{\ell}(\alpha) \cup k(\alpha) \unlhd \Phi^+$
and such that the quotient pattern group 
$P(\bar{\ell}(\alpha) \cup k(\alpha))/P(k(\alpha))$ is abelian. 
Among all choices of the arm $a(\alpha)$ such that condition (1) holds
and such that $\bar{\ell}(\alpha)$ has the properties just described
we take one such that $|\bar{\ell}(\alpha)|$ is minimal. For this
choice we show, by restriction to a suitable subgroup, that each 
$\chi \in \Irr(UY_n)_\alpha$ with  
$\bar{\ell}(\alpha) \setminus \ell(\alpha) \not\subseteq \rk(\chi)$ 
has a degree greater than~$q^{|a(\a)|}$. The Reduction Lemma \ref{nred}
then gives a one to one correspondence
\[
\Psi: \Irr^\lin(\overline{T}_\alpha) \times \Irr(X_\alpha)^* \to \Irr^\mida(UY_n)_\alpha.
\]
Hence for all irreducible root systems $\Phi$ (the classical ones and 
the exceptional ones) and all positive roots $\alpha$ we obtain a one
to one correspondence  
\[
\Psi: \Irr^\lin(\overline{T}_\alpha) \times \Irr(X_\alpha)^* \to \Irr^\mida(UY_n)_\alpha. 
\]
This leads to a construction of all midafis in $\Irr(UY_n)_\alpha$ and 
allows us to compute their number and their degrees if $\Phi$ is of
exceptional type. For each positive root~$\alpha$ we call the 
$q-1$ irreducible characters of $UY_n$ in 
$\Psi(\{\mathbf{1}_{\overline{T}_\alpha}\} \times \Irr(X_\alpha)^*)$ the 
\emph{standard mifadis} corresponding to $\alpha$.


\section{Hook subgroups and midafis in classical groups}
\label{sec:single_classical}

In this section we prove Theorem~\ref{thm:main_classical}. We will use
the explicit construction of the root systems $\Phi$ of types $A_n$,
$B_n$, $C_n$ and $D_n$ given in~\cite[12.1]{Humphreys}. As before, we
use the following convention to keep the notation simple: Suppose 
that $N \subseteq S \subseteq \Phi^+$ are closed patterns such that 
$N \unlhd S$. Since for each $\gamma \in S \setminus N$ the root
subgroup $X_\gamma$ is mapped injectively into the quotient pattern
group $P(S)/P(N)$ we will often identify $X_\gamma$ with its image 
$X_\gamma P(N)$.  

\medskip

\noindent\emph{We assume throughout this section that Hypothesis~\ref{hyp:p}
  holds.}


\subsection{Type A} \label{subsec:An}

Let $n$ be a positive integer. We construct a root system of
type~$A_n$ as in~\cite[Section~12.1]{Humphreys}: Let 
$e_1, e_2, \dots, e_{n+1} \in \R^{n+1}$ be the usual orthonormal unit
vectors which form a basis of $\R^{n+1}$. Then
$\Phi := \{ e_i-e_j \, | \, 1 \le i \neq j \le n+1 \}$ is a root system
of type $A_n$ and the set $\{\alpha_1, \dots, \alpha_n\}$, where 
$\alpha_i := e_i-e_{i+1}$, is a set of simple roots. The
corresponding set of positive roots is 
$\Phi^+ = \{ e_i-e_j \, | \, 1 \le i < j \le n+1 \}$. 

Let $\alpha = e_i-e_j \in \Phi^+$. Obviously, the hook corresponding to
$\alpha$ is
\begin{equation} \label{eq:hook_An}
h(\alpha) = \{\alpha\} \cup \{e_i-e_s, e_s-e_j \, | \,  i < s < j \}.
\end{equation}
The next lemma describes the closed patterns $n(\alpha)$ and
$k(\alpha)$ for all positive roots~$\alpha$. 

\begin{lem} \label{la:An_n_k}
Let $\Phi$ be a root system of type $A_n$ as described above. For all
positive roots $\alpha = e_i-e_j$ the following are true:
\begin{enumerate}
\item[(a)] $n(\alpha), k(\alpha) \unlhd \Phi^+$,

\item[(b)] $n(\alpha) = \{e_s-e_t \, | \, 1 \le s \le i \,\,
  \text{and} \,\, j \le t \le n+1\} \setminus \{\alpha\}$,

\item[(c)] $k(\alpha) = \{e_s-e_t \, | \, 1 \le s < i \,\,
  \text{or} \,\, j < t \le n+1\}$.
\end{enumerate}
\end{lem}

\pr
Part (a) was already shown in Section~\ref{sec:representable}.

\medskip

\noindent (b) By definition $n_0(\alpha) = \{e_i-e_j\}$.
Let $\gamma = e_k-e_l \in \Phi^+$. Then
$\alpha + \gamma  \in \Phi^+$ if and only if $l=i$ or $k=j$. 
In the first case we have $\alpha + \gamma = e_k-e_j$ and in the
second case we have $\alpha + \gamma = e_i-e_l$. Thus,
\[
n_1(\alpha) = \{e_k-e_j \, | \, 1 \le k < i\} \cup
\{e_i-e_l \, | \, j < l \le n+1\}.
\]
Now let $\gamma = e_s-e_t \in \Phi^+$ and $k < i$. Then 
$\gamma + e_k-e_j \in \Phi^+$ if and only if $t=k$ or $s=j$. In the
first case we have $\gamma + e_k-e_j = e_s-e_j \in n_1(\alpha)$ and in
the second case we have $\gamma + e_k-e_j = e_k-e_t$ where $k<i$ and
$t>j$. Similarly, $\gamma + e_i-e_l \in \Phi^+$ if and only if $t=i$
or $s=l$. In the first case we have $\gamma + e_i-e_l = e_s-e_l$ where
$s<i$ and $l>j$ and in the second case we have 
$\gamma + e_i-e_l = e_i-e_t \in n_1(\alpha)$. It follows that
\[
n_2(\alpha) \setminus n_1(\alpha) = \{e_s-e_t \, | \, 1 \le s < i \,\,
\text{and} \,\, j < t \le n+1\}.
\]
Finally, we see that for all $\gamma \in \Phi^+$ and 
$\beta \in n_2(\alpha) \setminus n_1(\alpha)$ we have 
$\beta+\gamma \in n_2(\alpha)$ and therefore
$n(\alpha) = n_1(\alpha) \cup n_2(\alpha)$ and (b) follows. 

\medskip

(c) Let $\beta = e_s-e_t \in \Phi^+$ where $i \le s < t \le j$. Then
we have $\beta \in h(e_i-e_t)$ and $e_i-e_t \in h(\alpha)$ and hence
$\beta \in w(\alpha)$. Hence $\{ e_s-e_t \, | \, i \le s < t \le j \}
\subseteq w(\alpha)$ and Lemma~\ref{la:k_w} implies that 
$k(\alpha) \subseteq \{e_s-e_t \, | \, 1 \le s < i \,\,
  \text{or} \,\, j < t \le n+1\} =: M$.

We claim that $M \unlhd \Phi^+$. Let $\beta = e_s-e_t \in M$,
$\gamma = e_k-e_l \in \Phi^+$ with $\beta+\gamma \in \Phi^+$. By the
definition of $M$ we have $s<i$ or $t>j$. Suppose that $s<i$. We have
$t=k$ or $s=l$. If $t=k$ then $\beta+\gamma=e_s-e_l$ with $s<i$ and
hence $\beta+\gamma \in M$. If $s=l$ then $\beta+\gamma=e_k-e_t$
with $k<l=s<i$ and again $\beta+\gamma \in M$. Suppose that $t>j$. 
Again we have $t=k$ or $s=l$. If $t=k$ then $\beta+\gamma=e_s-e_l$ 
with $l>k=t>j$ and hence $\beta+\gamma \in M$. If $s=l$ then 
$\beta+\gamma=e_k-e_t$ with $t>j$ and again $\beta+\gamma \in M$.
Thus we have $M \unlhd \Phi^+$. By the definition of $M$ we have
$\alpha \not\in M$. Hence Remark~\ref{rem:kSigma}~(a) implies that 
$M \subseteq k(\alpha)$. Thus $k(\alpha)=M$ and (c) follows.
\epr

\medskip

In the following we study the hooks for type $A_n$ more closely. We
define the arm $a(\alpha)$ and the leg $\ell(\alpha)$ of $h(\alpha)$
as follows: 
\[
   a(\alpha) := \{ e_i - e_s \, | \, i< s < j \} \quad \text{and} \quad 
\ell(\alpha) := \{ e_s - e_j \, | \, i < s < j \}.
\]
Recall that a $p$-group $P$ is \emph{special} if 
$\Phi(P) = [P,P] = Z(P)$ is elementary abelian or if $P$ itself is
elementary abelian. As in \cite{HML_D4} we say that~$P$ is 
\emph{special of type $q^{1+2a}$} if $P$ is special with $|P| =
q^{2a+1}$ and $|Z(P)|=q$.

We will show below that for each $\alpha \in \Phi^+$ the
hook $h(\alpha)$ is a closed pattern. We call the pattern subgroup 
$H_\alpha := P(h(\alpha))$ the \emph{hook subgroup} corresponding
to~$\alpha$. 

\begin{lem} \label{la:An_hooks}
Let $\Phi$ be a root system of type $A_n$ as described at the
beginning of this section. For all $\alpha = e_i-e_j \in \Phi^+$
the following are true:
\begin{enumerate}
\item[(a)] The hook $h(\alpha)$, the arm $a(\alpha)$ and the leg
  $\ell(\alpha)$ are closed patterns. 

\item[(b)] The canonical projection $\pi: UA_n \rightarrow
  UA_n/P(k(\alpha))$ maps the hook subgroup $H_\alpha=P(h(\alpha))$
  injectively into $UA_n/P(k(\alpha))$ and $\pi(H_\alpha)$ is normal
  in $UA_n/P(k(\alpha))$. 

\item[(c)] The pattern subgroups $P(a(\alpha))$ and $P(\ell(\alpha))$
  are elementary abelian.

\item[(d)] If $\height(\alpha) > 1$ then the hook subgroup $H_\alpha$ is
  special of type $q^{1+2(j-i-1)}$ and $[y, H_\alpha]=Z(H_\alpha)=X_\alpha$ 
  for all $y \in H_\alpha \setminus Z(H_\alpha)$. More specifically:
  For each $y \in H_\alpha \setminus Z(H_\alpha)$ there is some
  $\beta \in h(\alpha)$ such that 
  $\{[y,x_\beta(t)] \, | \, t \in \F_q\}=X_\alpha$.
\end{enumerate}
\end{lem}
\pr
(a), (c): Let $\beta, \gamma \in h(\alpha)$. Then we have
$\beta+\gamma \in \Phi^+$ if and only if 
$\{\beta, \gamma\} = \{e_i-e_s, e_s-e_j\}$ for some $i<s<j$ and in
this case $\beta+\gamma = \alpha \in h(\alpha)$. In particular, we
have $\beta+\gamma \not\in \Phi^+$ if $\beta, \gamma \in a(\alpha)$ or
if $\beta, \gamma \in \ell(\alpha)$. This implies (a) and (c).

\medskip

\noindent (b) By Lemma~\ref{la:k_w} we have 
$h(\alpha) \cap k(\alpha) \subseteq w(\alpha) \cap k(\alpha) =
\emptyset$. Since $H_\alpha$, $P(k(\alpha))$ are pattern subgroups we
get that the restriction of $\pi$ to $H_\alpha$ is injective. 

Let $\beta \in h(\alpha)$ and $\gamma=e_k-e_l \in \Phi^+$ so that
$k<l$. Suppose first that $\beta = e_i-e_s$ for some $i<s\le j$. Then
$\beta + \gamma \in \Phi^+$ if only if $s=k$ or $l=i$. In the first
case we have $\beta+\gamma=e_i-e_l \in h(\alpha) \cup k(\alpha)$. In
the second case we have $\beta+\gamma=e_k-e_s \in k(\alpha)$ because
$k<i$. Now suppose that $\beta = e_s-e_j$ for some $i \le s < j$. Then
$\beta + \gamma \in \Phi^+$ if only if $s=l$ or $k=j$. In the first
case we have $\beta+\gamma=e_k-e_j \in h(\alpha) \cup k(\alpha)$. In
the second case we have $\beta+\gamma=e_s-e_l \in k(\alpha)$ because
$l>j$. This proves (b).

\medskip

\noindent (d) Suppose that $\height(\alpha)>1$. We have seen in the
proof of (a) and (c) that for all $\beta, \gamma \in h(\alpha)$
we have $\beta+\gamma \in \Phi^+$ if and only if 
$\{\beta, \gamma\} = \{e_i-e_s, e_s-e_j\}$ for some $i<s<j$ and in
this case $\beta+\gamma = \alpha \in h(\alpha)$. It follows from
Lemmas~\ref{la:derquopat} and \ref{la:cenquopat} that
$Z(H_\alpha) = [H_\alpha, H_\alpha] = X_\alpha$. So $H_\alpha$ is
special. Now let $y \in H_\alpha \setminus Z(H_\alpha)$. Note that
$H_\alpha \setminus Z(H_\alpha) \neq \emptyset$ since $\height(\alpha)>1$.
We write $y=\prod_{\gamma \in h(\alpha)} x_\gamma(t_\gamma)$ as in
(\ref{eq:elemu}). Because $y \not\in X_\alpha$ there is some 
$\gamma \in h(\alpha) \setminus \{\alpha\}$ such that $t_\gamma \neq 0$.
Thus, $\beta := \alpha-\gamma \in h(\alpha)$ and we get from
Lemma~\ref{la:patnorm}~(a) that $\{[y,x_\beta(t)] \, | \, t \in \F_q\} = X_\alpha$.
\epr

\medskip

The following lemma prepares a reduction result for $\Irr(UA_n)_\alpha$, 
the set of single root characters of $UA_n$ corresponding to positive roots 
$\alpha$. 

\begin{lem} \label{la:An_source}
Let $\Phi$ be a root system of type $A_n$ as described at the
beginning of this section. For all $\alpha = e_i-e_j \in \Phi^+$
the following are true:
\begin{enumerate}
\item[(a)] The source $s(\alpha)$ is a closed pattern.

\item[(b)] $\ell(\alpha) \cup k(\alpha) \unlhd s(\alpha)$.

\item[(c)] $\ell(\alpha) \cup \{\alpha\} \cup k(\alpha) \unlhd \Phi^+$.
\end{enumerate}
\end{lem}
\pr
(a) Let $\beta = e_k-e_l, \gamma = e_{k'}-e_{l'} \in s(\alpha)$.
Suppose that $\beta+\gamma \in a(\alpha)$. Then $\beta+\gamma=e_i-e_s$ 
for some $i<s<j$. Thus $k=i$ or $k'=i$. We can assume that $k=i$ and
then $k'=l$ and $l'=s$. But then $\beta=e_i-e_l$ where $l=k'<l'=s<j$
and hence $\beta \in a(\alpha)$, a contradiction. Thus 
$\beta+\gamma \not\in \Phi^+$ or $\beta+\gamma \in s(\alpha)$.

\medskip

\noindent (b) Let $\beta = e_s-e_j \in \ell(\alpha)$ where $i<s<j$ and
let $\gamma = e_k-e_l \in s(\alpha)$. Then $\beta + \gamma \in \Phi^+$
if only if $s=l$ or $k=j$. In the first case we have 
$\beta+\gamma=e_k-e_j \in \ell(\alpha) \cup k(\alpha)$ because 
$k \neq i$. In the second case we have $\beta+\gamma=e_s-e_l \in k(\alpha)$ 
because $l>j$. This proves (b). 

\medskip

\noindent (c) We know from Definition~\ref{def:kSigma_almfaith} (a) and
Lemma~\ref{la:alpha_rep} that $\{\alpha\} \cup k(\alpha) \unlhd \Phi^+$
and we have seen in (b) that $s(\alpha)$ normalizes 
$\ell(\alpha) \cup k(\alpha)$. By Lemmas~\ref{la:derquopat} and
\ref{la:An_hooks}~(d) the arm $a(\alpha)$ normalizes 
$\ell(\alpha) \cup \{\alpha\}$. Since $\Phi^+ = s(\alpha) \cup a(\alpha)$ 
the claim in (c) follows.

\epr

\medskip

By Lemma~\ref{la:An_source} (a) and (b) we can consider the source
group $S_\alpha = P(s(\alpha))$ and its quotient pattern group
$S_\alpha/P(\ell(\alpha) \cup k(\alpha))$. It follows 
from Lemma~\ref{la:derquopat} that we can identify the root
subgroup~$X_\alpha$ with the quotient pattern group 
$S_\alpha/P(s(\alpha) \setminus \{\alpha\})$. 
Lemmas~\ref{la:derquopat} and \ref{la:representable} imply that
\[
S_\alpha/P(\ell(\alpha) \cup k(\alpha)) \cong \overline{T}_\alpha \times X_\alpha
\] 
where 
\[
\overline{T}_\alpha := P(s(\alpha) \setminus \{\alpha\})/P(\ell(\alpha)
\cup k(\alpha)).
\]
The group $\overline{T}_\alpha$ is generated by the images of
$X_\gamma$ in $UA_n/P(\ell(\alpha) \cup k(\alpha))$ where $\gamma$ is
an element of
\begin{eqnarray*}
(s(\alpha) \setminus \{\alpha\}) \setminus (\ell(\alpha) \cup k(\alpha)) 
& = & \Phi^+ \setminus (a(\alpha) \cup \{\alpha\} \cup \ell(\alpha) \cup k(\alpha))\\
& = & \Phi^+ \setminus (h(\alpha) \cup k(\alpha))
= \{e_s-e_t \mid i < s < t < j\}.
\end{eqnarray*}
The set $\{e_s-e_t \mid i < s < t < j \} \subseteq \Phi^+$ is a
closed pattern, and if $j-i-2>0$ it generates a root subsystem of 
type $A_{j-i-2}$ of $\Phi$. Thus 
\[
\overline{T}_\alpha \cong P(\{e_s-e_t \mid i < s < t < j \}) \cong UA_{j-i-2}.
\]
We can now give a recursive description for the single root 
characters of~$UA_n$:

\begin{prop} \label{prop:An_singlechars}
Let $\Phi$ be a root system of type $A_n$ as described above. For
each root $\alpha \in \Phi^+$ the map
$\Psi_\alpha: \Irr(\overline{T}_\alpha) \times \Irr(X_\alpha)^*
\to \Irr(UA_n)_\alpha$ with 
\[
(\mu, \lambda) \mapsto (\infl_{\overline{T}_\alpha}^{S_\alpha} \mu \cdot
 \infl_{X_\alpha}^{S_\alpha} \lambda)^{UA_n}
\]
is a one to one correspondence with the property
$\Psi_\alpha(\mu, \lambda)(1) = q^{j-i-1} \cdot \mu(1)$. 
We have $\overline{T}_\alpha \cong UA_{j-i-2}$ if $j-i-2>0$ and 
$\overline{T}_\alpha = \{1\}$ otherwise.
\end{prop}
\pr
The content of Lemma \ref{la:An_source} is that 
is that hypothesis (1), (2) and (3) of Proposition \ref{prop:red_to_thm}
are satisfied, and the content of Lemma~\ref{la:An_hooks} (d) is that 
hypothesis (4) of Proposition \ref{prop:red_to_thm}
is satisfied. This proves the correspondence. 
\epr


\subsection{Type B} \label{subsec:Bn}

Let $n \ge 2$ be an integer. We construct a root system of
type~$B_n$ as in~\cite[Section~12.1]{Humphreys}: Let 
$e_1, e_2, \dots, e_n \in \R^n$ be the usual orthonormal unit
vectors which form a basis of $\R^n$. Then
$\Phi := \{ \pm(e_i \pm e_j) \, | \, 1 \le i \neq j \le n \} \cup \{ \pm e_i \, | \, 1 \le i \le n\}$ is a
root system of type~$B_n$ and the set $\{\alpha_1, \dots, \alpha_n\}$,
where $\alpha_i := e_i-e_{i+1}$ for $i=1,2,\dots,n-1$ and 
$\alpha_n := e_n$, is a set of simple roots. The
corresponding set of positive roots is 
$\Phi^+ = \{ e_i \pm e_j \, | \, 1 \le i < j \le n \} \cup \{ e_i \, | \, 1 \le i \le n \}$. 

The highest long root with respect to this base is $e_1+e_2$ whereas
the highest short root is $e_1$. We note that the long roots of $\Phi$
form a $D_n$-subsystem. Also we recall that the Chevalley commutator
relations imply that if $\a$ is short and $\b$ is long with 
$\alpha+\beta \in \Phi^+$, then $[X_\a,X_\b] \subseteq X_{\a+\b}X_{2\a+b}$, 
where $\a+\b$ is short and $2\a + \b$ is long, and that when both 
$\a$ and $\b$ are short with $\alpha+\beta \in \Phi^+$, then 
$[X_\a,X_\b] \subseteq X_{\a +\b}$, where $\a+\b$ is long.

From the explicit description of the root systems (or from the Dynkin 
diagram) we see that $\{\alpha_2, \alpha_3, \dots, \alpha_n\}$
generates a root subsystem~$\Phi_{1'}$ of type $B_{n-1}$ and that 
$\{\alpha_1, \alpha_2, \dots, \alpha_{n-1}\}$ generates 
a root sub\-system~$\Phi_{n'}$ of type~$A_{n-1}$. We set 
$\Phi^+_{1'} := \Phi_{1'} \cap \Phi^+$, $\Psi^+_{1'} := \Phi^+
\setminus \Phi^+_{1'}$ and $\Phi^+_{n'} := \Phi_{n'} \cap \Phi^+$,
$\Psi^+_{n'} := \Phi^+ \setminus \Phi^+_{n'}$. The sets 
$\Psi^+_{1'}$ and $\Psi^+_{n'}$ are normal closed patterns. In fact,
$P(\Psi^+_{1'})$ and $P(\Psi^+_{n'})$ are the unipotent radicals of
the standard parabolic subgroups corresponding to 
$\{\alpha_2, \alpha_3, \dots, \alpha_n\}$ and
$\{\alpha_1, \alpha_2, \dots, \alpha_{n-1}\}$, respectively.
Furthermore, the sets $\Phi^+_{1'}$ and $\Phi^+_{n'}$ are closed
patterns and $P(\Phi^+_{1'}) \cong UB_n/P(\Psi^+_{1'}) \cong UB_{n-1}$
and $P(\Phi^+_{n'}) \cong UB_n/P(\Psi^+_{n'}) \cong UA_{n-1}$.

Let $\alpha \in \Phi^+$ and $\chi \in \Irr(UB_n)_\alpha$. Suppose that  
$\alpha \in \Phi^+_{1'}$. It follows from Lemma~\ref{la:alpha_rep} and
Remark~\ref{rem:kSigma} (a), (c) that $P(\Psi^+_{1'}) \subseteq \Ker(\chi)$. 
Thus, we can identify $\chi$ with a single root character of 
$P(\Phi^+_{1'}) \cong UB_n/P(\Psi^+_{1'}) \cong UB_{n-1}$. In this
way, the classification and construction of the elements of
$\Irr(UB_n)_\alpha$ are reduced to the case $B_{n-1}$. Similarly, if
$\alpha \in \Phi^+_{n'}$ we can identify $\chi$ with a single root
character of $P(\Phi^+_{n'})$ and thereby get a reduction to the case
$A_{n-1}$ which has already been treated in Section~\ref{subsec:An}. 
Hence, we only have to consider positive roots $\alpha$ which are not
contained in $\Phi^+_{1'} \cup \Phi^+_{n'}$, i.e., the roots
$e_1 + e_i$ where $1 < i \le n$ and the root $e_1$.

\begin{lem} \label{la:hook_Bn} Let $\Phi$ be a root system of type 
$B_n$ as described above. For all positive roots of the form $\alpha = e_1+e_i$ 
we have:
\[
h(\alpha) = \{e_1+e_i, e_1, e_i \} \cup \{e_1 -e_s, e_s+e_i \mid 1 < s < i \} 
\cup \{e_1 \pm e_s, e_i \mp e_s \mid i < s \leq n \}.
\]
\end{lem}

\pr
First we observe that clearly $e_1, e_i \in h(e_1+e_i)$, and that no
other combination of short roots can add to $\alpha.$

Next let $\beta=e_s \pm e_t$, $\gamma = e_{s'} \in \Phi^+$. We see
immediately that $\beta+\gamma \not \in \Phi^+$ or $\beta+\gamma$ is a
short root and therefore $\beta+\gamma \neq \alpha$. Hence, 
$\beta, \gamma \not\in h(\alpha)$.

Finally let $\beta=e_s \pm e_t$, $\gamma = e_{s'} \pm e_{t'} \in \Phi^+$ 
with $s<t$ and $s'<t'$ and $\beta+\gamma=\alpha$. We see immediately
that one of the two $\pm$ signs has to be a $+$~sign and that the
other one has to be a $-$~sign and that $t=t'$. Furthermore, we can
assume $s=1$ and $s'=i$. Hence 
\[
(\beta, \gamma) \in \{(e_1-e_s, e_s+e_i) \, | \, 1 < s
< i \,\, \text{or} \,\, i < s \le n\} \cup \{(e_1+e_s, e_i-e_s) \, | \, i
< s \le n\}
\]
and the claim follows.
\epr

\medskip

We are now able to determine the patterns $n(\alpha)$ and $k(\alpha)$.

\begin{lem} \label{la:Bn_n_k}
Let $\Phi$ be a root system of type $B_n$ as described above. For all
positive roots of the form $\alpha = e_1+e_i$ the following are true:
\begin{enumerate}
\item[(a)] $n(\alpha), k(\alpha) \unlhd \Phi^+$. 

\item[(b)] $n(\alpha) = \{ e_1+e_s \, | \, 1 < s < i\}$.

\item[(c)] $k(\alpha) = \{ e_s+e_t \, | \, 1 \le s < t < i\}$. 
\end{enumerate}
\end{lem}

\pr
Part (a) was already shown in Section~\ref{sec:representable}.

\medskip

\noindent (b) By definition $n_0(\alpha) = \{e_1+e_i\}$.
Let $\gamma \in \Phi^+$. For $\gamma = e_s$ or $\gamma = e_k+e_l$ we
have $\alpha + \gamma  \not \in \Phi^+$. 

Therefore, $\alpha + \gamma  \in \Phi^+$ if and only if $\gamma = e_s-e_i$ 
for some $1 < s <i$ and in this case we have $\alpha + \gamma = e_1+e_s$.
Thus, $n_1(\alpha) = \{e_1+e_s \, | \, 1 < s < i\}$. Again, if 
$\gamma \in \Phi^+$ is of the form $\gamma = e_k$ then 
$\gamma + (e_1 +e_s) \not \in \Phi^+$. But if 
$\gamma = e_k \pm e_l \in \Phi^+$ then $\gamma + (e_1+e_s)  \in \Phi^+$
if and only if $\gamma = e_t-e_s$ for some $1 < t <s$ and in this case
we have $\gamma + (e_1+e_s) = e_1+e_t$. Thus $n_2(\alpha) \subseteq n_1(\alpha)$ 
and then $n(\alpha) = n_1(\alpha)$ and (b) follows. 

\medskip

\noindent (c) Let $\beta = e_s$. If $s = 1$, then 
$\beta \in h(\alpha) \subseteq w(\alpha)$. If $s > 1$, then 
$\beta + (e_1 -e_s) = e_1 \in h(\alpha)$ implies that 
$\beta \in w(\alpha)$.

Let $\beta = e_s + e_t \in \Phi^+$ where $s<t$ and $t \ge i$. 
If $s=1$ then $\beta \in h(\alpha) \subseteq w(\alpha)$.
If $s>1$ then $\beta \in h(e_1+e_t)$ and $e_1+e_t \in h(\alpha)$,
hence $\beta \in w(\alpha)$.  

Now let $\beta = e_s - e_t \in \Phi^+$ where $s<t$. If $s=1$, 
then $\beta \in h(e_1)$ and $e_1 \in h(\alpha)$. Thus, 
$\beta \in w(\alpha)$. If $s>1$, then $\beta \in h(e_1-e_t)$ and 
$e_1-e_t \in h(e_1)$ and $e_1 \in h(\alpha)$. Hence 
$\beta \in w(\alpha)$. 

It follows from Lemma~\ref{la:k_w} that 
\[
k(\alpha) \subseteq \{ e_s+e_t \, | \, 1 \le s < t < i\} =: M.
\]
We claim that $M \unlhd \Phi^+$. Let $\beta = e_s+e_t \in M$,
$\gamma \in \Phi^+$ with $\beta+\gamma \in \Phi^+$. It follows that
$\gamma$ is of the form $\gamma = e_k - e_s$ with $k<s<t<i$ or 
$\gamma = e_k - e_t$ with $k<t<i$. In the first case we have
$\beta+\gamma = e_k + e_t$ with $k<t<i$, in the second case we have
$\beta+\gamma = e_k + e_s$ with $k,t<i$. Thus, $\beta+\gamma \in M$ in
both cases. Thus we have $M \unlhd \Phi^+$. By the definition of $M$
we have $\alpha \not\in M$. Hence Remark~\ref{rem:kSigma}~(a) implies
that $M \subseteq k(\alpha)$. Thus $k(\alpha)=M$ and (c) follows.
This completes the proof.
\epr

\medskip

For $\alpha = e_1+e_i \in \Phi^+$  we define the arm $a(\alpha)$ and
the leg $\ell(\alpha)$ of $h(\alpha)$ as follows:
\begin{eqnarray*}
a(\alpha) & := & \{ e_i \} \cup \{ e_i \pm e_s \mid i < s \} \cup \{e_1-e_s \mid 1 < s < i \} 
\quad \text{and}\\
\ell(\alpha) & := & \{e_1\} \cup \{ e_1 \pm e_s  \mid i < s\} \cup \{e_s+e_i \mid 1 < s < i \}. 
\end{eqnarray*}
As for type $A_n$, we will show below that for each $\alpha \in \Phi^+$
the hook $h(\alpha)$ is a closed pattern and call $H_\alpha := P(h(\alpha))$
the \emph{hook subgroup} corresponding to~$\alpha$. 

\begin{lem} \label{la:Bn_hooks}
Let $\Phi$ be a root system of type $B_n$ as described at the
beginning of this section. For all $\alpha = e_1+e_i \in \Phi^+$
the following are true:
\begin{enumerate}
\item[(a)] The hook $h(\alpha)$, the arm $a(\alpha)$ and the leg
  $\ell(\alpha)$ are closed patterns. 

\item[(b)] The canonical projection $\pi: UB_n \rightarrow
  UB_n/P(k(\alpha))$ maps the hook subgroup $H_\alpha=P(h(\alpha))$
  injectively into $UB_n/P(k(\alpha))$.

\item[(c)] The pattern subgroups $P(a(\alpha))$ and $P(\ell(\alpha))$
  are elementary abelian.

\item[(d)] If $\height(\alpha) > 1$ then the hook subgroup $H_\alpha$ is
  special of type $q^{1+2(2n-i-1)}$ and $[y, H_\alpha]=Z(H_\alpha)=X_\alpha$ 
  for all $y \in H_\alpha \setminus Z(H_\alpha)$. More specifically:
  For each $y \in H_\alpha \setminus Z(H_\alpha)$ there is some
  $\beta \in h(\alpha)$ such that 
  $\{[y,x_\beta(t)] \, | \, t \in \F_q\}=X_\alpha$.
\end{enumerate}
\end{lem}
\pr
(a), (c): Let $\beta, \gamma \in h(\alpha)$. Then we have
$\beta+\gamma \in \Phi^+$ if and only if 
$\{\beta, \gamma\} = \{e_1-e_s, e_s+e_i\}$ for some $s \neq i$, or
$\{\beta, \gamma\} = \{e_1+e_s, e_i-e_s \}$ for some $s>i$, or 
$\{\beta, \gamma\} = \{e_1, e_i \}$ and in all
of these cases $\beta+\gamma = \alpha \in h(\alpha)$. In particular,
we have $\beta+\gamma \not\in \Phi^+$ if $\beta, \gamma \in a(\alpha)$
or if $\beta, \gamma \in \ell(\alpha)$. This implies (a) and (c).

\medskip

\noindent (b) By Lemma~\ref{la:k_w} we have 
$h(\alpha) \cap k(\alpha) \subseteq w(\alpha) \cap k(\alpha) =
\emptyset$. Since $H_\alpha$, $P(k(\alpha))$ are pattern subgroups we
get that the restriction of $\pi$ to $H_\alpha$ is injective. 

\medskip 

\noindent (d) Suppose that $\height(\alpha)>1$. We have seen in the
proof of (a) and (b) that for all $\beta, \gamma \in h(\alpha)$
we have $\beta+\gamma \in \Phi^+$ if and only if 
$\{\beta, \gamma\} = \{e_1-e_s, e_s+e_i\}$ for some $s \neq i$ or
$\{\beta, \gamma\} = \{e_1+e_s, e_i-e_s \}$ for some $s>i$ or 
$\{\beta, \gamma\} = \{e_1, e_i \}$  and in all of these cases 
$\beta+\gamma = \alpha \in h(\alpha)$.  It follows from
Lemmas~\ref{la:derquopat} and \ref{la:cenquopat} that
$Z(H_\alpha) = [H_\alpha, H_\alpha] = X_\alpha$. So $H_\alpha$ is
special. The remaining part of the proof of (d) is analogous to the
proof of Lemma~\ref{la:An_hooks} (d). 
\epr

\medskip

Note that, in contrast to type $A_n$, we did not claim in
Lemma~\ref{la:Bn_hooks} that $\pi(H_\alpha)$ is normal in
$UB_n/P(k(\alpha))$. In fact, for all $2<i<n$ the patterns $h(\alpha)$ 
and $h(\alpha) \cup k(\alpha)$ are \emph{not} normal because 
$(e_2-e_i) + (e_i-e_n) \not\in h(\alpha) \cup k(\alpha)$. 

\begin{lem} \label{la:Bn_source}
Let $\Phi$ be a root system of type $B_n$ as described at the
beginning of this section. For all $\alpha = e_1+e_i \in \Phi^+$
the following are true:
\begin{enumerate}
\item[(a)] The source $s(\alpha)$ is a closed pattern.

\item[(b)] $\ell(\alpha) \cup k(\alpha) \unlhd s(\alpha)$.

\item[(c)] $\ell(\alpha) \cup \{\alpha\} \cup k(\alpha) \unlhd \Phi^+$.
\end{enumerate}
\end{lem}
\pr 
(a) Let $\beta \in a(\alpha)$ and $\gamma, \gamma' \in \Phi^+$ such
that $\beta=\gamma+\gamma'$. To prove that $s(\alpha)$ is a pattern it
suffices to show that $\gamma \in a(\alpha)$ or $\gamma' \in a(\alpha)$.

Suppose that $\beta=e_i$. Then $\{\gamma, \gamma'\} = \{e_i-e_s, e_s\}$ for some $i < s$ and
hence $\gamma \in a(\alpha)$ or $\gamma' \in a(\alpha)$. Suppose next that 
$\beta=e_1-e_s$ where $1 < s < i$. Then 
$\{\gamma, \gamma'\} = \{e_1-e_l, e_l-e_s\}$ for some $1<l<s<i$ and
hence $\gamma \in a(\alpha)$ or $\gamma' \in a(\alpha)$. Now suppose
that $\beta = e_i+e_s$ where $i < s \le n$. Then 
$\{\gamma, \gamma'\} = \{e_i+e_l, e_s-e_l\}$ for some $i < l \le n$ or 
$\{\gamma, \gamma'\} = \{e_s+e_l, e_i-e_l\}$ for some $i < l \le n$ or 
$\{\gamma, \gamma'\} = \{e_s, e_i\}$.
Hence in all three cases $\gamma \in a(\alpha)$ or $\gamma' \in a(\alpha)$.
Finally, suppose that $\beta = e_i-e_s$ where $i < s \le n$. Then 
$\{\gamma, \gamma'\} = \{e_i-e_l, e_l-e_s\}$ for some $i < l < s \le n$ 
and again $\gamma \in a(\alpha)$ or $\gamma' \in a(\alpha)$.

\medskip 

\noindent (b) Let $\beta \in \ell(\alpha)$ and $\gamma \in s(\alpha)$ 
such that $\beta+\gamma \in \Phi^+$. We have to show that 
$\beta + \gamma \in \ell(\alpha) \cup k(\alpha)$. If
$\beta = e_s+e_i$ where $1<s<i$ then $\gamma$ has to be of the form 
$\gamma=e_k-e_l$ where $k < l$. Thus $\beta+\gamma=e_k+e_i$ where
$1<k<i$ or $\beta+\gamma=e_s+e_k$ where $s,k<i$. In the first case 
$\beta+\gamma \in \ell(\alpha)$ and in the second case 
$\beta+\gamma \in k(\alpha)$. Now suppose that $\beta=e_1+e_s$ where
$s>i$. Then $\gamma$ has to be of the form $\gamma=e_k-e_l$ where 
$k < l$ and $\beta+\gamma=e_1+e_k$ where $k\neq1, i$. Hence
$\beta+\gamma \in k(\alpha) \cup \ell(\alpha)$.

Now suppose that
$\beta=e_1-e_s$ where $s>i$ and $\gamma=e_k+e_l$ where $k<l$ or $\gamma = e_l$. Then
$\beta+\gamma=e_1+e_k$ where $k\neq1,i$ or $\beta+\gamma=e_1+e_l$
where $l\neq1,i$ or $\beta + \gamma = e_1$. Hence
$\beta+\gamma \in k(\alpha) \cup \ell(\alpha)$. Finally, suppose
that $\beta=e_1-e_s$ where $s>i$ and $\gamma=e_k-e_l$ where
$k<l$. Then $\beta+\gamma=e_1-e_l$ where $l>i$. Thus 
$\beta+\gamma \in \ell(\alpha)$, proving (b).

\bigskip

\noindent (c) In light (b) and Proposition \ref{prop:red_to_thm} it suffices to 
prove that $\ell(\alpha) \cup \{\alpha \} \cup k(\alpha) \unlhd h(\alpha) \cup k(\alpha).$
Let $\beta \in \ell(\alpha)$. As $P(\ell(\alpha))$ is abelian we may assume without loss 
that $\gamma \in a(\alpha)$. 

If $\gamma = e_i$ or $e_i \pm e_s$ with $s > i$, then 
$\beta  \neq e_r + e_i$ with $r < i$, but if $\beta = e_1 \pm e_r$ with $s > i$ then 
$\gamma + \beta = e_1 +  e_i \pm e_r \pm e_s $. Thus 
$\gamma + \beta \in \Phi^+$ only if $r =s$ and with opposite signs; i.e., 
$\gamma + \beta = \alpha.$

If $\gamma = e_1 -e_r$ with $1 < r < i$, then $\beta \neq e_1 \pm e_s $ with $s > i$.
But if $\beta = e_s + e_i$, then $\gamma + \beta \in \Phi^+$ only if $r =s$ and in this 
case $\gamma + \beta = \alpha.$ This proves that 
$\ell(\alpha) \cup \{\alpha \} \cup k(\alpha) \unlhd h(\alpha) \cup k(\alpha),$ which implies (c).
\bigskip
\epr

\begin{rem}
Note that, in contrast to type $A_n$, we did not claim in
Lemma~\ref{la:Bn_hooks} that $\pi(H_\alpha)$ is normal in
$UB_n/P(k(\alpha))$. In fact, for all $2<i<n$ the patterns $h(\alpha)$ 
and $h(\alpha) \cup k(\alpha)$ are \emph{not} normal because 
$(e_2-e_i) + (e_i-e_n) \not\in h(\alpha) \cup k(\alpha)$. 
\end{rem}

We can now establish our desired correspondence for long root subgroups of root systems of 
type $B_n.$

\begin{prop} \label{prop:Bn_singlechars}
Let $\Phi$ be a root system of type $B_n$ as described above. For
each long root $\alpha \in \Phi^+$ the map
$\Psi_\alpha: \Irr(\overline{T}_\alpha) \times \Irr(X_\alpha)^*
\to \Irr(UB_n)_\alpha$ with 
\[
(\mu, \lambda) \mapsto (\infl_{\overline{T}_\alpha}^{S_\alpha} \mu \cdot
 \infl_{X_\alpha}^{S_\alpha} \lambda)^{UB_n}
\]
is a one to one correspondence with the property
$\Psi_\alpha(\mu, \lambda)(1) = q^{2n-i-1} \cdot \mu(1)$. 
\end{prop}
\pr
The content of Lemma \ref{la:Bn_source} is that 
is that hypothesis (1), (2) and (3) of Proposition \ref{prop:red_to_thm}
are satisfied, and the content of Lemma~\ref{la:Bn_hooks} (d) is that 
hypothesis (4) of Proposition \ref{prop:red_to_thm}
is satisfied. This proves the correspondence. 
\epr

Our next result gives structural information concerning the group $\overline{T}_\alpha$. To this end
we define patterns 
$$b_{i,j}:= \{ e_t, e_r \pm e_s \ | \ r,s,t \not \in \{ i,j\}, \ \mbox{and} \ r < s \}$$
where $1 \leq i<j \leq n$ and,   
$$ob(e_1+e_i):= \{ e_1-e_i ,\dots, e_{i-1}-e_i \}.$$ 

\begin{lem} \label{la:Bn_structure_sa}
Let $\Phi$ be a root system of type $B_n$ as described at the
beginning of this subsection. For all $\alpha=e_1+e_i \in \Phi^+$
the following is true:
\begin{enumerate}
\item [(a)] $S_\alpha/P(\ell(\alpha) \cup k(\alpha)) \cong
  \overline{S}_\alpha/P(\ell(\alpha))P(k(\alpha)) \cong \overline{T}_\alpha \times X_\alpha$.

\item[(b)] $\overline{T}_\alpha \cong P(b_{1,i})/P(b_{1,i} \cap k(\alpha))\ltimes
  P(\ob(\alpha))/P(\ob(\alpha) \cap k(\alpha))$. 

\item[(c)] If $i<n$ then $\kob(\alpha) := b_{1,i} \setminus \{e_r-e_s
  \, | \, 1<r<s<i \}$ is a pattern that is normal in $b_{1,i}$ and
  $P(\kob(\alpha))$ centralizes $P(\ob(\alpha))/P(\ob(\alpha) \cap k(\alpha))$.

\item[(d)] If $i \in \{2,3\}$, then $\kob(\alpha) = b_{1,i}$,
  $k(\alpha) = n(\alpha)$ and
  \[
  \overline{T}_\alpha \cong UB_{n-2} \times (UA_1)^{i-1}.
  \]

\item[(e)] If $i>3$ then $P(b_{1,i})/P(b_{1,i} \cap k(\alpha))$ is
  isomorphic to a quotient pattern group of $UB_{n-2}$.

\item[(f)] If $i>3$ then $P(b_{1,i})/P(\kob(\alpha) \cup (k(\alpha)
  \cap b_{1,i})) \cong UA_{i-3}$.
\end{enumerate}
\end{lem}

\pr
(a) follows from the definition of $\overline{T}_\alpha$ and
Lemma~\ref{la:Bn_source} (c) and (d). 

\medskip

\noindent (b) By inspection, 
$b_{1,i} \cap (\{\alpha\} \cup a(\alpha)) = \emptyset$ and
$\ob(\alpha) \cap (\{\alpha\} \cup a(\alpha)) = \emptyset$. Hence,
$b_{1,i} \cup \ob(\alpha) \cup \{\alpha\} \cup \ell(\alpha) \cup
k(\alpha) \subseteq s(\alpha) \setminus \{\alpha\}$. Thus
\begin{equation} \label{eq:semi_bn}
P(b_{1,i} \cup \ob(\alpha))P(\ell(\alpha) \cup
k(\alpha))/P(\ell(\alpha) \cup k(\alpha)) \subseteq \overline{T}_\alpha.
\end{equation}
By definition, $\ob(\alpha) \cup b_{1,i} \cup k(\alpha)
\cup h(\alpha) = \Phi^+$ and therefore
\[
b_{1,i} \cup \ob(\alpha) \cup k(\alpha) \cup \ell(\alpha) \supseteq
\Phi^+ \setminus (\{\alpha\} \cup a(\alpha)) = s(\alpha) \setminus
\{\alpha\},
\]
and it follows that we have equality in (\ref{eq:semi_bn}). For all
$\beta \in \ob(\alpha)$ and $\gamma \in b_{1,i}$ we have 
$\beta+\gamma \not\in \Phi^+$ or $\beta+\gamma \in \ob(\alpha)$ and
hence $P(b_{1,i})$ normalizes $P(\ob(\alpha))$. By definition of
$\ob(\alpha)$ we have $\ob(\alpha) \cap b_{1,i} = \emptyset$ and thus
$\overline{T}_\alpha$ is the semidirect product of $P(b_{1,i}) P(\ell(\alpha)
\cup k(\alpha))/ P(\ell(\alpha) \cup k(\alpha))$ and
$P(\ob(\alpha)) P(\ell(\alpha) \cup k(\alpha))/ P(\ell(\alpha) \cup
k(\alpha))$. Because $b_{1,i} \cap \ell(\alpha) = \emptyset$ and
$\ob(\alpha) \cap \ell(\alpha) = \emptyset$ we have
$P(b_{1,i}) \cap P(\ell(\alpha) \cup k(\alpha)) = P(b_{1,i} \cap
k(\alpha))$ and $P(\ob(\alpha)) \cap P(\ell(\alpha) \cup k(\alpha))
= P(\ob(\alpha) \cap k(\alpha))$ which implies (b).

\medskip

\noindent (c) The only way to write $e_r-e_s$, where $1<r<s<i$, as a
sum of two positive roots is $e_r-e_s = \beta+\gamma$ where
$\{\beta,\gamma\} = \{e_r-e_t, e_t-e_s\}$ for some $r<t<s$. Since
$d_{1,i}$ is a pattern and $\beta, \gamma \not\in \kob(\alpha)$ it
follows that $\kob(\alpha)$ is a pattern which is normal in
$b_{1,i}$. Suppose that $\beta = e_r \pm e_{r'} \in d_{1,i}$ and
$\gamma = e_s - e_i \in \ob(\alpha)$. Then
$\beta+\gamma \in \ob(\alpha)$ only if $\beta=e_r-e_s$ where
$1<r<s<i$, i.e., if $\beta \not\in \kob(\alpha)$. This proves (c). 

\noindent (d) For $i \in \{2,3\}$ we have 
$\{e_r-e_s \, | \, 1<r<s<i \} = \emptyset$. Hence
$\kob(\alpha)=b_{1,i}$ and the semidirect product in (b) is a direct
product. Now Lemma~\ref{la:Bn_n_k}
(b) and (c) imply that $k(\alpha)=n(\alpha)$. Furthermore  
$|\ob(\alpha)| = i-1$ and 
$\ob(\alpha) \cap k(\alpha) = \emptyset$. Since $P(\ob(\alpha))$ is 
elementary abelian we get 
$P(\ob(\alpha))/P(\ob(\alpha) \cap k(\alpha)) \cong (UA_1)^{i-1}$. 
Because $b_{1,i} \cap k(\alpha) = \emptyset$ and $b_{1,i}$ generates a
root subsystem of $\Phi$ of type $B_{n-2}$ we get
$P(b_{1,i})/P(b_{1,i} \cap k(\alpha)) \cong UB_{n-2}$.

\medskip

\noindent (e) follows from the fact that $b_{1,i}$ generates a root
subsystem of $\Phi$ of type $B_{n-2}$.

\medskip

\noindent (f) By definition of $\kob(\alpha)$ the set $b_{1,i}$ is the
disjoint union of $\kob(\alpha)$ and the set
$\{e_r-e_s \, | \, 1<r<s<i \}$. The latter generates a root
subsystem of $\Phi$ of type $A_{i-3}$. This completes the proof.
\epr

\medskip
Finally we consider the highest short root. A simple calculation shows that:

\begin{lem}
$h(e_1) = \{e_1, e_1 - e_i, e_i \ | \ 2 \leq i \leq n  \}$.
\end{lem}

We define 
\begin{eqnarray*} \ell(e_1) & := & \{ e_1 - e_i \ | \ 2 \leq i \leq n \} \,\, \text{   and}\\ 
 a(e_1) & := & \{ e_i \ | \ 2 \leq i \leq n  \},
\end{eqnarray*}  and remark that the leg is a pattern whereas the arm is not a pattern modulo $k(\alpha)$.

\begin{lem}\label{la:b_short}
The following  are true:
\begin{enumerate}
\item[(a)] $n(e_1) = \{ e_1+e_r \, | \, 1 < r \leq n\}$ and 
$ k(e_1) = \{ e_r + e_s \ | \ 1 \leq r < s \leq n \}$. 
\item[(b)] $k(e_1)$ is a normal pattern in $\Phi^+$.
\item[(c)] $\{ e_1 \} \cup \ell(e_1) \cup k(e_1) \unlhd h(e_1) \cup k(e_1)$.
\item[(d)] The group $P(h(e_1) \cup k(e_1))/P(k(e_1))$ is special of order $q^{1+2(n-1)}$. 
Moreover modulo $P(k(e_1))$
$[y, H_{e_1}]=Z(H_{e_1})=X_{e_1}$  for all $y \in H_{e_1}\setminus Z(H_{e_1})$.
\item[(e)] $s(e_1)$ is a closed pattern.
\item[(f)] $\ell(e_1) \unlhd s(e_1)$ and  $s(e_1)/(\{e_1\} \cup \ell(e_1) \cup k(e_1))$ is a pattern of type $A_{n-2}$.
\item[(g)] $\overline{T}_{e_1} \cong UA_{n-2}$, $S_{e_1}/ P(\ell(e_1) \cup k(e_1)) \cong \overline{T}_{e_1} \times X_{e_1}$.
\item[(h)] The map
$\Psi_{e_1}: \Irr(\overline{T}_{e_1}) \times \Irr(X_{e_1})^*
\to \Irr(UB_n)_{e_1}$ with 
\[
(\mu, \lambda) \mapsto (\infl_{\overline{T}_{e_1}}^{S_{e_1}} \mu \cdot
 \infl_{X_{e_1}}^{S_{e_1}} \lambda)^{UB_n}
\]
is a one to one correspondence with the property
$\Psi_{e_1}(\mu, \lambda)(1) = q^{n-1} \cdot \mu(1)$. 
\end{enumerate}

\end{lem}
\pr
\noindent (b) was 
already shown in Section~\ref{sec:representable}.
\medskip

\noindent (a) The fact that $e_1 + e_r$ is the sum of the two short roots $e_1$ and $e_r$ both 
of which lie in $h(e_1)$ shows that $n_0:=\{ e_1 + e_r \ | \ 1 < r \leq n \} \subset n(\alpha)$.
Now if $\gamma + \nu \in \Phi^+$ for $\gamma \in \Phi^+$ and $\nu \in n_0$, 
then $\gamma = e_r -e_s$ and $\gamma + \nu \in n_0$.

We note that $w(e_1) = \Phi^+ \setminus k(e_1)$. Now we already noted that 
the only way to express $e_r$ as a sum of positive roots is as $(e_r-e_s) + e_s$.
Similarly the only way to express $e_1 - e_s$ as a sum of positive roots is as 
$(e_1 - e_r) + (e_r -e_s)$ where $r < s$. Thus we see that 
$w(e_1) = \{ e_i \ | \ 1 \leq i \leq n\} \cup \{e_i - e_j \ | \ 1 \leq i < j \leq n \}$ and thus 
$k(e_1)$ is as claimed. 

\medskip 
\noindent (c) follows from the observation that $e_r + (e_1 - e_s) \in \Phi^+$ only if $r = s$
in which case the sum is equal to $e_1.$ 

\noindent (d) follows  from the 
observation that the sum of any pair of roots from $a(e_1)$ lies in $k(e_1)$ which means that $a(e_1)$ is 
a pattern modulo $k(e_1)$. The proof of the second statement of (d) is as the proof of Lemma~\ref{la:Bn_hooks} (d).  
\medskip 

\noindent To prove Part (e) we note that $s(e_1)$ contains all long roots of $\Phi^+$ and the root $e_1$.
The long roots are the positive roots of a $D_4$ root system and thus from a closed pattern. 
If $\gamma$ is a long root, then $e_1 + \gamma \not \in \Phi^+$. This shows (e).

\medskip 

\noindent
To see the second part of (f) note that 
$\{ e_r -e_s \, | \, 1 < r < s \leq n \} = \Phi^+ \setminus (k(e_1) \cup h(e_1)).$  
To see the normality of the leg note that $(e_1 -e_r) + \beta \in Phi^+$ only if 
$\beta \in \{ e_r- e_s, e_r + e_t \ | \  r < s, \ r \neq t \}$. Hence 
$(e_1 -e_r) + \beta = e_1 - e_s \  \mbox{or} \ e_1 + e_t \in \ell(e_1) \cup k(e_1)$, which is
our claim.

\medskip 

\noindent Now (g)  follows from (f).

\medskip 

\noindent Finally (e), (f), (c) and (d) are the hypotheses of 
Proposition \ref{prop:red_to_thm} and thus (h) follows.  

\epr


\subsection{Type D} \label{subsec:Dn}

Let $n \ge 4$ be an integer. We construct a root system of
type~$D_n$ as in~\cite[Section~12.1]{Humphreys}: Let 
$e_1, e_2, \dots, e_n \in \R^n$ be the usual orthonormal unit
vectors which form a basis of $\R^n$. Then
$\Phi := \{ \pm(e_i \pm e_j) \, | \, 1 \le i \neq j \le n \}$ is a
root system of type~$D_n$ and the set $\{\alpha_1, \dots, \alpha_n\}$,
where $\alpha_i := e_i-e_{i+1}$ for $i=1,2,\dots,n-1$ and 
$\alpha_n := e_{n-1}+e_n$, is a set of simple roots. The
corresponding set of positive roots is 
$\Phi^+ = \{ e_i \pm e_j \, | \, 1 \le i < j \le n \}$. 

From the explicit description of the root systems (or the Dynkin 
diagram) we see that $\{\alpha_2, \alpha_3, \dots, \alpha_n\}$
generates a root subsystem~$\Phi_{1'}$ of type $D_{n-1}$,  that 
$\{\alpha_1, \alpha_2, \dots, \alpha_{n-2},\alpha_{n}\}$ generates 
a root sub\-system~$\Phi_{(n-1)'}$ of type~$A_{n-1}$,
and that $\{\alpha_1, \alpha_2, \dots, \alpha_{n-1}\}$ generates 
a root sub\-system~$\Phi_{n'}$ of type~$A_{n-1}$. We set 
$\Phi^+_{1'} := \Phi_{1'} \cap \Phi^+$, $\Psi^+_{1'} := \Phi^+
\setminus \Phi^+_{1'}$, $\Phi^+_{(n-1)'} := \Phi_{(n-1)'} \cap \Phi^+$,
$\Psi^+_{(n-1)'} := \Phi^+ \setminus \Phi^+_{(n-1)'}$, and $\Phi^+_{n'} := \Phi_{n'} \cap \Phi^+$,
$\Psi^+_{n'} := \Phi^+ \setminus \Phi^+_{n'}$. The sets 
$\Psi^+_{1'}$, $\Psi^+_{(n-1)'}$, and $\Psi^+_{n'}$ are normal patterns. In fact,
$P(\Psi^+_{1'})$, $\Psi^+_{(n-1)'}$ and $P(\Psi^+_{n'})$ are the unipotent radicals of
the standard parabolic subgroups corresponding to 
$\{\alpha_2, \alpha_3, \dots, \alpha_n\}$, 
$\{\alpha_1, \alpha_2, \dots, \alpha_{n-2},\alpha_{n}\}$,  and
$\{\alpha_1, \alpha_2, \dots, \alpha_{n-1}\}$, respectively.
Furthermore, the sets $\Phi^+_{1'}$, $\Phi^+_{(n-1)'}$ and $\Phi^+_{n'}$ are patterns and
$P(\Phi^+_{1'}) \cong UD_n/P(\Psi^+_{1'}) \cong UD_{n-1}$
and $P(\Phi^+_{n'}) \cong UD_n/P(\Psi^+_{n'}) \cong UA_{n-1}$.

Let $\alpha \in \Phi^+$ and $\chi \in \Irr(UD_n)_\alpha$. Suppose that  
$\alpha \in \Phi^+_{1'}$. It follows from Lemma~\ref{la:alpha_rep} and
Remark~\ref{rem:kSigma} (a), (c) that 
$P(\Psi^+_{1'}) \subseteq \Ker(\chi)$. Thus, we can identify $\chi$
with a single root character of 
$P(\Phi^+_{1'}) \cong UD_n/P(\Psi^+_{1'}) \cong UD_{n-1}$. In this
way, the classification and construction of the elements of
$\Irr(UD_n)_\alpha$ is reduced to the case $D_{n-1}$.
Similarly, if $\alpha \in \Phi^+_{(n-1)'}$ or $\alpha \in \Phi^+_{n'}$ we can identify $\chi$ with a
single root character of $P(\Phi^+_{(n-1)'})$ respectively $P(\Phi^+_{n'})$ and thereby obtain a
reduction to the case $A_{n-1}$ which has already been treated in
Section~\ref{subsec:An}. 
Hence, we only have to consider positive roots $\alpha$ which are not
contained in $\Phi^+_{1'} \cup \Phi^+_{n'}$, i.e., the roots
$e_1+e_i$, where $1 < i < n$. Finally we observe that $\Phi$ is a subsystem of 
the root system of type $B_n$ described in Section~\ref{subsec:Bn}. Thus most of what follows 
is a simple consequence of restricting from type $B_n$ to type $D_n$.

\begin{lem} \label{la:hook_Dn}
Let $\Phi$ be a root system of type $D_n$ as described above. For all
positive roots of the form $\alpha = e_1+e_i$ we have:
\[
h(\alpha) = \{\alpha\} \cup \{e_1-e_s, e_s+e_i \, | \, 1<s<i\} \cup
\{e_1 \pm e_s, e_i \mp e_s \, | \, i < s \le n\}.
\]
\end{lem}

\pr
This follows from Lemma \ref{la:hook_Bn} via restriction.
\epr

\medskip

Next we determine the closed patterns $n(\alpha)$ and $k(\alpha)$.

\begin{lem} \label{la:Dn_n_k}
Let $\Phi$ be a root system of type $D_n$ as described above. For all
positive roots of the form $\alpha = e_1+e_i$ the following is true:
\begin{enumerate}
\item[(a)] The sets $n(\alpha)$ and $k(\alpha)$ are normal patterns in $\Phi^+$.

\item[(b)] $n(\alpha) = \{ e_1+e_s \, | \, 1 < s < i\}$.

\item[(c)] If $i<n$ then $k(\alpha) = \{ e_s+e_t \, | \, 1 \le s < t < i\}$. 
\end{enumerate}
\end{lem}

\pr
This follows from Lemma \ref{la:Bn_n_k} via restriction.
\epr

\medskip

Let $\alpha \in \Phi^+$. As already mentioned above, if 
$\alpha \in \Phi^+_{(n-1)'} \cup \Phi^+_{n'}$ then the hook $h(\alpha)$, the arm
$a(\alpha)$ and the leg $\ell(\alpha)$ are defined as for type $A_n$. 
Also if $\alpha \in \Phi^+_{1'}$ then we can assume recursively that
$h(\alpha)$, $a(\alpha)$ and $\ell(\alpha)$ are already defined. 
For $\alpha = e_1+e_i$, where $1 < i \le n$, we define the arm and the 
leg as in Section~\ref{subsec:Bn} via restriction:
\begin{eqnarray*}
a(\alpha) & := & \{e_1-e_s \, | \, 1 < s < i \} \cup \{ e_i \pm  e_s
\, | \, i < s \le n \} \,\, \text{   and}\\ 
\ell(\alpha) & := & \{e_s+e_i \, | \, 1 <  s < i\} \cup \{ e_1 \pm
e_s \, | \, i < s \le n\}. 
\end{eqnarray*}
As for types $A_n$ and $B_n$, we will see that for each $\alpha \in \Phi^+$
the hook $h(\alpha)$ is a pattern. We set $H_\alpha := P(h(\alpha))$
and call it the \emph{hook subgroup} corresponding to~$\alpha$. 

\begin{lem} \label{la:Dn_hooks}
Let $\Phi$ be a root system of type $D_n$ as described at the
beginning of this section. For all $\alpha = e_1+e_i \in \Phi^+$
the following are true:
\begin{enumerate}
\item[(a)] The hook $h(\alpha)$, the arm $a(\alpha)$ and the leg
  $\ell(\alpha)$ are patterns. 

\item[(b)] The pattern subgroups $P(a(\alpha))$ and $P(\ell(\alpha))$
  are elementary abelian.

\item[(c)] $|a(\alpha)| = |\ell(\alpha)| = 2n-i-2$.

\item[(d)] If $\height(\alpha) > 1$ then the hook subgroup $H_\alpha$ is
  special of type $q^{1+2(2n-i-2)}$ and $[y, H_\alpha]=Z(H_\alpha)=X_\alpha$ 
  for all $y \in H_\alpha \setminus Z(H_\alpha)$. More specifically:
  For each $y \in H_\alpha \setminus Z(H_\alpha)$ there is some
  $\beta \in h(\alpha)$ such that 
  $\{[y,x_\beta(t)] \, | \, t \in \F_q\}=X_\alpha$.
\end{enumerate}
\end{lem}

\pr
This follows from Lemma \ref{la:Bn_hooks} via restriction.
\epr

\medskip

\begin{lem} \label{la:Dn_source}
Let $\Phi$ be a root system of type $D_n$ as described at the
beginning of this section. For all $\alpha = e_1+e_i \in \Phi^+$
the following are true:
\begin{enumerate}
\item[(a)] The source $s(\alpha)$ is a closed pattern.

\item[(b)] $\ell(\alpha) \cup k(\alpha) \unlhd s(\alpha)$.

\item[(c)] $\ell(\alpha) \cup \{\alpha\} \cup k(\alpha) \unlhd \Phi^+$.
\end{enumerate}
\end{lem}

\pr
This follows from Lemma \ref{la:Bn_hooks} via restriction.
\epr

\medskip

\begin{prop} \label{prop:Dn_singlechars}
Let $\Phi$ be a root system of type $D_n$ as described above. For
each root $\alpha \in \Phi^+$ the map
$\Psi_\alpha: \Irr(\overline{T}_\alpha) \times \Irr(X_\alpha)^*
\to \Irr(UD_n)_\alpha$ with 
\[
(\mu, \lambda) \mapsto (\infl_{\overline{T}_\alpha}^{S_\alpha} \mu \cdot
 \infl_{X_\alpha}^{S_\alpha} \lambda)^{UD_n}
\]
is a one to one correspondence with the property
$\Psi_\alpha(\mu, \lambda)(1) = q^{2n-i-1} \cdot \mu(1)$. 
\end{prop}
\pr
The content of Lemma \ref{la:Dn_source} is that 
is that hypothesis (1), (2) and (3) of Proposition \ref{prop:red_to_thm}
are satisfied, and the content of Lemma~\ref{la:Dn_hooks} (d) is that 
hypothesis (4) of Proposition \ref{prop:red_to_thm}
is satisfied. This proves the correspondence. 
\epr

We conclude this section by giving structural information concerning the group $\overline{T}_\alpha$. 
For all $1 < i \le n$ and $\alpha = e_1+e_i$ we define 
\begin{eqnarray*}
d_{1,i} & := & \{ e_r \pm e_s \mid 1 < r < s \le n \ \mbox{and}
\ r \neq i \neq s \} \ \text{and}\\
\ob(\alpha) & := & \Phi^+ \setminus (n(\alpha) \cup h(\alpha) \cup
d_{1,i}) = \{e_r-e_i \mid 1 \le r < i \} \ \text{and}\\
 \kob(\alpha) & := & d_{1,i} \setminus \{e_r-e_s \mid 1<r<s<i \}.
\end{eqnarray*}
Note that $d_{1,i}$ and $\ob(e_1+e_i)$ are patterns for all 
$1 < i \le n$ and that $P(\ob(\alpha))$ is elementary abelian.
Furthermore, we have $\ob(\alpha) \cap k(\alpha) = \emptyset$ for
$i<n$ and $\ob(\alpha) \subseteq k(\alpha)$ for $i=n$.

\begin{lem} \label{la:Dn_structure_sa}
Let $\Phi$ be a root system of type $D_n$ as described at the
beginning of this section. For all $\alpha=e_1+e_i \in \Phi^+$
the following are true:
\begin{enumerate}
\item [(a)] $S_\alpha/P(\ell(\alpha) \cup k(\alpha)) \cong
  \overline{S}_\alpha/P(\ell(\alpha))P(k(\alpha)) \cong \overline{T}_\alpha \times X_\alpha$.

\item[(b)] $\overline{T}_\alpha \cong P(d_{1,i})/P(d_{1,i} \cap k(\alpha))\ltimes
  P(\ob(\alpha))/P(\ob(\alpha) \cap k(\alpha))$. 

\item[(c)] If $i<n$, then $\kob(\alpha)$ is a pattern that is normal in $d_{1,i}$ and
  $P(\kob(\alpha))$ centralizes $P(\ob(\alpha))/P(\ob(\alpha) \cap k(\alpha))$.

\item[(d)] If $i=n$ then $P(\ob(\alpha))/P(\ob(\alpha) \cap k(\alpha))
  = \{1\}$.

\item[(e)] If $i \in \{2,3\}$ and $i=n$ then $\overline{T}_\alpha = \{1\}$.

\item[(f)] If $i \in \{2,3\}$ and $i<n$ then $\kob(\alpha)=d_{1,i}$,
  $k(\alpha) = n(\alpha)$ and
  \[
  \overline{T}_\alpha \cong UD_{n-2} \times (UA_1)^{i-1}.
  \]

\item[(g)] If $i>3$ then $P(d_{1,i})/P(d_{1,i} \cap k(\alpha))$ is
  isomorphic to a quotient pattern group of $UD_{n-2}$.

\item[(h)] If $i>3$ then $P(d_{1,i})/P(\kob(\alpha) \cup (k(\alpha)
  \cap d_{1,i})) \cong UA_{i-3}$.
\end{enumerate}
\end{lem}
\pr
This follows from Lemma \ref{la:Bn_structure_sa} via restriction.
\epr

\begin{thm} \label{thm:Dn_singlechars}
Let $\Phi$ be a root system of type $D_n$ as described at the
beginning of this section. For all $\alpha = e_1+e_i$ there is a
one to one correspondence 
\[
\Psi: \Irr(\overline{T}_\alpha) \times \Irr(\F_q)^* \rightarrow \Irr(UD_n)_\alpha
\]
such that $\deg(\Psi(\chi, \lambda)) = q^{2n-i-2} \deg(\chi)$.
\end{thm}

\pr
The theorem follows from Proposition~\ref{prop:Dn_singlechars}
and Lemmas \ref{la:Dn_structure_sa} (a) and \ref{la:Dn_hooks}~(c). 

\epr


\subsection{Type C} \label{subsec:Cn}

Let $n \ge 3$ be an integer. We construct a root system of
type~$C_n$ as in~\cite[Section~12.1]{Humphreys}: Let 
$e_1, e_2, \dots, e_n \in \R^n$ be the usual orthonormal unit
vectors which form a basis of $\R^n$. Then
$\Phi := \{ \pm(e_i \pm e_j) \, | \, 1 \le i \neq j \le n \} \cup \{ \pm 2e_i \, | \, 1 \le i \le n\}$ is a
root system of type~$C_n$ and the set $\{\alpha_1, \dots, \alpha_n\}$,
where $\alpha_i := e_i-e_{i+1}$ for $i=1,2,\dots,n-1$ and 
$\alpha_n := 2e_n$, is a set of simple roots. The
corresponding set of positive roots is 
$\Phi^+ = \{ e_i \pm e_j \, | \, 1 \le i < j \le n \} \cup \{ 2e_i \, | \, 1 \le i \le n \}$.

The highest long root with respect to this base is $2e_1 = 2\alpha_1 + \dots 2\alpha_{n-1} + \alpha_n$ 
whereas the highest short root is $e_1+e_2$.  We note that the long roots of $\Phi$ form a $A_1^n$-subsystem.
Also we recall that the Chevalley commutator relations imply that if $\a$ is short and $\b$ is 
long, then $[X_\a,X_\b] \subset X_{\a +\b}X_{2\a+b}$, where $\a+\b$ is short and $2\a + \b$ is 
long, and that when both $\a$ and $\b$ are short, then $[X_\a,X_\b] \subset X_{\a +\b}$, where 
$\a+\b$ is long.

From the explicit description of the root systems (or from the Dynkin 
diagram) we see that $\{\alpha_2, \alpha_3, \dots, \alpha_n\}$
generates a root subsystem~$\Phi_{1'}$ of type $C_{n-1}$ and that 
$\{\alpha_1, \alpha_2, \dots, \alpha_{n-1}\}$ generates 
a root sub\-system~$\Phi_{n'}$ of type~$A_{n-1}$. We set 
$\Phi^+_{1'} := \Phi_{1'} \cap \Phi^+$, $\Psi^+_{1'} := \Phi^+
\setminus \Phi^+_{1'}$ and $\Phi^+_{n'} := \Phi_{n'} \cap \Phi^+$,
$\Psi^+_{n'} := \Phi^+ \setminus \Phi^+_{n'}$. The sets 
$\Psi^+_{1'}$ and $\Psi^+_{n'}$ are normal patterns. In fact,
$P(\Psi^+_{1'})$ and $P(\Psi^+_{n'})$ are the unipotent radicals of
the standard parabolic subgroups corresponding to 
$\{\alpha_2, \alpha_3, \dots, \alpha_n\}$ and
$\{\alpha_1, \alpha_2, \dots, \alpha_{n-1}\}$, respectively.
Furthermore, the sets $\Phi^+_{1'}$ and $\Phi^+_{n'}$ are patterns and
$P(\Phi^+_{1'}) \cong UC_n/P(\Psi^+_{1'}) \cong UC_{n-1}$
and $P(\Phi^+_{n'}) \cong UC_n/P(\Psi^+_{n'}) \cong UA_{n-1}$.

Let $\alpha \in \Phi^+$ and $\chi \in \Irr(UC_n)_\alpha$. Suppose that  
$\alpha \in \Phi^+_{1'}$. It follows from Lemma~\ref{la:alpha_rep} and
Remark~\ref{rem:kSigma} (a), (c) that
$P(\Psi^+_{1'}) \subseteq \Ker(\chi)$. Thus, we can identify $\chi$
with a single root center character of 
$P(\Phi^+_{1'}) \cong UC_n/P(\Psi^+_{1'}) \cong UC_{n-1}$. In this
way, the classification and construction of the elements of
$\Irr(UC_n)_\alpha$ is reduced to the case $C_{n-1}$.
Similarly, if $\alpha \in \Phi^+_{n'}$ we can identify $\chi$ with a
single root center character of $P(\Phi^+_{n'})$ and thereby get a
reduction to the case $A_{n-1}$ which has already been treated in
Subsection~\ref{subsec:An}. 
Hence, we only have to consider positive roots $\alpha$ which are not
contained in $\Phi^+_{1'} \cup \Phi^+_{n'}$, i.e., the roots
$e_1 + e_i$ where $1 < i \le n$ and the root $2e_1$.
We begin with single root characters lying of the longest root $2e_1$. 

\begin{lem}\label{la:hook_Clong}
Let $\Phi$ be a root system of type $C_n$ as described above. For
 $\alpha = 2e_1$ we have
\begin{equation} 
h(\alpha) = \{2e_1, e_1+e_i, e_1 - e_i \ | \ 2 \leq i \leq n  \}.
\end{equation}
\end{lem}
\pr We observe that the sum of two long roots is never in $\Phi^+$. 
Next we observe that the sum of a long and a short root is either short or 
not in $\Phi^+$. Thus if $2e_1 = \beta+\gamma$, then both $\beta$ and $\gamma$ are short. 
This implies that $\beta = e_1 \pm e_s$ and hence $\gamma = e_1 \mp e_s$.
This proves our claim. 
\epr
\medskip

We define the arm $a(2e_1)$ and the leg $\ell(2e_1)$ of $h(2e_1)$
as follows
\[
   a(2e_1) := \{ e_1 - e_s \, | \, 1 < s \leq n \} \quad \text{and} \quad 
\ell(2e_1) := \{ e_1 + e_j \, | \, 1 < s < \leq n \}.
\]

\begin{lem}\label{la:c_long_prep}
The following  are true:
\begin{enumerate}
\item[(a)] $n(2e_1) = k(2e_1) = \emptyset$. 
\item[(b)] $\ell(\alpha) \cup \{\alpha \} \unlhd h(\alpha)$.
\item[(c)] The group $P(h(2e_1))$ is special of order $q^{1+2(n-1)}$. Moreover 
$[y, H_{2e_1}]=Z(H_{2e_1})=X_{2e_1}$  for all $y \in H_{2e_1}\setminus Z(H_{2e_1})$.
\item[(d)] $s(2e_1)$ is a closed pattern.
\item[(e)] $\ell(2e_1) \unlhd s(2e_1)$ and $P(\ell(2e_1))$ is abelian. 
\item[(f)] $s(2e_1) \setminus (\{2e_1\} \cup \ell(2e_1))$ is a pattern of type $UC_{n-1}$,
$\overline{T}_{2e_1} \cong UC_{n-1}$, and $S_{2e_1}/ P(\ell(2e_1)) \cong \overline{T}_{2e_1} \times X_{2e_1}$.

\end{enumerate}

\end{lem}
\pr
As $2e_1$ is the highest root of $\Phi^+$ part (a) is clear. Also we note that 
$2e_1$ is the unique long root in $h(2e_1)$. The short roots in $h(\alpha)$ are all of the 
form $e_1 \pm e_r$ with $r > 1$ and thus sum of two of these lies in $\Phi^+$
if and only if that sum is $2e_1$. This combined with Lemma \ref{la:derquopat} 
proves part (b) and the first half of part (c). The proof of the Lemma
\ref{la:An_hooks} (d) carries over verbatim to prove the second part
of (c). 

Now let $\beta := e_1-e_s \in a(2e_1)$ and suppose that $\b = \gamma + \gamma'$ with 
$\gamma,\gamma' \in \Phi^+$. Then, as $\beta \in \Phi_{n'}$, we see that without loss 
$\gamma = e_1 - e_r$ and $\gamma' = e_r - e_s $ with $ 1 < r < s$.
As $\gamma \in a(2e_1)$ this proves (d).
Now if $e_1 + e_r, e_1+e_s \in \ell(2e_1)$, then $e_1 + e_r + e_1+e_s = 2e_1 + e_r + e_s \not \in \Phi^+$, 
and thus the Chevalley commutator relations imply that $P(\ell(2e_1))$ is elementary abelian. 
If $e_1 + e_s \in \ell(2e_1)$ and $\beta \in s(2e_1)$, then $e_1 + e_s + \beta \in \Phi^+$ only 
if $\beta = e_r -e_s$ where $r < s$ and then $ e_1 + e_s + \beta \in \ell(2e_1)$ proving (e).

Finally we observe that 
\[
s(2e_1) \setminus (\{2e_1\} \cup \ell(2e_1)) = \{ 2e_i  \, | \, 2 \leq i \leq n \} \cup \{ e_r \pm e_s  \, | \, 2 \leq r < s \leq n \}
\]
from which the claims in (f) follow. 
\epr

\begin{prop} \label{prop:Cn_singlechars_long}
Let $\Phi$ be a root system of type $C_n$ as described above. For
the root $\alpha:=2e_1 \in \Phi^+$ the map
$\Psi_\alpha: \Irr(\overline{T}_\alpha) \times \Irr(X_\alpha)^*
\to \Irr(UC_n)_\alpha$ with 
\[
(\mu, \lambda) \mapsto (\infl_{\overline{T}_\alpha}^{S_\alpha} \mu \cdot
 \infl_{X_\alpha}^{S_\alpha} \lambda)^{UC_n}
\]
is a one to one correspondence with the property
$\Psi_\alpha(\mu, \lambda)(1) = q^{n-1} \cdot \mu(1)$, 
where $\overline{T}_\alpha \cong UC_{n-1}$.
\end{prop}
\pr
The content of Lemma \ref{la:c_long_prep} (d), (e), (b) and (c) is that 
is that hypothesis (1), (2), (3) and (4) of Proposition \ref{prop:red_to_thm}
are satisfied. This proves the correspondence. The statement 
about the structure of $\overline{T}_\alpha$ is the content of 
Lemma \ref{la:c_long_prep} (f).
\epr

\medskip

Next we consider the short roots $\alpha = e_1+e_i $ which are not
contained in $\Phi^+_{1'} \cup \Phi^+_{n'}$.

\begin{lem} \label{la:hook_Cn_short}
Let $\Phi$ be a root system of type $C_n$ as described above. For all
positive roots of the form $\alpha = e_1+e_i$ we have:
\[
 h(\alpha) = \{\alpha\} \cup \{e_1 -e_s, e_s+e_i \ | \ 1 < s < i \} \cup 
\{e_1 \pm e_s, e_i \mp e_s \ | \ i < s \leq n \} \cup \{2e_i, (e_1-e_i) \}
\]
\end{lem}

\pr Let $\beta=2e_t$, $\beta' = e_{s'} \pm e_{t'} \in \Phi^+$ such
that $\beta+\beta'=\alpha$. Then $\beta'= e_1 - e_i$ and $\beta = 2e_i$. 

Next let $\beta=e_s \pm e_t$, $\beta' = e_{s'} \pm e_{t'} \in \Phi^+$ such
that $\beta+\beta'=\alpha$. We see immediately that one of the
two $\pm$ signs has to be a $+$~sign and that the other one has to be
a $-$~sign and that $t=t'$. Furthermore, we can assume $s=1$ and
$s'=i$. Hence 
\[
(\beta, \beta') \in \{(e_1-e_s, e_s+e_i) \, | \, 1 < s
< i \,\, \text{or} \,\, i < s \le n\} \cup \{(e_1+e_s, e_i-e_s) \, | \, i
< s \le n\}
\]
and the claim follows.
\epr

\medskip

Recall that the contents of Lemma~\ref{la:representable} and Lemma~\ref{la:alpha_rep} is that if 
$\chi \in \Irr(XU)_{\alpha}$, then $\rk(\chi) = k(\alpha).$ 
Also recall that Lemma \ref{la:k_w} states $w(\alpha) \cap k(\alpha)=\emptyset$ for each 
$\alpha \in \Phi^+$. 
Using these two facts, we can now describe the patterns $n(\alpha)$ and $k(\alpha)$.

\begin{lem} \label{la:Cn_n_k}
Let $\Phi$ be a root system of type $C_n$ as described above and let $\alpha = e_1+e_i$. 
The following are true:
\begin{enumerate}
\item[(a)] The sets $n(\alpha)$ and $k(\alpha)$ are normal patterns in $\Phi^+$. 

\item[(b)] $n(\alpha) = \{e_1+e_s \, | \,  1 < s < i \} \cup \{ 2e_1\} $.

\item[(c)] $k(\alpha) = n(\alpha) \cup \{ e_r + e_s \, | \, 1 < r < s < i \} \cup \{2e_r \, | \, 1 < r < i \} $.
\end{enumerate}
\end{lem}

\pr
Part (a) was established in Section~\ref{sec:representable}.

\medskip

\noindent (b) By definition $n_0(\alpha) = \{e_1+e_i\}$.
Let $\gamma = e_k \pm e_l \in \Phi^+$. Then
$\alpha + \gamma  \in \Phi^+$ if and only if $\gamma = e_s-e_i$ for
some $1 < s <i$ and in this case we have $\alpha + \gamma = e_1+e_s$.
Now if $\gamma' = 2e_r$, then $\gamma' + \delta \not \in \Phi^+$ for all $\delta = e_r + e_t$.
Thus, $n_1(\alpha) = \{e_1+e_s \, | \, 1 < s < i\}$. Again, if 
$\gamma = e_k \pm e_l \in \Phi^+$ then $\gamma + (e_1+e_s)  \in \Phi^+$
if and only if $\gamma = e_t-e_s$ for 
some $1 < t <s$ and in this case we have $\gamma + (e_1+e_s) = e_1+e_t$.
Thus $n_2(\alpha) = n_1(\alpha)$ and therefore
$n(\alpha) = n_1(\alpha)$ and (b) follows. 

\medskip

\noindent (c) Let $\beta = e_s + e_t \in \Phi^+$ where $s<t$ and 
 $t \ge i$. If $s=1$ then $\beta \in h(\alpha) \subseteq w(\alpha)$.
 If $s>1$ then $\beta \in h(e_1+e_t)$ and $e_1+e_t \in h(\alpha)$,
 hence $\beta \in w(\alpha)$.  

Let $\beta = 2e_s$ with $i < s$, then 
$\beta + (e_1 - e_s) = e_1 + e_s \in h(\alpha) \subseteq w(\alpha)$.

 Now let $\beta = e_s - e_t \in \Phi^+$ where $s<t$ and 
 $t \ge i$. If $s=1$ and $t \neq i$ then 
 $\beta \in h(\alpha) \subseteq w(\alpha)$. If $s=1$ and $t=i<n$ then
 $\beta \in h(e_1+e_n)$ and $e_1+e_n \in h(\alpha)$. Thus, 
 $\beta \in w(\alpha)$. If $s>1$ and $t \neq i$ then 
 $\beta \in h(e_1-e_t)$ and $e_1-e_t \in h(\alpha)$. Hence 
 $\beta \in w(\alpha)$. If $s>1$ and $t=i<n$ then $\beta \in
 h(e_1-e_t)$, $e_1-e_t \in h(e_1-e_n)$ and $e_1-e_n \in h(\alpha)$ so
 we have $\beta \in w(\alpha)$.

It follows from Lemma~\ref{la:k_w} that 
$k(\alpha) \subseteq \{ e_s+e_t \, | \, 1 \le s < t < i\} \cup \{2e_s \, | \, s < i \}.$

Let $\beta \in \Phi^+ \setminus \{\alpha\}$ and let $i',j'$ be
nonnegative integers. We can see from the explicit
description of the root system of type $C_n$ that for all 
$\gamma \in \Phi^+$ we have  $i' \beta + j' \gamma \in \Phi^+$ only if
$i',j' \le 2$. Thus, the definition of $k(\alpha)$,
Lemma~\ref{la:cenquopat} and the commutator relations imply that
$\beta \in k(\alpha)$ if and only if for all $\gamma \in \Phi^+$ we
have $\beta+\gamma \not\in \Phi^+$ or $\beta+\gamma \in k(\alpha)$. 
A straightforward induction on $s+t$ shows that 
$2e_s, e_s+e_t \in k(\alpha)$ for all $1 \le s < t < i$ (where the induction
beginning $s=1$, $t=2$ follows from $n(\alpha) \subseteq k(\alpha)$).

This completes the proof.
\epr

We define 
\[ 
\ell(e_1+e_i) := \{ e_1 \pm e_s  \ | \ i < s \leq n \} \cup \{ e_1 - e_i\} \cup \{e_s+e_i \ | \  1 < s < i \} \]
and
 \[
 a(e_1+e_i):= \{ e_i \mp e_s \ | \ i < s \leq n \} \cup \{ 2e_i \} \cup \{e_1-e_s \ | \ 1 < s < i  \} .\]

As before we will show that for each $\alpha \in \Phi^+$
the hook $h(\alpha)$ is a pattern modulo $k(\alpha)$. Call $H_\alpha := P(h(\alpha))$
the \emph{hook subgroup} corresponding to~$\alpha$, always bearing in mind that we
calculate modulo $P(k(\alpha))$. 

\begin{lem} \label{la:Cn_hooks}
Let $\Phi$ be a root system of type $C_n$ as described at the
beginning of this subsection. For all $\alpha = e_1+e_i \in \Phi^+$
the following are true:
\begin{enumerate}
\item[(a)] The hook $h(\alpha) \cup k(\alpha)$ and the leg
  $\ell(\alpha)\cup k(\alpha)$ are closed patterns.

\item[(b)] $\{\alpha \} \cup \ell(\alpha) \cup k(\alpha) \unlhd h(\alpha) \cup k(\alpha)$.

\item[(c)] $P(a(\alpha))$ is a pattern subgroup. It is isomorphic to a product of  
   an elementary abelian group of order $q^{i-2}$ and a special of order $q^{1+2(n-i)}$.
The center of the special group is $X_{2e_i}$. 

\item[(d)] $|a(\alpha)| = 2n-i-1$.

\item[(e)] 
If $\height(\alpha) > 1$, then 
for each $y \in (\prod_{\gamma \in a(\alpha)} X_\gamma) \setminus \{1\}$
there is some $\beta \in \ell(\alpha)$ such that $\{[y,x_\beta(t)] \, | \, t \in \F_q\}=X_\alpha$.
\end{enumerate}
\end{lem}

\pr
(a), (b): Let $\beta, \gamma \in h(\alpha)$. We have
$\beta+\gamma \in \Phi^+$ if and only if 
$\{\beta, \gamma\} = \{e_1-e_s, e_s+e_i\}$ for some $s \neq i$ or
$\{\beta, \gamma\} = \{e_1+e_s, e_i-e_s \}$ for some $s>i$ or 
$\{\beta, \gamma\} = \{e_i+e_s, e_i-e_s \}$ for some $s>i$ or 
$\{\beta, \gamma\} = \{e_1+e_s, e_1-e_s \}$ for some $s>i$ or 
$\{\beta, \gamma\} = \{2e_i, e_1-e_i \}$. 
In all of these cases $\beta+\gamma = \alpha \in h(\alpha) \cup \{2e_1\}$,
and (a) follows. 
\medskip 

If $\beta \in h(\alpha)$ and $ \gamma \in \ell(\alpha)$, then either 
$ \beta + \gamma \in \{\alpha, 2e_1\}$ or
$ \beta + \gamma \not \in \Phi^+ $, and (b) follows. 

\medskip

Also for $\beta, \gamma \in a(\alpha)$ we see that $\beta + \gamma \in \{2e_i \}$ if 
$\beta, \gamma \in a(\alpha) \setminus \{e_1 - e_s \ | \ s < i \}$ and that 
$\beta, \gamma \not \in \Phi^+$ if $\beta, \gamma \in \{e_1 - e_s \ | \ s < i \}$.
This shows (c).

\medskip 
(d) follows from (c).

\medskip 
Lastly suppose that $\height(\alpha)>1$. We have seen in the
proof of (a) and (b) that for all $\beta, \gamma \in h(\alpha)$
we have $\beta+\gamma \in \Phi^+$ if and only if 
$\{\beta, \gamma\} = \{e_1-e_s, e_s+e_i\}$ for some $s \neq i$ or
$\{\beta, \gamma\} = \{e_1+e_s, e_i-e_s \}$ for some $s>i$ or 
$\{\beta, \gamma\} = \{e_i+e_s, e_i-e_s \}$ for some $s>i$ or 
$\{\beta, \gamma\} = \{e_1+e_s, e_1-e_s \}$ for some $s>i$ or 
$\{\beta, \gamma\} = \{2e_i, e_1-e_i \}$. 
In all of these cases $\beta+\gamma = \alpha \in h(\alpha) \cup \{2e_1\} $.
It follows that
$Z(H_\alpha) = X_\alpha$ and $[H_\alpha, H_\alpha] \subseteq
X_\alpha \times X_{2e_i}$ modulo $P(k(\alpha))$. 
Moreover it follows that $P(\{\alpha \} \cup \ell(\alpha))$ is normal
and elementary abelian. Now let 
$y = \prod_{\gamma \in h(\alpha)} x_\gamma(t_\gamma) \in (\prod_{\gamma \in a(\alpha)} X_\gamma) \setminus \{1\}$. 
Because $y \neq 1$ there is some $\gamma \in a(\alpha) \setminus \{\alpha\}$ such $t_\gamma \neq 0$.
If $\gamma \neq 2e_i$ then we pick $\beta := \alpha-\gamma \in h(\alpha)$ and we get from
Lemma~\ref{la:patnorm}~(a) that 
$\{[y,x_\beta(t)] \mid t \in \F_q\} = X_\alpha = Z(H_\alpha).$
If $y \in X_{2e_i}P(\{\alpha \} \cup \ell(\alpha) \cup k(\alpha)) \setminus P(\{\alpha \} \cup \ell(\alpha) \cup k(\alpha))$, 
then we pick $\beta := e_1-e_i$ and again we get that
$\{[y,x_\beta(t)] \mid t \in \F_q\} = X_\alpha = Z(H_\alpha).$
\epr 

\medskip

Next we show 

\begin{lem} \label{la:Cn_source}
Let $\Phi$ be a root system of type $C_n$ as described at the
beginning of this subsection. For all $\alpha = e_1+e_i \in \Phi^+$
the following is true:
\begin{enumerate}
\item[(a)] The source $s(\alpha)$ is a closed pattern and $k(\alpha)
  \subseteq s(\alpha)$.

\item[(b)] $\ell(\alpha) \unlhd s(\alpha)$.

\item[(c)] The canonical projection $\pi: UC_n \rightarrow
  UC_n/P(k(\alpha))$ maps the hook subgroup $H_\alpha=P(h(\alpha))$
  injectively into $UC_n/P(k(\alpha))$ and 
  $\pi(P(\ell(\alpha)) \unlhd \overline{S}_\alpha$. 

\item[(d)] $X_\alpha \subseteq Z(\overline{S}_\alpha)$.
\end{enumerate}
\end{lem}
\pr 
(a) By Lemma~\ref{la:k_w} we have $k(\alpha) \cap a(\alpha) \subseteq
k(\alpha) \cap w(\alpha) = \emptyset$. Hence 
$k(\alpha) \subseteq \Phi^+ \setminus a(\alpha) = s(\alpha)$.
Let $\beta \in a(\alpha)$ and $\gamma, \gamma' \in \Phi^+$ such that
$\beta=\gamma+\gamma'$.

To prove that $s(\alpha)$ is a pattern it suffices to show that
$\gamma \in a(\alpha)$ or $\gamma' \in a(\alpha)$. Suppose that 
$\beta=e_1-e_s$ where $1 < s < i$. Then 
$\{\gamma, \gamma'\} = \{e_1-e_l, e_l-e_s\}$ for some $1<l<s<i$ and
hence $\gamma \in a(\alpha)$ or $\gamma' \in a(\alpha)$. Now suppose
that $\beta = e_i+e_s$ where $i < s \le n$. Then 
$\{\gamma, \gamma'\} = \{e_i+e_l, e_s-e_l\}$ for some $i < l \le n$ or 
$\{\gamma, \gamma'\} = \{e_s+e_l, e_i-e_l\}$ for some $i < l \le n$. 
Hence in both cases $\gamma \in a(\alpha)$ or $\gamma' \in a(\alpha)$.
Next, suppose that $\beta = e_i-e_s$ where $i < s \le n$. Then 
$\{\gamma, \gamma'\} = \{e_i-e_l, e_l-e_s\}$ for some $i < l < s \le n$ 
and again $\gamma \in a(\alpha)$ or $\gamma' \in a(\alpha)$. 
Finally, suppose that $\beta = 2e_i$, then $\{\gamma, \gamma'\} = \{e_i+e_s, e_i-e_s\}$ 
for some $i < s \leq n$. But then $\gamma, \gamma' \in a(\alpha)$.
This proves (a).

\medskip

\noindent 
Let $\beta \in \ell(\alpha)$ and $\gamma \in s(\alpha)$.
First suppose
that $\beta = e_1+e_s$ where $i < s \le n$. Then $\beta + \gamma \in \Phi^+$
only if $\gamma = e_r - e_s$ with $i \neq r < s$ (as $e_i - e_s \in a(\alpha)$) and thus 
$\beta + \gamma = e_1 + e_r \in k(\alpha) \cup \ell(\alpha)$.

Next suppose that $\beta = e_1-e_s$ where $i < s \le n$. Then $\beta + \gamma \in \Phi^+$
only if $\gamma = e_r + e_s$ with $i \neq r $ (as $e_i + e_s \in a(\alpha)$) and thus 
$\beta + \gamma = e_1 + e_r \in k(\alpha) \cup \ell(\alpha)$.

Next suppose that $\beta = e_s+e_i$ where $s < i$. Then $\beta + \gamma \in \Phi^+$
only if $\gamma = e_r - e_s$ with $r < s $ or  $\gamma = e_r - e_i$ with $r < i $. 
Thus 
$\beta + \gamma \in \{ e_r + e_s \ | \ r< s < i \} \cup \{e_r + e_i \ | \ r < i \} \subset  k(\alpha) \cup \ell(\alpha)$.

Finally suppose that $\beta = e_1-e_i$. Then $\beta + \gamma \in \Phi^+$ only 
if $\gamma = 2e_i$, $e_1+e_i$ or $e_i + e_r$ or $e_i-e_r$ with $r > i$.
As $\gamma \in s(\alpha)$ we see that $\gamma = e_1+e_i$ or $e_i + e_r$ with $r < i$.
Thus 
$\beta + \gamma \in \{ 2e_1 \} \cup \{ e_1 + e_r \ | \ r< i \}  \subset  k(\alpha) \cup \ell(\alpha)$.

Now (b) follows.

\medskip

\noindent (c) By Lemma~\ref{la:k_w} we have 
$h(\alpha) \cap k(\alpha) \subseteq w(\alpha) \cap k(\alpha) =
\emptyset$. Since $H_\alpha$, $P(k(\alpha))$ are pattern subgroups we
get that the restriction of $\pi$ to $H_\alpha$ is injective. 

Let $\beta \in \ell(\alpha)$ and $\gamma \in S_\alpha$ such that
$\beta+\gamma \in \Phi^+$. We have to show that 
$\beta + \gamma \in \ell(\alpha) \cup k(\alpha)$. If
$\beta = e_s+e_i$ where $1<s<i$ then $\gamma$ has to be of the form 
$\gamma=e_k-e_l$ where $k < l$. Thus $\beta+\gamma=e_k+e_i$ where
$1<k<i$ or $\beta+\gamma=e_s+e_k$ where $k<i$. In the first case we
have $\beta+\gamma \in \ell(\alpha)$ and in the second case 
$\beta+\gamma \in k(\alpha)$. Now suppose that $\beta=e_1+e_s$ where
$s>i$. Then $\gamma$ has to be of the form $\gamma=e_k-e_l$ where 
$k < l$ and $\beta+\gamma=e_1+e_k$ where $k\neq1, i$. Hence
$\beta+\gamma \in k(\alpha) \cup \ell(\alpha)$. 
Finally, suppose
that $\beta=e_1-e_s$ where $s \geq i$ and $\gamma=e_k-e_l$ where
$k<l$. Then $\beta+\gamma=e_1-e_l$ where $l>i$. Thus 
$\beta+\gamma \in \ell(\alpha)$.

\medskip

\noindent (d) By definition of $n(\alpha)$ we have 
$\alpha+\gamma \in n(\alpha) \subseteq k(\alpha)$ for all 
$\gamma \in \Phi^+$ and so (d) follows from the commutator relations. 
\epr

\medskip

We can now complete the proof of Theorem \ref{thm:main_classical}.

\begin{thm} \label{Cn_singlechars}
Let $\Phi$ be a root system of type $C_n$ as described above. For
each root $\alpha \in \Phi^+$ the map
$\Psi_\alpha: \Irr(\overline{T}_\alpha) \times \Irr(X_\alpha)^*
\to \Irr(UA_n)_\alpha$ with 
\[
(\mu, \lambda) \mapsto (\infl_{\overline{T}_\alpha}^{S_\alpha} \mu \cdot
 \infl_{X_\alpha}^{S_\alpha} \lambda)^{UA_n}
\]
is a one to one correspondence with the property
$\Psi_\alpha(\mu, \lambda)(1) = q^{2n-i-1} \cdot \mu(1)$. 
\end{thm}
\pr If $\alpha$ is long then our claim is the content of Proposition \ref{prop:Cn_singlechars_long}.
If $\alpha$ is short, then 
the content of Lemma \ref{la:Cn_source} (a) and (b) is that 
is that hypothesis (1) and (2) of Proposition \ref{prop:red_to_thm}
are satisfied, and the content of Lemma~\ref{la:Cn_hooks} (b) and (e) is that 
hypotheses (3) and (4) of Proposition \ref{prop:red_to_thm}
is satisfied. This proves the correspondence. 
\epr

\medskip

We now derive some structural information about $\overline{T}_\alpha$. To this end we define patterns 
$$c_{i,j}:= \{ 2e_t, e_r \pm e_s \ | \ r,s,t \not \in \{ i,j\}, \ \mbox{and} \ r < s \}$$
where $1 \leq i<j \leq n$ and
$$ob(e_1+e_i):= \{ e_2-e_i ,\dots e_{i-1}-e_i \} \ \mbox{and} \ kob = c_{1,i} \setminus \{e_r-e_s \, | \, 1 < r < s < i \}.$$
Also we define $\overline{T}_\alpha := S_\alpha/ P(\{\alpha\} \cup \ell(\alpha) \cup k(\alpha))$.

\begin{lem} \label{la:Cn_structure_sa}If $\alpha = e_1+e_i$, then 
$S_\alpha/P(\ell(\alpha) \cup k(\alpha)) \cong \overline{S}_\alpha/P(\ell(\alpha)) \cong \overline{T}_\alpha \times X_\alpha$ 
and $\overline{T}_\alpha$ is a semidirect product of
$P(ob(\alpha))$ with $P(c_{1,i})/P(c_{1,i} \cap k(\alpha))$, where $c_{1,i}$ is as above. 
Moreover the following are true. 
\begin{enumerate}
\item [(a)] The kernel of the action of $P(c_{1,i})$ on $P(ob)$ is $P(kob)$. 
\item [(b)] For $i> 2$, $P(c_{1,i}/(c_{1,i} \cap k(\alpha))) \cong  UC_{n-2}/Z_{i-2}(UC_{n-i})$, 
\item [(c)] For $i > 2$, $P(c_{1,i})/P(kob) \cong  UA_{i-3},$
\item [(d)] If $i = 2$, then $kob = \emptyset = ob$, $k(\alpha) = n(\alpha) = \{ 2e_1 \}$, and 
$\overline{T}_\alpha \cong UC_{n-2}.$
\item[(e)]  If $i = 3$, then $kob = \emptyset$, $k(\alpha) = n(\alpha) = \{ 2e_1,e_1+e_2 \}$, and \\
$\overline{T}_\alpha \cong UC_{n-2}/Z(UC_{n-2}) \times X_{e_2-e_3}$.
\end{enumerate}

\end{lem}
\pr
The fact that $\overline{S}_\alpha/P(\ell(\alpha)) \cong \overline{T}_\alpha \times X_\alpha$ follows 
from the definition of 
$\overline{T}_\alpha$ and parts (a) and (b) of Lemma \ref{la:Cn_source}.

\medskip

(a) The set $ob(\alpha)$ is the set of roots of $\Phi^+$ which are not contained in 
$n(\alpha) \cup h(\alpha) \cup c_{1,i}$. Notice that $ob(\alpha)$ is normalized 
by $c_{1,i}$; proving the first part of (a).

For the second part we observe that no root of $c_{1,i}$ has $e_i$ in its support. Thus only roots from 
$c_{1,i}$ that can be added to an element of $ob(\alpha)$ to yield an element of $\Phi^+$ are of the form 
$e_r - e_s$ where $1 < r < s < i$. This yields the second part of (a). 

\medskip
(b)  We note that $c_{1,i} \subset \Phi^+_{1'}$ and so our claim follows from 
Part (c) of Lemma \ref{la:Cn_n_k}.

\medskip
(c)  We observe that $r,s \neq 1,i$ then $(e_r - e_s) + (e_k - e_i) \in \Phi$ if and only if 
$1 < r < s = k$.
Next we note that $(e_r+e_s) + (e_k-e_i) \not \in \Phi$ for all $r,s \neq 1,i$ and 
that $e_r + (e_k-e_i) \not \in \Phi$ for all $r \neq 1,i$. 
Thus (c) follows. 

\medskip
(d) and (e) follow from (a), and Lemma \ref{la:Cn_n_k}.
\epr

\medskip


\section{Single root midafis of exceptional groups}
\label{sec:exc_midafis}

In this section we deal with the case that the root system $\Phi$
is irreducible of type $E_6$, $E_7$, $E_8$, $F_4$ or $G_2$ and prove
Theorems~\ref{thm:main_exceptional} and \ref{thm:main_exceptionalnormal}.
We will use the explicit construction of these root systems given
in~\cite[Section~12.1]{Humphreys}. We assume the setting and notation
from Sections~\ref{sec:nota}-\ref{sec:single}. 

\medskip

\emph{We assume throughout this section that Hypothesis~\ref{hyp:p} holds.}

\smallskip


\subsection{\texorpdfstring{Types $E_6, E_7, E_8$}{Types E6, E7, E8}} 
\label{subsec:E6E7E8}

We deal with the types $E_6$ and $E_7$ by considering suitable root
subsystems of a root system of type $E_8$. The root system of
type~$E_8$ is constructed as follows (see~\cite[Section~12.1]{Humphreys}): 
Let $e_1, e_2, \dots, e_8 \in \R^8$ be the usual orthonormal unit
vectors which form a basis of $\R^8$ and let $\Phi_8$ be the union of
the sets $\{ \pm(e_i \pm e_j) \, | \, 1 \le i \neq j \le 8 \}$ and 
\[
\{\frac{1}{2}(\pm e_1 \pm e_2 \pm e_3 \pm e_4 \pm e_5 \pm e_6 \pm e_7 \pm e_8) \, | \, 
\mbox{the number of minus signs is even} \}.
\]
Then $\Phi_8$ is is a root system of type~$E_8$ and the set
$\{\alpha_1, \dots, \alpha_8\}$, where 
\begin{align*}
\begin{split}
\alpha_1 &:= \frac{1}{2}(e_1+e_8 - e_2 - e_3 - e_4 - e_5 - e_6 - e_7), \
\alpha_2 := e_1+e_2, \ \alpha_3 := e_2-e_1,\\
\alpha_4 &:= e_3-e_2, \ \alpha_5 := e_4-e_3, \ \alpha_6 := e_5-e_4,
\ \alpha_7 := e_6-e_5, \ \alpha_8 := e_7-e_6,
\end{split}
\end{align*}
is a set of simple roots. The corresponding set of positive roots is 
\[
\Phi_8^+ = \{ e_i \pm e_j \, | \, 1 \le j < i \le 8 \} \cup
 \{\frac{1}{2}(\pm e_1 \pm e_2 \pm e_3 \pm e_4 \pm e_5 \pm e_6 \pm e_7 +
 e_8) \}
\]
where the number of minus signs of the coefficients of the vectors in
the second set is even. 

The subsystem of $\Phi_8$ generated by 
$\{\a_1, \dots , \a_6\}$ is a root system of type~$E_6$ which 
we denote by $\Phi_6$, whereas the subsystem generated by 
$\{\a_1, \dots , \a_7 \}$ is a root system of type~$E_7$ 
which we denote by $\Phi_7$. For $i \in \{6,7\}$, the set 
$\{\a_1, \dots , \a_i\}$ is a set of simple roots for $\Phi_i$ 
and $\Phi_i^+ := \Phi_8^+ \cap \Phi_i$ is the corresponding 
set of positive roots. For $\alpha \in \Phi_i^+$ we always take 
$k(\alpha)$ with respect to $\Phi_8$ in this section.

We number the positive roots of $\Phi_8$ according to
Table~\ref{tab:rootse8}. This table contains the following
information: The first column fixes the notation for the positive
roots of $\Phi_8$. The second column lists the coefficients $m_j$ when
the root $\alpha_i = \sum_{j=1}^8 m_j \alpha_j$ is written as a linear
combination of the simple roots $\alpha_1, \dots, \alpha_8$. The third
column expresses the root $\alpha_i$ as a linear combination of the
vectors $e_1, \dots, e_8$ and the last column contains the height
$\height(\alpha_i)$. For example, the positive root $\alpha_{69}$ is
\begin{eqnarray*}
\alpha_{69} & = & 1 \cdot \alpha_1 + 2 \cdot \alpha_2 + 2 \cdot \alpha_3
+ 3 \cdot \alpha_4 + 2 \cdot \alpha_5 + 1 \cdot \alpha_6 + 0 \cdot
\alpha_7 + 0 \cdot \alpha_8 \\
& = & \frac12(e_1+e_2+e_3+e_4+e_5-e_6-e_7+e_8)
\end{eqnarray*}
and we have $\height(\alpha_{69}) = 11$. In particular, we have
\[
UE_6 = \prod_{\alpha \in \Phi_6^+} X_\alpha, \qquad 
UE_7 = \prod_{\alpha \in \Phi_7^+} X_\alpha \quad \text{and} \quad
UE_8 = \prod_{\alpha \in \Phi_8^+} X_\alpha.
\]

\vspace{-0.2cm}

Let $i \in \{6,7\}$. By factoring out the unipotent radical
$\prod_{\gamma \in \Phi_8^+ \setminus \Phi_i^+} X_\gamma$ we can
identify the group $UE_i$ with a factor group of $UE_8$ in a natural way and it
follows from Remark~\ref{rem:kSigma} that for each $\alpha \in \Phi_i^+$ 
we have $\Phi_8^+ \setminus \Phi_i^+ \subseteq k(\alpha)$.

To formulate the next result we introduce the following set of positive roots:
\begin{eqnarray*}
R^{\norml}_{6/7/8} & := & \Phi_8^+ \setminus \{\alpha_{45},
\alpha_{53}, \alpha_{57}, \alpha_{59}, \alpha_{60}, 
\alpha_{64}, \alpha_{67}, \alpha_{70}, \alpha_{71}, \alpha_{72},
\alpha_{73}, \alpha_{76}, \alpha_{77},\\
&& \alpha_{78}, \alpha_{80}, \alpha_{83}, \alpha_{84}, \alpha_{85},
\alpha_{86}, \alpha_{87}, \alpha_{88}, \alpha_{89}, \alpha_{90},
\alpha_{91}, \alpha_{92}, \alpha_{94}, \alpha_{95}, \alpha_{98},\\
&& \alpha_{99}, \alpha_{100}, \alpha_{102}, \alpha_{103},
\alpha_{104}, \alpha_{105}, \alpha_{106}, \alpha_{107}, \alpha_{108},
\alpha_{109}, \alpha_{110}, \alpha_{111},\\
&& \alpha_{113}, \alpha_{ 114}, \alpha_{115}, \alpha_{116},
\alpha_{117}, \alpha_{118}\}. 
\end{eqnarray*}
The next proposition includes a construction of all single root
midafis of the groups $UE_6$, $UE_7$, $UE_8$ and proves
Theorems~\ref{thm:main_exceptional} and \ref{thm:main_exceptionalnormal} 
for root systems $\Phi$ of type $E_6$, $E_7$ and $E_8$. The proof of
the proposition is based on computer calculations. These calculations
are carried out with the help of computer programs which we have
implemented in CHEVIE~\cite{CHEVIE}, \cite{CHEVIEdev}.

We proceed as follows: Let $i \in \{6,7,8\}$ and let $\alpha \in \Phi_i^+$. 
We choose an arm $a(\alpha) = \{\alpha_{i_1}, \alpha_{i_2}, \dots, \alpha_{i_r}\}$ 
of the hook $h(\alpha)$ such that the indices $i_1, i_2, \dots, i_r$ are 
given by the second column of Table~\ref{tab:armse8} in the appendix. 
The corresponding leg is $\ell(\alpha) = \{\alpha-\gamma \mid \gamma
\in a(\alpha)\}$ and the source is $s(\alpha) = \Phi_i^+ \setminus a(\alpha)$.
By computer calculations using CHEVIE we will show that $s(\alpha)$ is
a closed pattern. We know from Lemma~\ref{la:derquopat} that we can
identify the root subgroup~$X_\alpha$ with the quotient pattern group
$S_\alpha/P(s(\alpha) \setminus \{\alpha\})$ of
the source group $S_\alpha = P(s(\alpha))$.

Suppose that $\alpha \in \Phi_i^+ \setminus R^{\norml}_{6/7/8}$. Let
$\bar{\ell}(\alpha) := \ell(\alpha) \cup \{a_{j_1},
\alpha_{j_2},\dots, \alpha_{j_t}\}$ where the indices $j_1, j_2,
\dots, j_t$ are the maxima of the sets of integers given in the 
third column of Table~\ref{tab:armse8}. 
By computer calculations using CHEVIE we will show
that $\bar{\ell}(\alpha) \cup k(\alpha) \unlhd s(\alpha)$,
that the quotient pattern group 
$P(\bar{\ell}(\alpha) \cup k(\alpha))/P(k(\alpha))$
is abelian and that $\bar{\ell}(\alpha) \cup \{\alpha\} \cup k(\alpha) \unlhd \Phi_i^+$.
Hence, we obtain the quotient pattern group 
\[
\overline{T}_\alpha := S_\alpha/P(\{\alpha\} \cup \bar{\ell}(\alpha) \cup k(\alpha)).
\]

Now suppose that $\alpha \in \Phi_i^+ \cap R^{\norml}_{6/7/8}$. 
Again using CHEVIE we will show that $\ell(\alpha) \cup k(\alpha) \unlhd s(\alpha)$
and that $\ell(\alpha) \cup \{\alpha\} \cup k(\alpha) \unlhd \Phi_i^+$
so that we can consider the quotient pattern group 
\[
\overline{T}_\alpha := S_\alpha/P(\{\alpha\} \cup \ell(\alpha) \cup k(\alpha)).
\]

\begin{prop} \label{prop:single_mida_e678}
Let $i \in \{6,7,8\}$ and let $\Phi_i$ be a root system of type $E_i$
as described above. For each positive root $\alpha \in \Phi_i^+$ the
following are true:
\begin{enumerate}
\item[(a)] If $\alpha \in \Phi_i^+ \setminus R^{\norml}_{6/7/8}$ then 
  $\Psi_\alpha: \Irr^\lin(\overline{T}_\alpha) \times \Irr(X_\alpha)^*
  \to \Irr^\mida(UE_i)_\alpha$ with 
  \[
  (\mu, \lambda) \mapsto (\infl_{\overline{T}_\alpha}^{S_\alpha} \mu \cdot
   \infl_{X_\alpha}^{S_\alpha} \lambda)^{UE_i}
  \]
 is a one to one correspondence.

\item[(b)] If $\alpha \in R^{\norml}_{6/7/8}$ then 
  $\Psi_\alpha: \Irr(\overline{T}_\alpha) \times \Irr(X_\alpha)^*
  \to \Irr(UE_i)_\alpha$ with 
  \[
  (\mu, \lambda) \mapsto (\infl_{\overline{T}_\alpha}^{S_\alpha} \mu \cdot
   \infl_{X_\alpha}^{S_\alpha} \lambda)^{UE_i}
  \]
 is a one to one correspondence.
\end{enumerate} 
The number $|\Irr^\mida(UE_i)_\alpha|$ of midafis for $\alpha$ is
given in the second and the fifth column of Table~\ref{tab:mida_e8}
and the degree $\chi(1)$ for $\chi \in \Irr^\mida(UE_i)_\alpha$ is
given in the third and the sixth column of Table~\ref{tab:mida_e8}.   
\end{prop}

\begin{rem}
Let $i \in \{6,7,8\}$ and let $\alpha \in \Phi_i^+$ be a positive
root. By definition, all midafis $\chi \in \Irr^\mida(UE_i)_\alpha$
have the same degree; this is number given in the third and the 
sixth column of Table~\ref{tab:mida_e8}. For example: 
$|\Irr^\mida(UE_8)_{\alpha_{115}}| = q^8 (q-1)$ and each 
$\chi \in \Irr^\mida(UE_8)_{\alpha_{115}}$ has the degree $\chi(1) = q^{23}$.

The roots in $\Phi_6^+$, $\Phi_7^+ \setminus \Phi_6^+$,
$\Phi_8^+ \setminus \Phi_7^+$, respectively, are separated from 
each other by horizontal lines in Table~\ref{tab:mida_e8} and
Table~\ref{tab:armse8}. 
\end{rem}

\begin{center}
 \begin{longtable}{|c|l|l||c|l|l|}
 \caption{Numbers and degrees of the midafis for roots $\alpha_i \in \Phi_8^+$.} \label{tab:mida_e8} \\

 \hline
 \hspace{-0.2cm} Root \hspace{-0.2cm} & Number of midafis & Degree
 & \hspace{-0.2cm} Root \hspace{-0.2cm} & Number of midafis &
 Degree\\ \hline 
 \endfirsthead

 \multicolumn{6}{c}{{\tablename\ \thetable{} (cont.)}} \\
 \hline
 \hspace{-0.2cm} Root \hspace{-0.2cm} & Number of midafis & Degree
 & \hspace{-0.2cm} Root \hspace{-0.2cm} & Number of midafis &
 Degree\\ \hline 
 \endhead

 \hline 
 \endfoot

 \endlastfoot

$\alpha_{1}$ & 
$q-1$ & 
$1$&

$\alpha_{2}$ & 
$q-1$ & 
$1$\\

$\alpha_{3}$ & 
$q-1$ & 
$1$&

$\alpha_{4}$ & 
$q-1$ & 
$1$\\

$\alpha_{5}$ & 
$q-1$ & 
$1$&

$\alpha_{6}$ & 
$q-1$ & 
$1$\\

$\alpha_{9}$ & 
$q-1$ & 
$q$&

$\alpha_{10}$ & 
$q-1$ & 
$q$\\

$\alpha_{11}$ & 
$q-1$ & 
$q$&

$\alpha_{12}$ & 
$q-1$ & 
$q$\\

$\alpha_{13}$ & 
$q-1$ & 
$q$&

$\alpha_{16}$ & 
$q (q-1)$ & 
$q^2$\\

$\alpha_{17}$ & 
$q (q-1)$ & 
$q^2$&

$\alpha_{18}$ & 
$q (q-1)$ & 
$q^2$\\

$\alpha_{19}$ & 
$q (q-1)$ & 
$q^2$&

$\alpha_{20}$ & 
$q (q-1)$ & 
$q^2$\\

$\alpha_{23}$ & 
$q^2 (q-1)$ & 
$q^3$&

$\alpha_{24}$ & 
$q^2 (q-1)$ & 
$q^3$\\

$\alpha_{25}$ & 
$q^4 (q-1)$ & 
$q^3$&

$\alpha_{26}$ & 
$q^2 (q-1)$ & 
$q^3$\\

$\alpha_{27}$ & 
$q^2 (q-1)$ & 
$q^3$&

$\alpha_{30}$ & 
$q^5 (q-1)$ & 
$q^4$\\

$\alpha_{31}$ & 
$q^3 (q-1)$ & 
$q^4$&

$\alpha_{32}$ & 
$q^3 (q-1)$ & 
$q^4$\\

$\alpha_{33}$ & 
$q^5 (q-1)$ & 
$q^4$&

$\alpha_{37}$ & 
$q^5 (q-1)$ & 
$q^5$\\

$\alpha_{38}$ & 
$q^6 (q-1)$ & 
$q^5$&

$\alpha_{40}$ & 
$q^5 (q-1)$ & 
$q^5$\\

$\alpha_{44}$ & 
$q^4 (q-1)$ & 
$q^6$&

$\alpha_{45}$ & 
$q^7 (q-1)$ & 
$q^6$\\

$\alpha_{48}$ & 
$q^4 (q-1)$ & 
$q^6$&

$\alpha_{51}$ & 
$q^6 (q-1)$ & 
$q^7$\\

$\alpha_{52}$ & 
$q^6 (q-1)$ & 
$q^7$&

$\alpha_{57}$ & 
$q^6 (q-1)$ & 
$q^8$\\

$\alpha_{63}$ & 
$q^6 (q-1)$ & 
$q^9$&

$\alpha_{69}$ & 
$q^5 (q-1)$ & 
$q^{10}$\\
  \hline

$\alpha_{7}$ & 
$q-1$ & 
$1$ &

$\alpha_{14}$ & 
$q-1$ & 
$q$\\

$\alpha_{21}$ & 
$q (q-1)$ & 
$q^2$&

$\alpha_{28}$ & 
$q^2 (q-1)$ & 
$q^3$\\

$\alpha_{34}$ & 
$q^3 (q-1)$ & 
$q^4$&

$\alpha_{35}$ & 
$q^3 (q-1)$ & 
$q^4$\\

$\alpha_{39}$ & 
$q^4 (q-1)$ & 
$q^5$&

$\alpha_{41}$ & 
$q^6 (q-1)$ & 
$q^5$\\

$\alpha_{46}$ & 
$q^7 (q-1)$ & 
$q^6$&

$\alpha_{49}$ & 
$q^6 (q-1)$ & 
$q^6$\\

$\alpha_{53}$ & 
$q^8 (q-1)$ & 
$q^7$&

$\alpha_{55}$ & 
$q^6 (q-1)$ & 
$q^7$\\

$\alpha_{58}$ & 
$q^7 (q-1)$ & 
$q^8$&

$\alpha_{59}$ & 
$q^8 (q-1)$ & 
$q^8$\\

$\alpha_{61}$ & 
$q^5 (q-1)$ & 
$q^8$&

$\alpha_{64}$ & 
$q^8 (q-1)$ & 
$q^9$\\

$\alpha_{66}$ & 
$q^7 (q-1)$ & 
$q^9$&

$\alpha_{70}$ & 
$q^8 (q-1)$ & 
$q^{10}$\\

$\alpha_{71}$ & 
$q^7 (q-1)$ & 
$q^{10}$&

$\alpha_{75}$ & 
$q^7 (q-1)$ & 
$q^{11}$\\

$\alpha_{76}$ & 
$q^8 (q-1)$ & 
$q^{11}$&

$\alpha_{80}$ & 
$q^7 (q-1)$ & 
$q^{12}$\\

$\alpha_{82}$ & 
$q^7 (q-1)$ & 
$q^{12}$&

$\alpha_{85}$ & 
$q^7 (q-1)$ & 
$q^{13}$\\

$\alpha_{89}$ & 
$q^7 (q-1)$ & 
$q^{14}$&

$\alpha_{93}$ & 
$q^7 (q-1)$ & 
$q^{15}$\\

$\alpha_{97}$ & 
$q^6 (q-1)$ & 
$q^{16}$&&&\\
  \hline

$\alpha_{8}$ & $q-1$ & $1$ & 

$\alpha_{15}$ & 
$q-1$ & 
$q$\\

$\alpha_{22}$ & 
$q (q-1)$ & 
$q^2$&

$\alpha_{29}$ & 
$q^2 (q-1)$ & 
$q^3$\\
$\alpha_{36}$ & 
$q^3 (q-1)$ & 
$q^4$&

$\alpha_{42}$ & 
$q^4 (q-1)$ & 
$q^5$\\

$\alpha_{43}$ & 
$q^4 (q-1)$ & 
$q^5$&

$\alpha_{47}$ & 
$q^5 (q-1)$ & 
$q^6$\\

$\alpha_{50}$ & 
$q^7 (q-1)$ & 
$q^6$&

$\alpha_{54}$ & 
$q^8 (q-1)$ & 
$q^7$\\

$\alpha_{56}$ & 
$q^7 (q-1)$ & 
$q^7$&

$\alpha_{60}$ & 
$q^9 (q-1)$ & 
$q^8$\\

$\alpha_{62}$ & 
$q^7 (q-1)$ & 
$q^8$&

$\alpha_{65}$ & 
$q^8 (q-1)$ & 
$q^9$\\

$\alpha_{67}$ & 
$q^9 (q-1)$ & 
$q^9$&

$\alpha_{68}$ & 
$q^7 (q-1)$ & 
$q^9$\\

$\alpha_{72}$ & 
$q^9 (q-1)$ & 
$q^{10}$&

$\alpha_{73}$ & 
$q^9 (q-1)$ & 
$q^{10}$\\

$\alpha_{74}$ & 
$q^6 (q-1)$ & 
$q^{10}$&

$\alpha_{77}$ & 
$q^9 (q-1)$ & 
$q^{11}$\\

$\alpha_{78}$ & 
$q^9 (q-1)$ & 
$q^{11}$&

$\alpha_{79}$ & 
$q^8 (q-1)$ & 
$q^{11}$\\

$\alpha_{81}$ & 
$q^8 (q-1)$ & 
$q^{12}$&

$\alpha_{83}$ & 
$q^{10} (q-1)$ & 
$q^{12}$\\

$\alpha_{84}$ & 
$q^8 (q-1)$ & 
$q^{12}$&

$\alpha_{86}$ & 
$q^9 (q-1)$ & 
$q^{13}$\\

$\alpha_{87}$ & 
$q^9 (q-1)$ & 
$q^{13}$&

$\alpha_{88}$ & 
$q^9 (q-1)$ & 
$q^{13}$\\

$\alpha_{90}$ & 
$q^9 (q-1)$ & 
$q^{14}$&

$\alpha_{91}$ & 
$q^8 (q-1)$ & 
$q^{14}$\\

$\alpha_{92}$ & 
$q^9 (q-1)$ & 
$q^{14}$&

$\alpha_{94}$ & 
$q^9 (q-1)$ & 
$q^{15}$\\

$\alpha_{95}$ & 
$q^9 (q-1)$ & 
$q^{15}$&

$\alpha_{96}$ & 
$q^8 (q-1)$ & 
$q^{15}$\\

$\alpha_{98}$ & 
$q^9 (q-1)$ & 
$q^{16}$&

$\alpha_{99}$ & 
$q^9 (q-1)$ & 
$q^{16}$\\

$\alpha_{100}$ & 
$q^8 (q-1)$ & 
$q^{16}$&

$\alpha_{101}$ & 
$q^8 (q-1)$ & 
$q^{17}$\\

$\alpha_{102}$ & 
$q^9 (q-1)$ & 
$q^{17}$&

$\alpha_{103}$ & 
$q^9 (q-1)$ & 
$q^{17}$\\

$\alpha_{104}$ & 
$q^8 (q-1)$ & 
$q^{18}$&

$\alpha_{105}$ & 
$q^9 (q-1)$ & 
$q^{18}$\\

$\alpha_{106}$ & 
$q^8 (q-1)$ & 
$q^{18}$&

$\alpha_{107}$ & 
$q^8 (q-1)$ & 
$q^{19}$\\

$\alpha_{108}$ & 
$q^9 (q-1)$ & 
$q^{19}$&

$\alpha_{109}$ & 
$q^8 (q-1)$ & 
$q^{20}$\\

$\alpha_{110}$ & 
$q^9 (q-1)$ & 
$q^{20}$&

$\alpha_{111}$ & 
$q^9 (q-1)$ & 
$q^{21}$\\

$\alpha_{112}$ & 
$q^8 (q-1)$ & 
$q^{21}$&

$\alpha_{113}$ & 
$q^8 (q-1)$ & 
$q^{22}$\\

$\alpha_{114}$ & 
$q^8 (q-1)$ & 
$q^{22}$&

$\alpha_{115}$ & 
$q^8 (q-1)$ & 
$q^{23}$\\

$\alpha_{116}$ & 
$q^8 (q-1)$ & 
$q^{24}$&

$\alpha_{117}$ & 
$q^8 (q-1)$ & 
$q^{25}$\\

$\alpha_{118}$ & 
$q^8 (q-1)$ & 
$q^{26}$&

$\alpha_{119}$ & 
$q^8 (q-1)$ & 
$q^{27}$\\

$\alpha_{120}$ & 
$q^7 (q-1)$ & 
$q^{28}$&&&\\
  \hline
 \end{longtable}
\end{center}

\pr (of Proposition~\ref{prop:single_mida_e678})
The proof is carried out by computer programs which we have
implemented in CHEVIE. In particular, for all computations with roots
we use these CHEVIE programs.

\smallskip

(a) We demonstrate the proof only for $\alpha = \alpha_{115} \in \Phi_8^+ \setminus R^{\norml}_{6/7/8}$. 
The proof for the other roots in $\Phi_i^+ \setminus R^{\norml}_{6/7/8}$ is similar. 
Let $\alpha=\alpha_{115}$. Using the data in Table~\ref{tab:rootse8} we see
that the hook corresponding to $\alpha$ is
\begin{align*}
h(\alpha) &= \{\gamma \in \Phi_8^+ \mid \text{There is } \gamma' \in
\Phi_8^+ \text{ such that } \gamma+\gamma' = \alpha_{115}.\} \\
&= \{ \alpha_{2}, \alpha_{3}, \alpha_{9}, \alpha_{17}, \alpha_{23},
\alpha_{25}, \alpha_{30}, \alpha_{33}, \alpha_{38}, \alpha_{41},
\alpha_{44}, \alpha_{46}, \alpha_{50}, \alpha_{51}, \alpha_{54},\\
& \hspace{0.5cm} \alpha_{57}, \alpha_{58}, \alpha_{64}, \alpha_{65}, \alpha_{69},
\alpha_{71}, \alpha_{72}, \alpha_{75}, \alpha_{78}, \alpha_{80},
\alpha_{81}, \alpha_{84}, \alpha_{85}, \alpha_{86}, \alpha_{90},\\
& \hspace{0.5cm} \alpha_{91}, \alpha_{93}, \alpha_{95}, \alpha_{97}, \alpha_{98},
\alpha_{100}, \alpha_{101}, \alpha_{102}, \alpha_{104}, \alpha_{105},
\alpha_{107}, \alpha_{108}, \alpha_{109},\\
& \hspace{0.5cm} \alpha_{112}, \alpha_{113}, \alpha_{114}, \alpha_{115}\}.
\end{align*}
According to Table~\ref{tab:armse8} we make the following choice for
the arm and the leg of the hook $h(\alpha)$:
\begin{eqnarray*}
a(\alpha) & = & \{\alpha_{2}, \alpha_{3}, \alpha_{9}, \alpha_{17},
\alpha_{23}, \alpha_{25}, \alpha_{30}, \alpha_{33}, \alpha_{38},
\alpha_{41}, \alpha_{44}, \alpha_{46}, \alpha_{50}, \alpha_{51},
\alpha_{54}, \alpha_{57},\\
&& \alpha_{58}, \alpha_{64}, \alpha_{65}, \alpha_{71}, \alpha_{72},
\alpha_{78}, \alpha_{84}\}, \\ 
\ell(\alpha) & = & \{\alpha_{69}, \alpha_{75}, \alpha_{80},
\alpha_{81}, \alpha_{85}, \alpha_{86}, \alpha_{90}, \alpha_{91},
\alpha_{93}, \alpha_{95}, \alpha_{97}, \alpha_{98}, \alpha_{100},
\alpha_{101}, \alpha_{102},\\
&& \alpha_{104}, \alpha_{105}, \alpha_{107}, \alpha_{108},
\alpha_{109}, \alpha_{112}, \alpha_{113}, \alpha_{114}\}. 
\end{eqnarray*}
Using CHEVIE we verify that $s(\alpha) = \Phi_8^+ \setminus a(\alpha)$
is a closed pattern and using Definition~\ref{def:k} we get 
$k(\alpha) = \{\alpha_{116}, \alpha_{117}, \alpha_{118}, \alpha_{119}, \alpha_{120}\}$.
In particular, we have $\ell(\alpha) \cup k(\alpha) \subseteq s(\alpha)$. 
The normal closure of $\ell(\alpha) \cup k(\alpha)$ in $s(\alpha)$ is
$\bar{\ell}(\alpha) \cup k(\alpha)$ where 
\[
\bar{\ell}(\alpha) = \{\alpha_{89}, \alpha_{94}, \alpha_{99},
\alpha_{103}, \alpha_{106}, \alpha_{110}, \alpha_{111}\} \cup \ell(\alpha)
\]
and using CHEVIE we verify that 
$\bar{\ell}(\alpha) \cup \{\alpha\} \cup k(\alpha) \unlhd \Phi_8^+$. 
Also by a direct calculation we see that $h(\alpha) \cup k(\alpha)$ is
a closed pattern. We set
\[
S_\alpha := P(s(\alpha)) \quad \text{and} \quad 
H_\alpha := P(h(\alpha) \cup k(\alpha))/P(k(\alpha)).
\]
We claim that $H_\alpha$ is special of type 
$q^{1+2|a(\alpha)|} = q^{1+2 \cdot 23}$ and that $[x, H_\alpha]=Z(H_\alpha)=X_\alpha$  
for all $x \in H_\alpha \setminus Z(H_\alpha)$. The proof is analogous
to that of Lemma~\ref{la:An_hooks}~(d): A direct calculation shows
that for all $\gamma, \gamma' \in h(\alpha)$ with 
$\gamma+\gamma' \in \Phi_8^+$ we have $\gamma+\gamma' \in \{\alpha\}
\cup k(\alpha)$. It follows that $X_\alpha \subseteq Z(H_\alpha)$ and
$[H_\alpha, H_\alpha] \subseteq X_\alpha$. Now let 
$x \in H_\alpha \setminus X_\alpha$. We write 
$x=\prod_{\gamma \in h(\alpha)} x_\gamma(t_\gamma)$ as in~(\ref{eq:elemu}). 
Because $x \not\in X_\alpha$ there is some 
$\gamma \in h(\alpha) \setminus \{\alpha\}$ such that $t_\gamma \neq 0$.
Thus, $\gamma' := \alpha-\gamma \in h(\alpha)$ and we get from
Lemma~\ref{la:patnorm}~(a) that 
$\{[x,x_{\gamma'}(t)] \, | \, t \in \F_q\} = X_\alpha = Z(H_\alpha)
= [x, X_\alpha] = [H_\alpha, H_\alpha] = \Phi(H_\alpha)$. So 
$H_\alpha$ is special of type $q^{1+2|a(\alpha)|} = q^{1+2 \cdot 23}$.
Note that this argument also shows that $a(\alpha) \cup k(\alpha)$ 
is a closed pattern and that the quotient pattern group
$P(a(\alpha) \cup k(\alpha))/P(k(\alpha))$ is an abelian group.

To apply the Reduction Lemma~\ref{nred} and Lemma~\ref{la:twohook} we
introduce the following notation (only for this proof):
\begin{itemize}
\item $U := UE_8/P(k(\alpha))$,

\item $H := \overline{S}_\alpha := S_\alpha/P(k(\alpha)) \subseteq U$,

\item $Z := X_\alpha \subseteq U$,

\item $X := \prod_{\gamma \in a(\alpha)} X_\gamma \subseteq U$,

\item $Y := P(\bar{\ell}(\alpha) \cup k(\alpha))/P(k(\alpha))
  \subseteq U$,

\item $\Irr(U)_\alpha := \{\chi \in \Irr(U) \mid X_\alpha
  \not\subseteq \Ker(\chi)\}$,
\end{itemize}
where $X_\gamma$ denotes the image of $X_\gamma$ in $U$.

We have just seen that $X$ is an abelian group and it follows from
\eqref{eq:elemu} that $X$ is a set of representatives for $U/H$. 
A direct calculation shows that for all 
$\gamma, \gamma' \in \bar{\ell}(\alpha)$ such that 
$\gamma+\gamma' \in \Phi_8^+$ we have $\gamma+\gamma' \in k(\alpha)$. 
Thus $Y$ is abelian. By Lemma~\ref{la:alpha_rep} and
Definition~\ref{def:kSigma_almfaith} (a) we have 
$Z \subseteq Z(U)$. We have $Y \unlhd H$ and $ZY \unlhd U$ 
because $\bar{\ell}(\alpha) \cup k(\alpha) \unlhd s(\alpha)$ and
$\bar{\ell}(\alpha) \cup \{\alpha\} \cup k(\alpha) \unlhd \Phi_8^+$.
Furthermore we have $Y \cap Z = \{1\}$ because $\alpha \not\in \bar{\ell}(\alpha)$.
Hence, the conditions (a)-(d) in Lemma~\ref{nred} are satisfied. 
Suppose that $\lambda \in \Irr(Z)^*$ and let $\tilde{\lambda}$ be 
the inflation of $\lambda$ to $ZY = Z \times Y$. 
Suppose that there is $x \in X \setminus \{1\}$ such that  
${}^x\widetilde{\lambda} = \widetilde{\lambda}$. Hence we have
$\tilde{\lambda}(u^x)=\tilde{\lambda}(u)$ for all $u \in ZY$.
In particular, $\tilde{\lambda}(u^x)=\tilde{\lambda}(u)$ for
all $u \in L_\alpha := P(\ell(\alpha) \cup k(\alpha))/P(k(\alpha))$
and therefore $\tilde{\lambda}([u,x]) = 1$ for all $u \in L_\alpha$. 
Hence $x$ commutes with $L_\alpha$ modulo $\Ker(\tilde{\lambda})$. 
Since $X$ is abelian and $\langle X, L_\alpha\rangle = H_\alpha$ it
follows that $x$ commutes with $H_\alpha$ modulo $\Ker(\tilde{\lambda})$. 
Thus $1=\tilde{\lambda}([x,H_\alpha])=\tilde{\lambda}(X_\alpha)=
\tilde{\lambda}(Z)=\lambda(Z)$ which is impossible since $\lambda$ is
nontrivial. Hence, also condition~(e) of Lemma~\ref{nred}
is satisfied.

It follows from Lemmas~\ref{la:derquopat} and \ref{la:representable} 
that $H/Y \cong \overline{T}_\alpha \times X_\alpha$ where 
\[
\overline{T}_\alpha := S_\alpha/P(\{\alpha\} \cup \bar{\ell}(\alpha)
\cup k(\alpha)),
\]
and the Reduction Lemma~\ref{nred} gives a one to one
correspondence
\[
\tilde{\Psi}_\alpha: \Irr(\overline{T}_\alpha) \times \Irr(X_\alpha)^* \to
\Irr(U)_\alpha \cap \Irr(U, \mathbf{1}_Y), \quad 
(\mu, \lambda) \mapsto (\infl_{\overline{T}_\alpha}^H \mu \cdot
   \infl_{X_\alpha}^H \lambda)^U.
\]
We are interested in the characters $\chi \in \Irr(U)_\alpha$ such
that the degree $\chi(1)$ is minimal. Let $d := \min\{\chi(1) \mid
\chi \in \Irr(U)_\alpha\}$ and $\Irr^\mida(U)_\alpha := \{\chi \in
\Irr(U)_\alpha \mid \chi(1) = d\}$. For every $\lambda \in \Irr(X_\alpha)^*$ we have
\begin{equation} \label{eq:degq23}
\tilde{\Psi}_\alpha(\mu,\lambda)(1) \begin{cases}= [U:H] = q^{|a(\alpha)|} = q^{23} &
  \text{if } \mu \text{ is linear},\\
  >q^{23} & \text{if } \mu \text{ is nonlinear}.\end{cases}
\end{equation}
Let $\chi \in \Irr(U)_\alpha \setminus \Irr(U, \mathbf{1}_Y)$.
We claim that $\chi(1)>q^{23}$. There are exactly $q-1$ irreducible
characters of $H_\alpha$ not having $X_\alpha$ in their kernel, namely
the characters $(\infl_{X_\alpha}^{X_\alpha \times L_\alpha} \lambda)^{H_\alpha}$
where $\lambda \in \Irr(X_\alpha)^*$. Each such character has
degree $q^{|a(\alpha)|} = q^{23}$. Since $X_\alpha \not\subseteq \Ker(\chi)$ 
it follows that $\chi \in \Irr(U, \mathbf{1}_{L_\alpha})$. Thus the
restriction $\chi|_{ZY}$ has an irreducible constituent 
$\tilde{\mu} \in \Irr(ZY, \mathbf{1}_{L_\alpha})$. Note that
$\tilde{\mu}$ is a linear character since $Y$ and $ZY=Z \times Y$ 
are abelian. Since $\chi \in \Irr(U)_\alpha$, $ZY \unlhd U$ and 
$Z \subseteq Z(U)$ we have $X_\alpha \not\subseteq \Ker(\tilde{\mu})$
and since $\chi \not\in \Irr(U, \mathbf{1}_Y)$ we have 
$\tilde{\mu}|_Y \neq \mathbf{1}_Y$. It follows that there is some root 
\[
\beta \in \bar{\ell}(\alpha) \setminus \ell(\alpha) = \{\alpha_{89},
\alpha_{94}, \alpha_{99}, \alpha_{103}, \alpha_{106}, \alpha_{110},
\alpha_{111}\}
\]
such that $X_\beta \not\subseteq \Ker(\tilde{\mu})$. We only
demonstrate the case $\beta=\alpha_{106}$. The other cases are
similar.

We choose the subhook $h'(\beta) = \{\alpha_5, \alpha_{103}, \beta=\alpha_{106}\}$ 
according to Table~\ref{tab:armse8}. A direct calculation shows that 
$h'(\beta) \cup k(\alpha)$ is a closed pattern so that condition (a)
of Lemma~\ref{la:twohook} is satisfied. The fact $\alpha_5+\alpha_{103} = \beta$ 
implies that the group $H'_\beta := P(h'(\beta) \cup k(\alpha))/P(k(\alpha))$ 
satisfies condition (e) of Lemma~\ref{la:twohook}. We have 
already seen above that $H_\alpha$ satisfies condition (d) of
Lemma~\ref{la:twohook} and also condition~(c) is obviously satisfied.
As in Lemma~\ref{la:twohook} let $h_{\alpha\beta}'$ be the closed
pattern generated by $h(\alpha) \cup h'(\beta) \cup k(\alpha)$ so that
\begin{eqnarray*}
h_{\alpha\beta}' & = & \{\alpha_{2}, \alpha_{3}, \alpha_{5},
\alpha_{9}, \alpha_{17}, \alpha_{23}, \alpha_{25}, \alpha_{30},
\alpha_{33}, \alpha_{38}, \alpha_{41}, \alpha_{44}, \alpha_{46},
\alpha_{50}, \alpha_{51}, \alpha_{54},\\
&& \alpha_{57}, \alpha_{58}, \alpha_{64}, \alpha_{65}, \alpha_{69},
\alpha_{71}, \alpha_{72}, \alpha_{75}, \alpha_{78}, \alpha_{80},
\alpha_{81}, \alpha_{84}, \alpha_{85}, \alpha_{86}, \alpha_{90},\\
&& \alpha_{91}, \alpha_{93}, \alpha_{95}, \alpha_{97}, \alpha_{98},
\alpha_{100}, \alpha_{101}, \alpha_{102}, \alpha_{103}, \alpha_{104},
\alpha_{105}, \alpha_{106}, \alpha_{107},\\
&& \alpha_{108}, \alpha_{109}, \alpha_{112}, \alpha_{113},
\alpha_{114}, \alpha_{115}, \alpha_{116}, \alpha_{117}, \alpha_{118},
\alpha_{119}, \alpha_{120}\}. 
\end{eqnarray*}
A direct calculation shows that $h_{\alpha\beta}' \setminus a(\alpha)$ 
normalizes $\ell(\alpha) \cup k(\alpha)$ so that also condition (b) of
Lemma~\ref{la:twohook} is satisfied. Define
$\ell'_{\alpha\beta} := \{\alpha, \beta\} \cup \ell(\alpha) \cup k(\alpha)$
and $L'_{\alpha\beta} := P(\ell'_{\alpha\beta})/P(k(\alpha))$ as in Lemma~\ref{la:twohook}.
Let $\mu := \tilde{\mu}|_{L'_{\alpha\beta}}$. Since $\tilde{\mu}$ is a
linear character we have $\mu \in \Irr(L'_{\alpha\beta})$ and the 
properties of $\tilde{\mu}$ imply that $\mu$ satisfies the
assumptions in Lemma~\ref{la:twohook}. Because $\mu$ is
a constituent of $\chi|_{L'_{\alpha\beta}}$ there is a constituent
$\psi \in \Irr(H'_{\alpha\beta}, \mu)$ of $\chi|_{H'_{\alpha\beta}}$
and Lemma~\ref{la:twohook} gives us 
$\chi(1) \ge \psi(1) \ge q^{|a(\alpha)|+1} > q^{23}$. In particular,
we get that the degree $\chi(1)$ is not minimal among the degrees of
the irreducible characters in $\Irr(U)_\alpha$.

Using a similar argument for the other roots 
$\beta \in \bar{\ell}(\alpha) \setminus \ell(\alpha)$ we get that
$\chi(1)>q^{23}$ for all $\chi \in \Irr(U)_\alpha \setminus \Irr(U, \mathbf{1}_Y)$.
It follows that $d=q^{23}$ and
\[
\Irr^\mida(U)_\alpha = \{\chi \in \Irr(U)_\alpha \mid \chi(1) =
q^{23}\} \subseteq \Irr(U)_\alpha \cap \Irr(U, \mathbf{1}_Y). 
\]
Hence $\tilde{\Psi}_\alpha$ maps $\Irr^\lin(\overline{T}_\alpha) \times
\Irr(X_\alpha)^*$ one to one onto the set of irreducible characters in
$\Irr(U)_\alpha$ of minimal degree. Identifying the irreducible
characters of~$U$ with their inflations to $UE_8$ we see that the map
$\Psi_\alpha$ defined in part (a) of the proposition maps 
$\Irr^\lin(\overline{T}_\alpha) \times \Irr(\F_q)^*$ one to one onto
$\Irr^\mida(UE_8)_\alpha$. This completes the proof of part (a) of the
proposition.  

\medskip

(b) We demonstrate the proof only for $\alpha = \alpha_{112} \in \Phi_8^+ \cap R^{\norml}_{6/7/8}$. 
The proof for the other roots in $\Phi_i^+ \cap R^{\norml}_{6/7/8}$ is similar. 
Let $\alpha=\alpha_{112}$. Using the data in Table~\ref{tab:rootse8} we see
that the hook corresponding to $\alpha$ is
\begin{eqnarray*}
h(\alpha) & = & \{\gamma \in \Phi_8^+ \mid \text{There is } \gamma' \in
\Phi_8^+ \text{ such that } \gamma+\gamma' = \alpha_{112}.\} \\
& = & \{\alpha_{2}, \alpha_{10}, \alpha_{17}, \alpha_{18},
\alpha_{25}, \alpha_{26}, \alpha_{32}, \alpha_{33}, \alpha_{34},
\alpha_{40}, \alpha_{41}, \alpha_{42}, \alpha_{48}, \alpha_{49},
\alpha_{50}, \alpha_{55},\\ 
&& \alpha_{56}, \alpha_{61}, \alpha_{62}, \alpha_{68}, \alpha_{69},
\alpha_{74}, \alpha_{75}, \alpha_{80}, \alpha_{81}, \alpha_{85},
\alpha_{86}, \alpha_{89}, \alpha_{90}, \alpha_{91}, \alpha_{93},
\alpha_{94},\\
&& \alpha_{95}, \alpha_{98}, \alpha_{99}, \alpha_{100}, \alpha_{102},
\alpha_{103}, \alpha_{105}, \alpha_{106}, \alpha_{108}, \alpha_{110},
\alpha_{112}\}.
\end{eqnarray*}
According to Table~\ref{tab:armse8} we make the following choice for
the arm and the leg of the hook $h(\alpha)$:
\begin{eqnarray*}
a(\alpha) & = & \{\alpha_{2}, \alpha_{10}, \alpha_{17}, \alpha_{18},
\alpha_{25}, \alpha_{26}, \alpha_{32}, \alpha_{33}, \alpha_{34},
\alpha_{40}, \alpha_{41}, \alpha_{42}, \alpha_{48}, \alpha_{49},
\alpha_{50},\\
&& \alpha_{55}, \alpha_{56}, \alpha_{61}, \alpha_{62},
\alpha_{68}, \alpha_{74}\}, \\  
\ell(\alpha) & = & \{\alpha_{69}, \alpha_{75}, \alpha_{80},
\alpha_{81}, \alpha_{85}, \alpha_{86}, \alpha_{89}, \alpha_{90},
\alpha_{91}, \alpha_{93}, \alpha_{94}, \alpha_{95}, \alpha_{98},
\alpha_{99}, \alpha_{100},\\
&& \alpha_{102}, \alpha_{103}, \alpha_{105}, \alpha_{106}, \alpha_{108}, \alpha_{110}\}.
\end{eqnarray*}
Using CHEVIE we verify that $s(\alpha) = \Phi_8^+ \setminus a(\alpha)$
is a closed pattern and using Definition~\ref{def:k} we get 
\begin{eqnarray*}
k(\alpha) & = & \{\alpha_{97}, \alpha_{101}, \alpha_{104}, \alpha_{107},
\alpha_{109}, \alpha_{111}, \alpha_{113}, \alpha_{114}, \alpha_{115},
\alpha_{116}, \alpha_{117}, \alpha_{118},\\
&& \alpha_{119}, \alpha_{120}\}. 
\end{eqnarray*}
In particular, $\ell(\alpha) \cup k(\alpha) \subseteq s(\alpha)$. 
Using CHEVIE we verify that $\ell(\alpha) \cup k(\alpha) \unlhd s(\alpha)$,
that $\ell(\alpha) \cup \{\alpha\} \cup k(\alpha) \unlhd \Phi_8^+$, that
for all $\gamma, \gamma' \in \ell(\alpha)$ with
$\gamma+\gamma' \in \Phi_8^+$ we have $\gamma+\gamma' \in k(\alpha)$
and that for all $\gamma \in a(\alpha)$ and $\gamma' \in \ell(\alpha)$ we
have $\gamma+\gamma' \neq \alpha-\gamma$. We set $S_\alpha := P(s(\alpha))$. 
To apply the Reduction Lemma~\ref{nred} we introduce the following
notation (only for this proof):
\begin{itemize}
\item $U := UE_8/P(k(\alpha))$,

\item $H := \overline{S}_\alpha := S_\alpha/P(k(\alpha)) \subseteq U$,

\item $Z := X_\alpha \subseteq U$,

\item $X := \prod_{\gamma \in a(\alpha)} X_\gamma \subseteq U$,

\item $Y := P(\ell(\alpha) \cup k(\alpha))/P(k(\alpha))
  \subseteq U$,

\item $\Irr(U)_\alpha := \{\chi \in \Irr(U) \mid X_\alpha
  \not\subseteq \Ker(\chi)\}$,
\end{itemize}
where $X_\gamma$ denotes the image of $X_\gamma$ in $U$.

It follows from \eqref{eq:elemu} that $X$ is a set of representatives
for $U/H$. 
By Lemma~\ref{la:alpha_rep} and Definition~\ref{def:kSigma_almfaith} (a)
we have $Z \subseteq Z(U)$. We have $Y \unlhd H$ and $ZY \unlhd U$ 
because $\ell(\alpha) \cup k(\alpha) \unlhd s(\alpha)$ and
$\ell(\alpha) \cup \{\alpha\} \cup k(\alpha) \unlhd \Phi_8^+$.
Furthermore we have $Y \cap Z = \{1\}$ because $\alpha \not\in \ell(\alpha)$.
Hence, the conditions (a)-(d) in Lemma~\ref{nred} are satisfied. 

Suppose that $\lambda \in \Irr(Z)^*$. Let $\tilde{\lambda}$ be 
the inflation of $\lambda$ to $ZY = Z \times Y$. Note
that~$\tilde{\lambda}$ is a linear character. Suppose that there is 
$x \in X \setminus \{1\}$ such that ${}^x\widetilde{\lambda} = \widetilde{\lambda}$.
Because~$\lambda$ is nontrivial on $X_\alpha$ there is $t \in \F_q$ 
such that $\widetilde{\lambda}(x_\alpha(t)) = \lambda(x_\alpha(t)) \neq 1$.
It follows from Lemma~\ref{la:patnorm}~(a) that there is some 
root $\gamma \in \ell(\alpha)$ and an element $t' \in \F_q$ such that 
$[x, x_\gamma(t')^{-1}] = x_\alpha(t)$. Hence
\[
{}^x\widetilde{\lambda}(x_\gamma(t')) = 
\widetilde{\lambda}(x^{-1} \cdot x_\gamma(t') \cdot x \cdot
x_\gamma(t')^{-1} \cdot x_\gamma(t')) =
\widetilde{\lambda}( x_\alpha(t)) \widetilde{\lambda}( x_\gamma(t'))
\neq \widetilde{\lambda}(x_\gamma(t')),
\]
contradicting ${}^x\widetilde{\lambda} = \widetilde{\lambda}$. 
Hence, the condition~(e) of Lemma~\ref{nred} holds. The
condition~(f) of Lemma~\ref{nred} is also satisfied because 
$|X| = q^{|a(\alpha)|} = |Y|$. It follows from
Lemmas~\ref{la:derquopat} and \ref{la:representable}  
that $H/Y \cong \overline{T}_\alpha \times X_\alpha$ where 
\[
\overline{T}_\alpha := S_\alpha/P(\{\alpha\} \cup \bar{\ell}(\alpha) \cup k(\alpha)).
\]
Now the Reduction Lemma~\ref{nred} gives the one to one
correspondence
\[
\tilde{\Psi}_\alpha: \Irr(\overline{T}_\alpha) \times \Irr(X_\alpha)^* \to
\Irr(U)_\alpha, \quad (\mu, \lambda) \mapsto (\infl_{\overline{T}_\alpha}^H \mu \cdot
   \infl_{X_\alpha}^H \lambda)^U.
\]
Identifying the irreducible characters of~$U$ with their inflations to
$UE_8$ gives the one to one correspondence $\Psi_\alpha$ defined in
part (b) of the proposition. Note that the bijection $\Psi_\alpha$ 
maps $\Irr^\lin(\overline{T}_\alpha) \times \Irr(X_\alpha)^*$ onto
$\Irr^\mida(UE_i)_\alpha$. This completes the proof of part~(b) 
of the proposition.  

\medskip

Using Lemma~\ref{la:derquopat} we can easily compute
$|\Irr^\lin(\overline{T}_\alpha)| = 
|\overline{T}_\alpha/[\overline{T}_\alpha,\overline{T}_\alpha]|$
and get
\begin{eqnarray*}
|\Irr^\mida(UE_8)_\alpha| & = &
|\Psi_\alpha(\Irr^\lin(\overline{T}_\alpha) \times \Irr(X_\alpha)^*)|= 
|\Irr^\lin(\overline{T}_\alpha)| \cdot |\Irr(\F_q)^*|\\
& = & |\Irr^\lin(\overline{T}_\alpha)| \cdot (q-1).
\end{eqnarray*}
For every $\chi \in \Irr^\mida(UE_8)_\alpha$ we have 
$\chi(1) = \Psi_\alpha(\mu,\lambda)(1)$ for some linear character
$\mu$ of $\overline{T}_\alpha$ and $\lambda \in \Irr(X_\alpha)^*$ and hence
\[
\chi(1) = \Psi_\alpha(\mu,\lambda)(1) = [UE_8:S_\alpha] \cdot \mu(1)
\cdot \lambda(1) = q^{|a(\alpha)|}.
\]
This gives the entries in Table~\ref{tab:mida_e8} and completes the
proof of Proposition~\ref{prop:single_mida_e678}.
\epr

\begin{rem} \label{rem:choiceandstruct}
Let $i \in \{6,7,8\}$ and let $\alpha \in \Phi_i^+$ be a positive
root. 
\begin{enumerate}
\item[(a)] We could not find a canonical choice for the arm $a(\alpha)$. 
  However, in some sense the choices in Table~\ref{tab:armse8}
  are best possible; see Remark~\ref{rem:choicearm}. 

\item[(b)] For roots $\alpha \in \Phi_i^+ \cap R^\norml_{6/7/8}$
  Proposition~\ref{prop:single_mida_e678} (b) reduces the 
  classification of the midafis in $\Irr(UE_i)_\alpha$ to the
  character theory of the subquotient $\overline{T}_\alpha$. We
  illustrate by the example 
  $\alpha = \alpha_{112} = \frac{1}{2}(e_1 + e_2 + e_3 + e_4 + e_5 + e_6 + e_7 +e_8)$
  for $i=8$ how information on the structure of $\overline{T}_\alpha$
  can be obtained. The following are true:
  \begin{enumerate}
  \item[(1)] $n(\alpha) =  \{e_8 + e_j \mid 1 \leq j \leq
    7\} = \{\alpha_{113}, \alpha_{115}, \alpha_{116}, \alpha_{117},
    \alpha_{118}, \alpha_{119}, \alpha_{120}\}$,
  \item[(2)] $ k(\alpha) = \{e_8 \pm e_j \mid 1 \leq j \leq 7\} = 
    \{\alpha_{97}, \alpha_{101}, \alpha_{104}, \alpha_{107},
    \alpha_{109}, \alpha_{111}, \alpha_{113},$\newline
    \hspace*{1.2cm} $\alpha_{114}, \alpha_{115}, \alpha_{116}, \alpha_{117}, \alpha_{118},
    \alpha_{119}, \alpha_{120}\}$,
  \item[(3)] $a(\alpha) = \{ e_i + e_j \mid 1 \leq j < i \leq 7 \} =
  \{\alpha_{2}, \alpha_{10}, \alpha_{17}, \alpha_{18}, \alpha_{25},
  \alpha_{26}, \alpha_{32}, \alpha_{33},$\newline
  \hspace*{1.2cm} $\alpha_{34}, \alpha_{40}, \alpha_{41}, \alpha_{42},
  \alpha_{48}, \alpha_{49}, \alpha_{50}, \alpha_{55}, \alpha_{56},
  \alpha_{61}, \alpha_{62}, \alpha_{68}, \alpha_{74}\}$,
  \item[(4)] $\ell(\alpha) = \{\frac{1}{2}(\pm e_1 \pm e_2 \pm e_3 \pm e_4
    \pm e_5 \pm e_6 \pm e_7 + e_8) \mid \text{exactly two signs}$\newline
  \hspace*{1.2cm} $\text{are negative}\} = \{\alpha_{69}, \alpha_{75},
  \alpha_{80}, \alpha_{81}, \alpha_{85}, \alpha_{86}, \alpha_{89},
  \alpha_{90}, \alpha_{91}, \alpha_{93},$\newline
  \hspace*{1.2cm} $\alpha_{94}, \alpha_{95}, \alpha_{98}, \alpha_{99},
  \alpha_{100}, \alpha_{102}, \alpha_{103}, \alpha_{105},
  \alpha_{106}, \alpha_{108}, \alpha_{110}\}$.
  \end{enumerate}
  Equations (3) and (4) are the result of the calculation recorded
  in Table~\ref{tab:armse8}. The other equations can be confirmed by
  our CHEVIE programs, but they can also be proved by hand: (1)
  follows from the equation  
  \begin{eqnarray*}
  e_8 \pm e_j & = & \frac{1}{2}(\pm e_1 \pm \dots \pm e_{j-1} \pm e_j
  \pm e_{j+1} \dots \pm e_7 + e_8) \\
  && + \frac{1}{2}(\mp e_1 \mp \dots \mp e_{j-1} \pm e_j \mp e_{j+1} \dots
  \mp e_7 + e_8).
  \end{eqnarray*}
  For (2) we observe that $Z(UE_8/P(n(\alpha))) = X_\alpha, X_{e_8-e_1}$.
  Factoring out $X_{e_8-e_1}$ and repeating yields the claim after
  seven iterations.

  We now define the following subsets of $\Phi_8^+$: 
  \begin{eqnarray*}
  a_6 & := & \{ e_i - e_j \mid 1 \leq j < i \leq 7\} = \{\alpha_{3},
  \alpha_{4}, \alpha_{5}, \alpha_{6}, \alpha_{7}, \alpha_{8},
  \alpha_{11}, \alpha_{12}, \alpha_{13}, \alpha_{14},\\
  && \ \ \alpha_{15}, \alpha_{19}, \alpha_{20}, \alpha_{21},
  \alpha_{22}, \alpha_{27}, \alpha_{28}, \alpha_{29}, \alpha_{35},
  \alpha_{36}, \alpha_{43}\}, \\
  v & := & \{\frac{1}{2}(\pm e_1 \pm e_2 \pm e_3 \pm e_4 \pm e_5 \pm
  e_6 \pm e_7 + e_8) \mid \text{exactly two} + \text{signs}\}\\
  & = & \{\alpha_{1}, \alpha_{9}, \alpha_{16}, \alpha_{24},
  \alpha_{31}, \alpha_{39}, \alpha_{47}\},\\
  f & := & \{\frac{1}{2}(\pm e_1 \pm e_2 \pm e_3 \pm e_4 \pm e_5 \pm
  e_6 \pm e_7 + e_8) \mid \text{exactly four} + \text{signs}\}\\ 
  & = & \{\alpha_{23}, \alpha_{30}, \alpha_{37}, \alpha_{38},
  \alpha_{44}, \alpha_{45}, \alpha_{46}, \alpha_{51}, \alpha_{52},
  \alpha_{53}, \alpha_{54}, \alpha_{57}, \alpha_{58}, \alpha_{59}, \alpha_{60},\\
  && \ \ \alpha_{63}, \alpha_{64}, \alpha_{65}, \alpha_{66}, \alpha_{67},
  \alpha_{70}, \alpha_{71}, alpha_{72}, \alpha_{73}, \alpha_{76},
  \alpha_{77}, \alpha_{78}, \alpha_{79}, \alpha_{82},\\
  && \ \ \alpha_{83}, \alpha_{84}, \alpha_{87}, \alpha_{88}, \alpha_{92}, \alpha_{96}\}.
  \end{eqnarray*}
  With this notation we observe that $s(\alpha) = k(\alpha) \cup
  \{\alpha\} \cup a_6 \cup v \cup f \cup \ell(\alpha)$. Also let 
  $K$ denote the standard $D_7$-parabolic subgroup of $E_8(q)$ with
  Levi decomposition $K = Q \rtimes D$ parabolic subgroup (with respect
  to our choice of the root datum). Then 
  \[
  Z(Q) = P(k(\alpha)), \ \ Q/Z(Q) = P(\{\alpha\} \cup v \cup f \cup
  \ell(\alpha) \cup k(\alpha))/P(k(\alpha)) 
  \]
  and $\langle X_\gamma \mid \gamma \in a(\alpha) \cup a_6 \cup (-a_6)\rangle$ 
  is a maximal parabolic subgroup of $D$ with $A_6$-Levi factor 
  $A := \langle X_\gamma \mid \gamma \in a_6 \cup (-a_6)\rangle$ and unipotent radical 
  $P(a(\alpha))$. We remark that $Z(Q)$ is the natural module and $Q/Z(Q)$ is a
  half-spin module for $D$. Hence both are elementary abelian.  
  With these facts at our disposal we can now observe:
  \begin{enumerate}
  \item[(1)] $H_\alpha = P(h(\alpha))$ is a special group which is normalized by $P(a_6)$,
  \item[(2)] $\overline{P(v)} := P(v \cup k(\alpha))/P(k(\alpha))$
    and $\overline{P(f)} := P(f \cup k(\alpha))/P(k(\alpha))$ are
    elementary abelian and centralize each other,
  \item[(3)] $P(\ell(\alpha))$ is normalized by $P(a_6)$ and
    centralized by $\overline{P(v)} \times \overline{P(f)}$,
  \item[(4)] $A$ normalizes $\overline{P(v)}$ and $\overline{P(f)}$,
  \item[(5)] $\overline{P(v)}$ is the natural module for $A$ and $\overline{P(f)}$
    is the alternating cube of the natural module,
  \item[(6)] $\overline{T}_\alpha = (\overline{P(v)} \times \overline{P(f)}) \rtimes P(a_6)$,
  \item[(7)] $P(v \cup a_6) \cong UA_7(q)$.
  \end{enumerate}
\end{enumerate}
\end{rem}


\subsection{\texorpdfstring{Type $F_4$}{Type F4}} 
\label{subsec:F4}

We construct a root system of type~$F_4$ as in~\cite[Section~12.1]{Humphreys}: 
Let $e_1, e_2, e_3, e_4 \in \R^4$ be the usual orthonormal unit
vectors which form a basis of $\R^4$. Then
$\Phi_4 := \{ \pm(e_i \pm e_j) \, | \, 1 \le i \neq j \le 4 \} \cup \{
\pm e_i \, | \, 1 \le i \le 4\} \cup \{\frac{1}{2}(\pm e_1 \pm e_2 \pm
e_3 \pm e_4) \}$ is a root system of type~$F_4$ and the set 
$\{\alpha_1, \dots, \alpha_4\}$, where $\alpha_1 := e_2-e_3$,
$\alpha_2 := e_3-e_4$, $\alpha_3 := e_4$, and
$\alpha_4:=\frac{1}{2}(e_1 - e_2 -e_3 - e_4)$, is a set of simple
roots. The corresponding set of positive roots is  
\[
\Phi_4^+ = \{ e_i \pm e_j \, | \, 1 \le i < j \le 4 \} \cup \{ e_i \, |
\, 1 \le i \le 4 \} \cup \{\frac{1}{2}(e_1 \pm e_2 \pm e_3 \pm e_4)
\}.  
\]
We number the positive roots of $\Phi_4$ according to
Table~\ref{tab:rootsf4}. This table contains the following
information: The first column fixes the notation for the positive
roots of~$\Phi_4$. The second column lists the coefficients $m_j$ when
the root $\alpha_i = \sum_{j=1}^4 m_j \alpha_j$ is written as a linear
combination of the simple roots $\alpha_1, \alpha_2, \alpha_3, \alpha_4$. 
The third column expresses the root $\alpha_i$ as a linear combination
of the vectors $e_1, e_2, e_3, e_4$ and the last column contains the height
$\height(\alpha_i)$. For example, the positive root $\alpha_{19}$ is
\[
\alpha_{19} = 1 \cdot \alpha_1 + 2 \cdot \alpha_2 + 3 \cdot \alpha_3
+ 1 \cdot \alpha_4 = \frac12(e_1+e_2+e_3+e_4)
\]
and we have $\height(\alpha_{19}) = 7$.

The construction of the single root midafis of the group $UF_4$ is
similar to the one for $UE_i$: Let 
$R^{\norml}_4 := \Phi_4^+ \setminus \{\alpha_{20}, \alpha_{22}\}$ and
let $\alpha \in \Phi_4^+$. We choose an arm 
$a(\alpha) = \{\alpha_{i_1}, \alpha_{i_2}, \dots, \alpha_{i_r}\}$ 
of the hook $h(\alpha)$ such that the indices $i_1, i_2, \dots, i_r$ are 
given by the second column of Table~\ref{tab:armsf4} in the appendix. 
We define and construct the leg $\ell(\alpha)$, the source
$s(\alpha)$, the source group $S_\alpha$, the enlarged leg
$\bar{\ell}(\alpha)$ and the quotient pattern group
$\overline{T}_\alpha$ of $S_\alpha$ in the same way as for $UE_i$
(distinguishing the two cases $\alpha \not\in R^\norml_4$ and 
$\alpha \in R^\norml_4$).

\begin{prop} \label{prop:single_mida_f4}
Let $\Phi_4$ be a root system of type $F_4$ as described above. For
each positive root $\alpha \in \Phi_4^+$ the following are true:
\begin{enumerate}
\item[(a)] If $\alpha \in \Phi_4^+ \setminus R^{\norml}_4$ then 
  $\Psi_\alpha: \Irr^\lin(\overline{T}_\alpha) \times \Irr(X_\alpha)^*
  \to \Irr^\mida(UF_4)_\alpha$ with 
  \[
  (\mu, \lambda) \mapsto (\infl_{\overline{T}_\alpha}^{S_\alpha} \mu \cdot
   \infl_{X_\alpha}^{S_\alpha} \lambda)^{UF_4}
  \]
 is a one to one correspondence.

\item[(b)] If $\alpha \in R^{\norml}_4$ then 
  $\Psi_\alpha: \Irr(\overline{T}_\alpha) \times \Irr(X_\alpha)^*
  \to \Irr(UF_4)_\alpha$ with 
  \[
  (\mu, \lambda) \mapsto (\infl_{\overline{T}_\alpha}^{S_\alpha} \mu \cdot
   \infl_{X_\alpha}^{S_\alpha} \lambda)^{UF_4}
  \]
 is a one to one correspondence.
\end{enumerate} 
The number $|\Irr^\mida(UF_4)_\alpha|$ of midafis for $\alpha$ is
given in the second and the fifth column of Table~\ref{tab:mida_f4}
and the degree $\chi(1)$ for $\chi \in \Irr^\mida(UF_4)_\alpha$ is
given in the third and the sixth column of Table~\ref{tab:mida_f4}.   
\end{prop}

\bigskip

\begin{center}
 \begin{longtable}{|c|l|l||c|l|l|}
 \caption{Numbers and degrees of the midafis of ${UF}_4$.} \label{tab:mida_f4} \\

 \hline
 \hspace{-0.2cm} Root \hspace{-0.2cm} & Number of midafis & Degree & \hspace{-0.2cm} Root \hspace{-0.2cm} & Number of midafis & Degree\\ \hline
 \endfirsthead

 \multicolumn{6}{c}{{\tablename\ \thetable{} (cont.)}} \\
 \hline
 \hspace{-0.2cm} Root \hspace{-0.2cm} & Number of midafis & Degree & \hspace{-0.2cm} Root \hspace{-0.2cm} & Number of midafis & Degree\\ \hline
 \endhead

 \hline 
 \endfoot

 \endlastfoot

$\alpha_{1}$ &
 $q-1$ & 
$1$ &

$\alpha_{2}$ &
 $q-1$ & 
$1$\\

$\alpha_{3}$ &
 $q-1$ & 
$1$ &

$\alpha_{4}$ &
 $q-1$ & 
$1$\\

$\alpha_{5}$ &
 $q-1$ & 
$q$ &

$\alpha_{6}$ &
 $q-1$ & 
$q$\\

$\alpha_{7}$ &
 $q-1$ & 
$q$ &

$\alpha_{8}$ &
 $q(q-1)$ & 
$q^2$\\

$\alpha_{9}$ &
 $q(q-1)$ & 
$q$ &

$\alpha_{10}$ &
 $q(q-1)$ & 
$q^2$\\

$\alpha_{11}$ &
 $q^3 (q-1)$ & 
$q^2$ &

$\alpha_{12}$ &
 $q^2 (q-1)$ & 
$q^3$\\

$\alpha_{13}$ &
 $q (q-1)$ & 
$q^3$ &

$\alpha_{14}$ &
 $q^2 (q-1)$ & 
$q^3$\\

$\alpha_{15}$ &
 $q^5 (q-1)$ & 
$q^4$ &

$\alpha_{16}$ &
 $q^2 (q-1)$ & 
$q^2$\\

$\alpha_{17}$ &
 $q^4 (q-1)$ & 
$q^5$ &

$\alpha_{18}$ &
 $q^4 (q-1)$ & 
$q^3$\\

$\alpha_{19}$ &
 $q^3 (q-1)$ & 
$q^6$ &

$\alpha_{20}$ &
 $q^4 (q-1)$ & 
$q^4$\\

$\alpha_{21}$ &
 $q^3 (q-1)$ & 
$q^7$ &

$\alpha_{22}$ &
 $q^4 (q-1)$ & 
$q^5$\\

$\alpha_{23}$ &
 $q^4 (q-1)$ & 
$q^6$ &

$\alpha_{24}$ &
 $q^3 (q-1)$ & 
$q^7$\\
  \hline
 \end{longtable}
\end{center}

\pr 
The proof consists of computer calculations carried out by the
CHEVIE programs mentioned in Section~\ref{subsec:E6E7E8}. The
proof is analogous to the proof of Proposition~\ref{prop:single_mida_e678}  
with Table~\ref{tab:armse8} replaced by Table~\ref{tab:armsf4}.
\epr


\subsection{\texorpdfstring{Type $G_2$}{Type G2}} 
\label{subsec:G2}

We construct a root system of type~$G_2$ as in~\cite[Section~12.1]{Humphreys}:  
Let $e_1, e_2, e_3 \in \R^3$ be the usual orthonormal unit vectors
which form a basis of $\R^3$. Then 
$\Phi_2 := \pm \{ e_1-e_2, e_2-e_3, e_1-e_3, 2e_1-e_2-e_3,
2e_2-e_1-e_3, 2e_3-e_1-e_2 \}$ is a root system of type~$G_2$ 
and the set 
$\{\alpha_1, \alpha_2\}$, where $\alpha_1 := e_1-e_2$ and 
$\alpha_2 := -2e_1+e_2+e_3$, is a set of simple
roots. In particular, $\alpha$ is a short root and $\beta$ is 
a long root. The corresponding set of positive roots is  
\[
\Phi_2^+ = \{ \alpha_1, \alpha_2, \alpha_3 := \alpha_1+\alpha_2,
\alpha_4 := 2\alpha_1+\alpha_2, \alpha_5 := 3\alpha_1+\alpha_2,
\alpha_6 := 3\alpha_1+2\alpha_2\}.
\]
The construction of the single root midafis of the group $UG_2$ is
similar to the case $\alpha \in R^\norml_i$ for $UE_i$: We define arms
of the hooks for the positive roots as follows
\[
a(\alpha_1) := a(\alpha_2) := \emptyset, \quad 
a(\alpha_3) := a(\alpha_4) := a(\alpha_5) := \{\alpha_1\}, \quad
a(\alpha_6) := \{\alpha_2, \alpha_3\}
\]
with corresponding legs $\ell(\alpha_i)$. Let $\alpha \in \Phi_2^+$. 
It is easy to see that $s(\alpha) := \Phi_2^+ \setminus a(\alpha)$ is
a closed pattern and that 
$\ell(\alpha) \cup \{\alpha\} \cup k(\alpha) \unlhd \Phi_2^+$
so that we can consider the quotient pattern group 
$\overline{T}_\alpha := P(s(\alpha))/P(\{\alpha\} \cup \ell(\alpha) \cup k(\alpha))$.

\begin{prop} \label{prop:single_mida_g2}
Let $\Phi_2$ be a root system of type $G_2$ as described above. For
each positive root $\alpha \in \Phi_2^+$ the map
$\Psi_\alpha: \Irr(\overline{T}_\alpha) \times \Irr(X_\alpha)^* \to \Irr(UG_2)_\alpha$ 
with 
\[
(\mu, \lambda) \mapsto (\infl_{\overline{T}_\alpha}^{S_\alpha} \mu \cdot
   \infl_{X_\alpha}^{S_\alpha} \lambda)^{UE_i}
\]
is a one to one correspondence.
The number $|\Irr^\mida(UG_2)_\alpha|$ of midafis for the root $\alpha$ is
given in the second and the fifth column of Table~\ref{tab:mida_g2}
and the degree $\chi(1)$ for $\chi \in \Irr^\mida(UG_2)_\alpha$ is
given in the third and the sixth column of Table~\ref{tab:mida_g2}.   
\end{prop}

\begin{center}
 \begin{longtable}{|c|l|l||c|l|l|}
 \caption{Numbers and degrees of the midafis of ${UG}_2$.} \label{tab:mida_g2} \\

 \hline
 \hspace{-0.2cm} Root \hspace{-0.2cm} & Number of midafis & Degree & \hspace{-0.2cm} Root \hspace{-0.2cm} & Number of midafis & Degree\\ \hline
 \endfirsthead

 \multicolumn{6}{c}{{\tablename\ \thetable{} (cont.)}} \\
 \hline
 \hspace{-0.2cm} Root \hspace{-0.2cm} & Number of midafis & Degree & \hspace{-0.2cm} Root \hspace{-0.2cm} & Number of midafis & Degree\\ \hline
 \endhead

 \hline 
 \endfoot

 \endlastfoot

$\alpha_{1}$ & 
$q-1$ & 
$1$&

$\alpha_{2}$ & 
$q-1$ & 
$1$\\

$\alpha_{3}$ & 
$q-1$ & 
$q$&

$\alpha_{4}$ & 
$q(q-1)$ & 
$q$\\

$\alpha_{5}$ & 
$q^2 (q-1)$ & 
$q$&

$\alpha_{6}$ & 
$q(q-1)$ & 
$q^2$\\
  \hline
 \end{longtable}
\end{center}

\pr
The proof is analogous to the proof of Propositions
\ref{prop:single_mida_e678} and \ref{prop:single_mida_f4} (but the
calculations can be carried out by hand).
\epr

\begin{rem} \label{rem:mida_g2}
It turns out that the group $\overline{T}_\alpha$ is elementary
abelian for all $\alpha \in \Phi_2^+$ so that we have
$\Irr(UG_2)_\alpha = \Irr^\mida(UG_2)_\alpha$ for all $\alpha \in \Phi_2^+$.
It follows that the only characters $\chi \in \Irr(UG_2)$ which 
are not single root characters and are not covered by
Proposition~\ref{prop:single_mida_g2} are the linear characters 
$\chi \in \Irr(UG_2)$ with $|\rs(\chi)|=2$. 
\end{rem}

\medskip

\begin{rem} \label{rem:choicearm}
Let $i \in \{2,4,6,7,8\}$ and $\Phi_i$ a root system of type $G_2$,
$F_4$, $E_6$, $E_7$, $E_8$, respectively, with set of positive roots
$\Phi_i^+$ as in~\ref{subsec:E6E7E8}-\ref{subsec:G2}. We use the
notation from~\ref{subsec:E6E7E8}-\ref{subsec:G2}. The arms
$a(\alpha)$ in Tables~\ref{tab:armse8}, \ref{tab:armsf4} and
Section~\ref{subsec:G2} are chosen such that for all $\alpha \in \Phi_i^+$
the following condition is satisfied:
\begin{enumerate}
\item The source $s(\alpha) = \Phi_i^+ \setminus a(\alpha)$ is a 
  closed pattern. 
\end{enumerate}
For all roots $\alpha \in \Phi_i^+ \cap R^\norml_i$ the choice of
$a(\alpha)$ in Tables~\ref{tab:armse8}, \ref{tab:armsf4} and
Section~\ref{subsec:G2} implies that the corresponding leg $\ell(\alpha)$
satisfies the following condition
\begin{enumerate}
\item[(2)] $\ell(\alpha) \cup k(\alpha) \unlhd s(\alpha)$.
\end{enumerate}
Now suppose that $\alpha \in \Phi_i^+ \setminus R^\norml_i$ (where 
$R^\norml_2 := \Phi_2^+$). In this case there is no choice of
$a(\alpha)$ such that (1) and~(2) are satisfied simultaneously. For
each choice of $a(\alpha)$ with corresponding leg $\ell(\alpha)$ 
such that condition (1) is satisfied let $M$ be the normal closure 
of $\ell(\alpha) \cup k(\alpha)$ in $s(\alpha)$ and define 
$\bar{\ell}(\alpha) := M \setminus k(\alpha)$
so that $\bar{\ell}(\alpha) \supsetneq \ell(\alpha)$. For all 
$\alpha \in \Phi_i^+ \setminus R^\norml_i$ the choice of $a(\alpha)$ in
Tables~\ref{tab:armse8} and \ref{tab:armsf4} implies that
$\bar{\ell}(\alpha)$ has the following properties:
\begin{enumerate}
\item[(3)] $\bar{\ell}(\alpha) \cup k(\alpha) \unlhd s(\alpha)$, 
$\{\alpha\} \cup \bar{\ell}(\alpha) \cup k(\alpha) \unlhd \Phi^+$
and the quotient pattern group $P(\bar{\ell}(\alpha) \cup k(\alpha))/P(k(\alpha))$ 
is abelian. 
\end{enumerate}
Among all choices of the arm $a(\alpha)$ such that conditions (1) and
(3) hold the choice in Tables~\ref{tab:armse8} and \ref{tab:armsf4}
minimizes $|\bar{\ell}(\alpha)|$. This is achieved as follows:
Let
\[
\Pi_\alpha := \{(\alpha_j, \alpha_{j'}) \in \Phi_i^+ \times \Phi_i^+
\mid \alpha_j+\alpha_{j'} = \alpha \ \text{and} \ j<j'\}.
\]
Among all pairs $(\alpha_j, \alpha_{j'}) \in \Pi_\alpha$ we choose the
(unique) pair where the first index $j$ is maximal (since the roots of each
root system are labeled by increasing height this guarantees that both 
$\height(\alpha_i)$ and $\height(\alpha_j)$ are ``not too small''). 
Write $(\gamma, \gamma')$ for this pair of roots. 

For each choice of $a(\alpha)$ we have 
${\rm heart}(\alpha) := \Phi_i^+ \setminus h(\alpha) \subseteq \Phi_i^+ \setminus a(\alpha) = s(\alpha)$. 
Let $N_\gamma$ be the normal closure of $\{\gamma\}$ in the closed
pattern generated by ${\rm heart}(\alpha) \cup \{\gamma\}$ and $N_{\gamma'}$ 
the normal closure of $\{\gamma'\}$ in the closed pattern generated by
${\rm heart}(\alpha) \cup \{\gamma'\}$. Suppose that the arm $a(\alpha)$ is
chosen such that the conditions (1) and~(2) or that the conditions (1) 
and (3) hold. Then we have either $\gamma \in \ell(\alpha)$ or 
$\gamma' \in \ell(\alpha)$. In the first case we have 
$N_\gamma \cap h(\alpha) \subseteq \ell(\alpha)$ and in the second case
$N_{\gamma'} \cap h(\alpha) \subseteq \ell(\alpha)$. Hence in both
cases we have $N_\gamma \cap N_{\gamma'} \cap h(\alpha) \subseteq \ell(\alpha)$.  
This reduces the number of possible choices for $\ell(\alpha)$ and
hence $a(\alpha)$ such that the conditions (1) and (2) or the
conditions (1) and (3) are satisfied considerably. 

Then we do an exhaustive search over the remaining
possibilities to filter out those choices which satisfy conditions (1)
and (2). If there is such a choice then we know that 
$\alpha \in R^\norml_i$ and we choose $a(\alpha)$
such that conditions (1) and~(2) hold. If there is no such choice then
we know that $\alpha \in \Phi_i^+ \setminus R^\norml_i$ and we run through all
possibilities satisfying conditions (1) and (3) and choose $a(\alpha)$
so that $|\bar{\ell}(\alpha)|$ is minimal.

We consider the example $i=8$ (that is, $\Phi_i$ is of type $E_8$)
and $\alpha = \alpha_{115}$. The hook $h(\alpha)$ was already
determined in the proof of Proposition~\ref{prop:single_mida_e678}. 
We have $|h(\alpha)| = 47$. Hence there are 
$2^{(|h(\alpha)|-1)/2} = 2^{23}= 8388608$ possible choices for 
the arm~$a(\alpha)$. The pair $(\alpha_j,\alpha_{j'}) \in \Pi_\alpha$
with maximal first index $j$ is $(\gamma,\gamma') = (\alpha_{75}, \alpha_{78})$.
Using the CHEVIE programs we get
\begin{eqnarray*}
N_\gamma \cap N_{\gamma'} \cap h(\alpha) & = & \{ \alpha_{86},
\alpha_{90}, \alpha_{91}, \alpha_{95}, \alpha_{98}, \alpha_{100},
\alpha_{101}, \alpha_{102}, \alpha_{104}, \alpha_{105},
\alpha_{107},\\
&& \alpha_{108}, \alpha_{109}, \alpha_{112}, \alpha_{113}, \alpha_{114}\}.
\end{eqnarray*}
It follows that there are at most $2^{23-|N_\gamma \cap N_{\gamma'}
 \cap h(\alpha)|} = 2^{23-16} = 2^7 =128$ possible choices such that
conditions (1) and (2) or (1) and (3) are satisfied. Testing these
$128$ possibilities we see that there is no choice of $a(\alpha)$ such
that the conditions (1) and (2) are satisfied simultaneously. Thus
$\alpha \in \Phi_8^+ \setminus R^\norml_{6/7/8}$. Furthermore we see that for all choices
of $a(\alpha)$ such that the conditions (1) and (3) are satisfied
simultaneously we have $|\bar{\ell}(\alpha) \setminus \ell(\alpha)|
\ge 7$. Hence the choice of $a(\alpha)$ in Table~\ref{tab:armse8}
minimizes $|\bar{\ell}(\alpha)|$.
\end{rem}

\subsection{Proof of Theorems \ref{thm:main_exceptional} and \ref{thm:main_exceptionalnormal}}
\label{subsec:proofmainthmexc}

We can now complete the proof of the main results stated in the
introduction:

\medskip

\pr (of Theorem~\ref{thm:main_exceptional})
Let $\Phi_i$ be a root system of type $E_6$, $E_7$, $E_8$, $F_4$ or
$G_2$ as in Sections~\ref{subsec:E6E7E8}-\ref{subsec:G2} and 
$\alpha \in \Phi_i^+$. If $i \neq 2$ and $\alpha \in \Phi_i^+ \setminus R^\norml_i$ 
then the statement of the theorem follows from
Proposition~\ref{prop:single_mida_e678}~(a) and 
Proposition \ref{prop:single_mida_f4} (a). 

Suppose that $i=2$ or that $\alpha \in R^\norml_i$. Considering
degrees we see that the one to one correspondences $\Psi_\alpha$ 
in Propositions~\ref{prop:single_mida_e678}~(b) and
\ref{prop:single_mida_f4}~(b) and in Proposition~\ref{prop:single_mida_g2} 
map $\Irr^\lin(\overline{T}_\alpha) \times \Irr(X_\alpha)^*$ onto
$\Irr^\mida(UY_i)_\alpha$. This completes the proof of
Theorem~\ref{thm:main_exceptional}.
\epr

\bigskip

\pr (of Theorem~\ref{thm:main_exceptionalnormal})
The theorem follows from Propositions~\ref{prop:single_mida_e678}~(b),  
\ref{prop:single_mida_f4}~(b) and Proposition~\ref{prop:single_mida_g2}.
\epr


\section*{Appendix}
\label{sec:appendix}

\setcounter{table}{0}
\renewcommand{\thetable}{\textrm{A.\arabic{table}}}

\begin{center}
 \begin{longtable}{|c|l|l|c|}
 \caption{Positive roots in the root system $\Phi_8$ of type $E_8$.} \label{tab:rootse8} \\

 \hline
 Root & Linear combination & Linear combination of $e_1$, \dots, $e_8$ & Height\\
      & of simple roots    && \\
      & $\alpha_1 \alpha_2 \alpha_3 \alpha_4 \alpha_5 \alpha_6 \alpha_7 \alpha_8$ && \\ \hline
 \endfirsthead

 \multicolumn{4}{c}{{\tablename\ \thetable{} (cont.)}} \\
 \hline
 Root & $\alpha_1 \alpha_2 \alpha_3 \alpha_4 \alpha_5 \alpha_6
 \alpha_7 \alpha_8$ & Linear combination of $e_1$, \dots, $e_8$ & Height\\ \hline
 \endhead

 \hline 
 \endfoot

 \endlastfoot

$\alpha_{1}$ & $1 \ \  0 \ \ 0 \ \ 0 \ \ 0 \ \ 0 \ \ 0 \ \ 0$ & $\frac12(e_1-e_2-e_3-e_4-e_5-e_6-e_7+e_8)$ & $1$\\
$\alpha_{2}$ & $0 \ \ 1 \ \ 0 \ \ 0 \ \ 0 \ \ 0 \ \ 0 \ \ 0$ & $e_1+e_2$ & $1$\\
$\alpha_{3}$ & $0 \ \ 0 \ \ 1 \ \ 0 \ \ 0 \ \ 0 \ \ 0 \ \ 0$ & $-e_1+e_2$ & $1$\\
$\alpha_{4}$ & $0 \ \ 0 \ \ 0 \ \ 1 \ \ 0 \ \ 0 \ \ 0 \ \ 0$ & $-e_2+e_3$ & $1$\\
$\alpha_{5}$ & $0 \ \ 0 \ \ 0 \ \ 0 \ \ 1 \ \ 0 \ \ 0 \ \ 0$ & $-e_3+e_4$ & $1$\\
$\alpha_{6}$ & $0 \ \ 0 \ \ 0 \ \ 0 \ \ 0 \ \ 1 \ \ 0 \ \ 0$ & $-e_4+e_5$ & $1$\\
$\alpha_{7}$ & $0 \ \ 0 \ \ 0 \ \ 0 \ \ 0 \ \ 0 \ \ 1 \ \ 0$ & $-e_5+e_6$ & $1$\\
$\alpha_{8}$ & $0 \ \ 0 \ \ 0 \ \ 0 \ \ 0 \ \ 0 \ \ 0 \ \ 1$ & $-e_6+e_7$ & $1$\\
$\alpha_{9}$ & $1 \ \ 0 \ \ 1 \ \ 0 \ \ 0 \ \ 0 \ \ 0 \ \ 0$ &   $\frac12(-e_1+e_2-e_3-e_4-e_5-e_6-e_7+e_8)$  & $2$\\
$\alpha_{10}$ & $0 \ \ 1 \ \ 0 \ \ 1 \ \ 0 \ \ 0 \ \ 0 \ \ 0$ & $e_1+e_3$ & $2$\\
$\alpha_{11}$ & $0 \ \ 0 \ \ 1 \ \ 1 \ \ 0 \ \ 0 \ \ 0 \ \ 0$ & $-e_1+e_3$ & $2$\\
$\alpha_{12}$ & $0 \ \ 0 \ \ 0 \ \ 1 \ \ 1 \ \ 0 \ \ 0 \ \ 0$ & $-e_2+e_4$ & $2$\\
$\alpha_{13}$ & $0 \ \ 0 \ \ 0 \ \ 0 \ \ 1 \ \ 1 \ \ 0 \ \ 0$ & $-e_3+e_5$ & $2$\\
$\alpha_{14}$ & $0 \ \ 0 \ \ 0 \ \ 0 \ \ 0 \ \ 1 \ \ 1 \ \ 0$ & $-e_4+e_6$ & $2$\\
$\alpha_{15}$ & $0 \ \ 0 \ \ 0 \ \ 0 \ \ 0 \ \ 0 \ \ 1 \ \ 1$ & $-e_5+e_7$ & $2$\\
$\alpha_{16}$ & $1 \ \ 0 \ \ 1 \ \ 1 \ \ 0 \ \ 0 \ \ 0 \ \ 0$ &   $\frac12(-e_1-e_2+e_3-e_4-e_5-e_6-e_7+e_8)$  & $3$\\
$\alpha_{17}$ & $0 \ \ 1 \ \ 1 \ \ 1 \ \ 0 \ \ 0 \ \ 0 \ \ 0$ & $e_2+e_3$ & $3$\\
$\alpha_{18}$ & $0 \ \ 1 \ \ 0 \ \ 1 \ \ 1 \ \ 0 \ \ 0 \ \ 0$ & $e_1+e_4$ & $3$\\
$\alpha_{19}$ & $0 \ \ 0 \ \ 1 \ \ 1 \ \ 1 \ \ 0 \ \ 0 \ \ 0$ & $-e_1+e_4$ & $3$\\
$\alpha_{20}$ & $0 \ \ 0 \ \ 0 \ \ 1 \ \ 1 \ \ 1 \ \ 0 \ \ 0$ & $-e_2+e_5$ & $3$\\
$\alpha_{21}$ & $0 \ \ 0 \ \ 0 \ \ 0 \ \ 1 \ \ 1 \ \ 1 \ \ 0$ & $-e_3+e_6$ & $3$\\
$\alpha_{22}$ & $0 \ \ 0 \ \ 0 \ \ 0 \ \ 0 \ \ 1 \ \ 1 \ \ 1$ & $-e_4+e_7$ & $3$\\
$\alpha_{23}$ & $1 \ \ 1 \ \ 1 \ \ 1 \ \ 0 \ \ 0 \ \ 0 \ \ 0$ &   $\frac12(e_1+e_2+e_3-e_4-e_5-e_6-e_7+e_8)$  & $4$\\
$\alpha_{24}$ & $1 \ \ 0 \ \ 1 \ \ 1 \ \ 1 \ \ 0 \ \ 0 \ \ 0$ &   $\frac12(-e_1-e_2-e_3+e_4-e_5-e_6-e_7+e_8)$  & $4$\\
$\alpha_{25}$ & $0 \ \ 1 \ \ 1 \ \ 1 \ \ 1 \ \ 0 \ \ 0 \ \ 0$ & $e_2+e_4$ & $4$\\
$\alpha_{26}$ & $0 \ \ 1 \ \ 0 \ \ 1 \ \ 1 \ \ 1 \ \ 0 \ \ 0$ & $e_1+e_5$ & $4$\\
$\alpha_{27}$ & $0 \ \ 0 \ \ 1 \ \ 1 \ \ 1 \ \ 1 \ \ 0 \ \ 0$ & $-e_1+e_5$ & $4$\\
$\alpha_{28}$ & $0 \ \ 0 \ \ 0 \ \ 1 \ \ 1 \ \ 1 \ \ 1 \ \ 0$ & $-e_2+e_6$ & $4$\\
$\alpha_{29}$ & $0 \ \ 0 \ \ 0 \ \ 0 \ \ 1 \ \ 1 \ \ 1 \ \ 1$ & $-e_3+e_7$ & $4$\\
$\alpha_{30}$ & $1 \ \ 1 \ \ 1 \ \ 1 \ \ 1 \ \ 0 \ \ 0 \ \ 0$ &   $\frac12(e_1+e_2-e_3+e_4-e_5-e_6-e_7+e_8)$  & $5$\\
$\alpha_{31}$ & $1 \ \ 0 \ \ 1 \ \ 1 \ \ 1 \ \ 1 \ \ 0 \ \ 0$ &   $\frac12(-e_1-e_2-e_3-e_4+e_5-e_6-e_7+e_8)$  & $5$\\
$\alpha_{32}$ & $0 \ \ 1 \ \ 1 \ \ 2 \ \ 1 \ \ 0 \ \ 0 \ \ 0$ & $e_3+e_4$ & $5$\\
$\alpha_{33}$ & $0 \ \ 1 \ \ 1 \ \ 1 \ \ 1 \ \ 1 \ \ 0 \ \ 0$ & $e_2+e_5$ & $5$\\
$\alpha_{34}$ & $0 \ \ 1 \ \ 0 \ \ 1 \ \ 1 \ \ 1 \ \ 1 \ \ 0$ & $e_1+e_6$ & $5$\\
$\alpha_{35}$ & $0 \ \ 0 \ \ 1 \ \ 1 \ \ 1 \ \ 1 \ \ 1 \ \ 0$ & $-e_1+e_6$ & $5$\\
$\alpha_{36}$ & $0 \ \ 0 \ \ 0 \ \ 1 \ \ 1 \ \ 1 \ \ 1 \ \ 1$ & $-e_2+e_7$ & $5$\\
$\alpha_{37}$ & $1 \ \ 1 \ \ 1 \ \ 2 \ \ 1 \ \ 0 \ \ 0 \ \ 0$ &   $\frac12(e_1-e_2+e_3+e_4-e_5-e_6-e_7+e_8)$  & $6$\\
$\alpha_{38}$ & $1 \ \ 1 \ \ 1 \ \ 1 \ \ 1 \ \ 1 \ \ 0 \ \ 0$ &   $\frac12(e_1+e_2-e_3-e_4+e_5-e_6-e_7+e_8)$  & $6$\\
$\alpha_{39}$ & $1 \ \ 0 \ \ 1 \ \ 1 \ \ 1 \ \ 1 \ \ 1 \ \ 0$ &   $\frac12(-e_1-e_2-e_3-e_4-e_5+e_6-e_7+e_8)$  & $6$\\
$\alpha_{40}$ & $0 \ \ 1 \ \ 1 \ \ 2 \ \ 1 \ \ 1 \ \ 0 \ \ 0$ & $e_3+e_5$ & $6$\\
$\alpha_{41}$ & $0 \ \ 1 \ \ 1 \ \ 1 \ \ 1 \ \ 1 \ \ 1 \ \ 0$ & $e_2+e_6$ & $6$\\
$\alpha_{42}$ & $0 \ \ 1 \ \ 0 \ \ 1 \ \ 1 \ \ 1 \ \ 1 \ \ 1$ & $e_1+e_7$ & $6$\\
$\alpha_{43}$ & $0 \ \ 0 \ \ 1 \ \ 1 \ \ 1 \ \ 1 \ \ 1 \ \ 1$ & $-e_1+e_7$ & $6$\\
$\alpha_{44}$ & $1 \ \ 1 \ \ 2 \ \ 2 \ \ 1 \ \ 0 \ \ 0 \ \ 0$ &   $\frac12(-e_1+e_2+e_3+e_4-e_5-e_6-e_7+e_8)$  & $7$\\
$\alpha_{45}$ & $1 \ \ 1 \ \ 1 \ \ 2 \ \ 1 \ \ 1 \ \ 0 \ \ 0$ &   $\frac12(e_1-e_2+e_3-e_4+e_5-e_6-e_7+e_8)$  & $7$\\
$\alpha_{46}$ & $1 \ \ 1 \ \ 1 \ \ 1 \ \ 1 \ \ 1 \ \ 1 \ \ 0$ &   $\frac12(e_1+e_2-e_3-e_4-e_5+e_6-e_7+e_8)$  & $7$\\
$\alpha_{47}$ & $1 \ \ 0 \ \ 1 \ \ 1 \ \ 1 \ \ 1 \ \ 1 \ \ 1$ &   $\frac12(-e_1-e_2-e_3-e_4-e_5-e_6+e_7+e_8)$  & $7$\\
$\alpha_{48}$ & $0 \ \ 1 \ \ 1 \ \ 2 \ \ 2 \ \ 1 \ \ 0 \ \ 0$ & $e_4+e_5$ & $7$\\
$\alpha_{49}$ & $0 \ \ 1 \ \ 1 \ \ 2 \ \ 1 \ \ 1 \ \ 1 \ \ 0$ & $e_3+e_6$ & $7$\\
$\alpha_{50}$ & $0 \ \ 1 \ \ 1 \ \ 1 \ \ 1 \ \ 1 \ \ 1 \ \ 1$ & $e_2+e_7$ & $7$\\
$\alpha_{51}$ & $1 \ \ 1 \ \ 2 \ \ 2 \ \ 1 \ \ 1 \ \ 0 \ \ 0$ &   $\frac12(-e_1+e_2+e_3-e_4+e_5-e_6-e_7+e_8)$  & $8$\\
$\alpha_{52}$ & $1 \ \ 1 \ \ 1 \ \ 2 \ \ 2 \ \ 1 \ \ 0 \ \ 0$ &   $\frac12(e_1-e_2-e_3+e_4+e_5-e_6-e_7+e_8)$  & $8$\\
$\alpha_{53}$ & $1 \ \ 1 \ \ 1 \ \ 2 \ \ 1 \ \ 1 \ \ 1 \ \ 0$ &   $\frac12(e_1-e_2+e_3-e_4-e_5+e_6-e_7+e_8)$  & $8$\\
$\alpha_{54}$ & $1 \ \ 1 \ \ 1 \ \ 1 \ \ 1 \ \ 1 \ \ 1 \ \ 1$ &   $\frac12(e_1+e_2-e_3-e_4-e_5-e_6+e_7+e_8)$  & $8$\\
$\alpha_{55}$ & $0 \ \ 1 \ \ 1 \ \ 2 \ \ 2 \ \ 1 \ \ 1 \ \ 0$ & $e_4+e_6$ & $8$\\
$\alpha_{56}$ & $0 \ \ 1 \ \ 1 \ \ 2 \ \ 1 \ \ 1 \ \ 1 \ \ 1$ & $e_3+e_7$ & $8$\\
$\alpha_{57}$ & $1 \ \ 1 \ \ 2 \ \ 2 \ \ 2 \ \ 1 \ \ 0 \ \ 0$ &   $\frac12(-e_1+e_2-e_3+e_4+e_5-e_6-e_7+e_8)$  & $9$\\
$\alpha_{58}$ & $1 \ \ 1 \ \ 2 \ \ 2 \ \ 1 \ \ 1 \ \ 1 \ \ 0$ &   $\frac12(-e_1+e_2+e_3-e_4-e_5+e_6-e_7+e_8)$  & $9$\\
$\alpha_{59}$ & $1 \ \ 1 \ \ 1 \ \ 2 \ \ 2 \ \ 1 \ \ 1 \ \ 0$ &   $\frac12(e_1-e_2-e_3+e_4-e_5+e_6-e_7+e_8)$  & $9$\\
$\alpha_{60}$ & $1 \ \ 1 \ \ 1 \ \ 2 \ \ 1 \ \ 1 \ \ 1 \ \ 1$ &   $\frac12(e_1-e_2+e_3-e_4-e_5-e_6+e_7+e_8)$  & $9$\\
$\alpha_{61}$ & $0 \ \ 1 \ \ 1 \ \ 2 \ \ 2 \ \ 2 \ \ 1 \ \ 0$ & $e_5+e_6$ & $9$\\
$\alpha_{62}$ & $0 \ \ 1 \ \ 1 \ \ 2 \ \ 2 \ \ 1 \ \ 1 \ \ 1$ & $e_4+e_7$ & $9$\\
$\alpha_{63}$ & $1 \ \ 1 \ \ 2 \ \ 3 \ \ 2 \ \ 1 \ \ 0 \ \ 0$ &   $\frac12(-e_1-e_2+e_3+e_4+e_5-e_6-e_7+e_8)$  & $10$\\
$\alpha_{64}$ & $1 \ \ 1 \ \ 2 \ \ 2 \ \ 2 \ \ 1 \ \ 1 \ \ 0$ &   $\frac12(-e_1+e_2-e_3+e_4-e_5+e_6-e_7+e_8)$  & $10$\\
$\alpha_{65}$ & $1 \ \ 1 \ \ 2 \ \ 2 \ \ 1 \ \ 1 \ \ 1 \ \ 1$ &   $\frac12(-e_1+e_2+e_3-e_4-e_5-e_6+e_7+e_8)$  & $10$\\
$\alpha_{66}$ & $1 \ \ 1 \ \ 1 \ \ 2 \ \ 2 \ \ 2 \ \ 1 \ \ 0$ &   $\frac12(e_1-e_2-e_3-e_4+e_5+e_6-e_7+e_8)$  & $10$\\
$\alpha_{67}$ & $1 \ \ 1 \ \ 1 \ \ 2 \ \ 2 \ \ 1 \ \ 1 \ \ 1$ &   $\frac12(e_1-e_2-e_3+e_4-e_5-e_6+e_7+e_8)$  & $10$\\
$\alpha_{68}$ & $0 \ \ 1 \ \ 1 \ \ 2 \ \ 2 \ \ 2 \ \ 1 \ \ 1$ & $e_5+e_7$ & $10$\\
$\alpha_{69}$ & $1 \ \ 2 \ \ 2 \ \ 3 \ \ 2 \ \ 1 \ \ 0 \ \ 0$ &   $\frac12(e_1+e_2+e_3+e_4+e_5-e_6-e_7+e_8)$  & $11$\\
$\alpha_{70}$ & $1 \ \ 1 \ \ 2 \ \ 3 \ \ 2 \ \ 1 \ \ 1 \ \ 0$ &   $\frac12(-e_1-e_2+e_3+e_4-e_5+e_6-e_7+e_8)$  & $11$\\
$\alpha_{71}$ & $1 \ \ 1 \ \ 2 \ \ 2 \ \ 2 \ \ 2 \ \ 1 \ \ 0$ &   $\frac12(-e_1+e_2-e_3-e_4+e_5+e_6-e_7+e_8)$  & $11$\\
$\alpha_{72}$ & $1 \ \ 1 \ \ 2 \ \ 2 \ \ 2 \ \ 1 \ \ 1 \ \ 1$ &   $\frac12(-e_1+e_2-e_3+e_4-e_5-e_6+e_7+e_8)$  & $11$\\
$\alpha_{73}$ & $1 \ \ 1 \ \ 1 \ \ 2 \ \ 2 \ \ 2 \ \ 1 \ \ 1$ &   $\frac12(e_1-e_2-e_3-e_4+e_5-e_6+e_7+e_8)$  & $11$\\
$\alpha_{74}$ & $0 \ \ 1 \ \ 1 \ \ 2 \ \ 2 \ \ 2 \ \ 2 \ \ 1$ & $e_6+e_7$ & $11$\\
$\alpha_{75}$ & $1 \ \ 2 \ \ 2 \ \ 3 \ \ 2 \ \ 1 \ \ 1 \ \ 0$ &   $\frac12(e_1+e_2+e_3+e_4-e_5+e_6-e_7+e_8)$  & $12$\\
$\alpha_{76}$ & $1 \ \ 1 \ \ 2 \ \ 3 \ \ 2 \ \ 2 \ \ 1 \ \ 0$ &   $\frac12(-e_1-e_2+e_3-e_4+e_5+e_6-e_7+e_8)$  & $12$\\
$\alpha_{77}$ & $1 \ \ 1 \ \ 2 \ \ 3 \ \ 2 \ \ 1 \ \ 1 \ \ 1$ &   $\frac12(-e_1-e_2+e_3+e_4-e_5-e_6+e_7+e_8)$  & $12$\\
$\alpha_{78}$ & $1 \ \ 1 \ \ 2 \ \ 2 \ \ 2 \ \ 2 \ \ 1 \ \ 1$ &   $\frac12(-e_1+e_2-e_3-e_4+e_5-e_6+e_7+e_8)$  & $12$\\
$\alpha_{79}$ & $1 \ \ 1 \ \ 1 \ \ 2 \ \ 2 \ \ 2 \ \ 2 \ \ 1$ &   $\frac12(e_1-e_2-e_3-e_4-e_5+e_6+e_7+e_8)$  & $12$\\
$\alpha_{80}$ & $1 \ \ 2 \ \ 2 \ \ 3 \ \ 2 \ \ 2 \ \ 1 \ \ 0$ &   $\frac12(e_1+e_2+e_3-e_4+e_5+e_6-e_7+e_8)$  & $13$\\
$\alpha_{81}$ & $1 \ \ 2 \ \ 2 \ \ 3 \ \ 2 \ \ 1 \ \ 1 \ \ 1$ &   $\frac12(e_1+e_2+e_3+e_4-e_5-e_6+e_7+e_8)$  & $13$\\
$\alpha_{82}$ & $1 \ \ 1 \ \ 2 \ \ 3 \ \ 3 \ \ 2 \ \ 1 \ \ 0$ &   $\frac12(-e_1-e_2-e_3+e_4+e_5+e_6-e_7+e_8)$  & $13$\\
$\alpha_{83}$ & $1 \ \ 1 \ \ 2 \ \ 3 \ \ 2 \ \ 2 \ \ 1 \ \ 1$ &   $\frac12(-e_1-e_2+e_3-e_4+e_5-e_6+e_7+e_8)$  & $13$\\
$\alpha_{84}$ & $1 \ \ 1 \ \ 2 \ \ 2 \ \ 2 \ \ 2 \ \ 2 \ \ 1$ &   $\frac12(-e_1+e_2-e_3-e_4-e_5+e_6+e_7+e_8)$  & $13$\\
$\alpha_{85}$ & $1 \ \ 2 \ \ 2 \ \ 3 \ \ 3 \ \ 2 \ \ 1 \ \ 0$ &   $\frac12(e_1+e_2-e_3+e_4+e_5+e_6-e_7+e_8)$  & $14$\\
$\alpha_{86}$ & $1 \ \ 2 \ \ 2 \ \ 3 \ \ 2 \ \ 2 \ \ 1 \ \ 1$ &   $\frac12(e_1+e_2+e_3-e_4+e_5-e_6+e_7+e_8)$  & $14$\\
$\alpha_{87}$ & $1 \ \ 1 \ \ 2 \ \ 3 \ \ 3 \ \ 2 \ \ 1 \ \ 1$ &   $\frac12(-e_1-e_2-e_3+e_4+e_5-e_6+e_7+e_8)$  & $14$\\
$\alpha_{88}$ & $1 \ \ 1 \ \ 2 \ \ 3 \ \ 2 \ \ 2 \ \ 2 \ \ 1$ &   $\frac12(-e_1-e_2+e_3-e_4-e_5+e_6+e_7+e_8)$  & $14$\\
$\alpha_{89}$ & $1 \ \ 2 \ \ 2 \ \ 4 \ \ 3 \ \ 2 \ \ 1 \ \ 0$ &   $\frac12(e_1-e_2+e_3+e_4+e_5+e_6-e_7+e_8)$  & $15$\\
$\alpha_{90}$ & $1 \ \ 2 \ \ 2 \ \ 3 \ \ 3 \ \ 2 \ \ 1 \ \ 1$ &   $\frac12(e_1+e_2-e_3+e_4+e_5-e_6+e_7+e_8)$  & $15$\\
$\alpha_{91}$ & $1 \ \ 2 \ \ 2 \ \ 3 \ \ 2 \ \ 2 \ \ 2 \ \ 1$ &   $\frac12(e_1+e_2+e_3-e_4-e_5+e_6+e_7+e_8)$  & $15$\\
$\alpha_{92}$ & $1 \ \ 1 \ \ 2 \ \ 3 \ \ 3 \ \ 2 \ \ 2 \ \ 1$ &   $\frac12(-e_1-e_2-e_3+e_4-e_5+e_6+e_7+e_8)$  & $15$\\
$\alpha_{93}$ & $1 \ \ 2 \ \ 3 \ \ 4 \ \ 3 \ \ 2 \ \ 1 \ \ 0$ &   $\frac12(-e_1+e_2+e_3+e_4+e_5+e_6-e_7+e_8)$  & $16$\\
$\alpha_{94}$ & $1 \ \ 2 \ \ 2 \ \ 4 \ \ 3 \ \ 2 \ \ 1 \ \ 1$ &   $\frac12(e_1-e_2+e_3+e_4+e_5-e_6+e_7+e_8)$  & $16$\\
$\alpha_{95}$ & $1 \ \ 2 \ \ 2 \ \ 3 \ \ 3 \ \ 2 \ \ 2 \ \ 1$ &   $\frac12(e_1+e_2-e_3+e_4-e_5+e_6+e_7+e_8)$  & $16$\\
$\alpha_{96}$ & $1 \ \ 1 \ \ 2 \ \ 3 \ \ 3 \ \ 3 \ \ 2 \ \ 1$ &   $\frac12(-e_1-e_2-e_3-e_4+e_5+e_6+e_7+e_8)$  & $16$\\
$\alpha_{97}$ & $2 \ \ 2 \ \ 3 \ \ 4 \ \ 3 \ \ 2 \ \ 1 \ \ 0$ & $-e_7+e_8$ & $17$\\
$\alpha_{98}$ & $1 \ \ 2 \ \ 3 \ \ 4 \ \ 3 \ \ 2 \ \ 1 \ \ 1$ &   $\frac12(-e_1+e_2+e_3+e_4+e_5-e_6+e_7+e_8)$  & $17$\\
$\alpha_{99}$ & $1 \ \ 2 \ \ 2 \ \ 4 \ \ 3 \ \ 2 \ \ 2 \ \ 1$ &   $\frac12(e_1-e_2+e_3+e_4-e_5+e_6+e_7+e_8)$  & $17$\\
$\alpha_{100}$ & $1 \ \ 2 \ \ 2 \ \ 3 \ \ 3 \ \ 3 \ \ 2 \ \ 1$ &   $\frac12(e_1+e_2-e_3-e_4+e_5+e_6+e_7+e_8)$  & $17$\\
$\alpha_{101}$ & $2 \ \ 2 \ \ 3 \ \ 4 \ \ 3 \ \ 2 \ \ 1 \ \ 1$ & $-e_6+e_8$ & $18$\\
$\alpha_{102}$ & $1 \ \ 2 \ \ 3 \ \ 4 \ \ 3 \ \ 2 \ \ 2 \ \ 1$ &   $\frac12(-e_1+e_2+e_3+e_4-e_5+e_6+e_7+e_8)$  & $18$\\
$\alpha_{103}$ & $1 \ \ 2 \ \ 2 \ \ 4 \ \ 3 \ \ 3 \ \ 2 \ \ 1$ &   $\frac12(e_1-e_2+e_3-e_4+e_5+e_6+e_7+e_8)$  & $18$\\
$\alpha_{104}$ & $2 \ \ 2 \ \ 3 \ \ 4 \ \ 3 \ \ 2 \ \ 2 \ \ 1$ & $-e_5+e_8$ & $19$\\
$\alpha_{105}$ & $1 \ \ 2 \ \ 3 \ \ 4 \ \ 3 \ \ 3 \ \ 2 \ \ 1$ &   $\frac12(-e_1+e_2+e_3-e_4+e_5+e_6+e_7+e_8)$  & $19$\\
$\alpha_{106}$ & $1 \ \ 2 \ \ 2 \ \ 4 \ \ 4 \ \ 3 \ \ 2 \ \ 1$ &   $\frac12(e_1-e_2-e_3+e_4+e_5+e_6+e_7+e_8)$  & $19$\\
$\alpha_{107}$ & $2 \ \ 2 \ \ 3 \ \ 4 \ \ 3 \ \ 3 \ \ 2 \ \ 1$ & $-e_4+e_8$ & $20$\\
$\alpha_{108}$ & $1 \ \ 2 \ \ 3 \ \ 4 \ \ 4 \ \ 3 \ \ 2 \ \ 1$ &   $\frac12(-e_1+e_2-e_3+e_4+e_5+e_6+e_7+e_8)$  & $20$\\
$\alpha_{109}$ & $2 \ \ 2 \ \ 3 \ \ 4 \ \ 4 \ \ 3 \ \ 2 \ \ 1$ & $-e_3+e_8$ & $21$\\
$\alpha_{110}$ & $1 \ \ 2 \ \ 3 \ \ 5 \ \ 4 \ \ 3 \ \ 2 \ \ 1$ &   $\frac12(-e_1-e_2+e_3+e_4+e_5+e_6+e_7+e_8)$  & $21$\\
$\alpha_{111}$ & $2 \ \ 2 \ \ 3 \ \ 5 \ \ 4 \ \ 3 \ \ 2 \ \ 1$ & $-e_2+e_8$ & $22$\\
$\alpha_{112}$ & $1 \ \ 3 \ \ 3 \ \ 5 \ \ 4 \ \ 3 \ \ 2 \ \ 1$ &   $\frac12(e_1+e_2+e_3+e_4+e_5+e_6+e_7+e_8)$  & $22$\\
$\alpha_{113}$ & $2 \ \ 3 \ \ 3 \ \ 5 \ \ 4 \ \ 3 \ \ 2 \ \ 1$ & $e_1+e_8$ & $23$\\
$\alpha_{114}$ & $2 \ \ 2 \ \ 4 \ \ 5 \ \ 4 \ \ 3 \ \ 2 \ \ 1$ & $-e_1+e_8$ & $23$\\
$\alpha_{115}$ & $2 \ \ 3 \ \ 4 \ \ 5 \ \ 4 \ \ 3 \ \ 2 \ \ 1$ & $e_2+e_8$ & $24$\\
$\alpha_{116}$ & $2 \ \ 3 \ \ 4 \ \ 6 \ \ 4 \ \ 3 \ \ 2 \ \ 1$ & $e_3+e_8$ & $25$\\
$\alpha_{117}$ & $2 \ \ 3 \ \ 4 \ \ 6 \ \ 5 \ \ 3 \ \ 2 \ \ 1$ & $e_4+e_8$ & $26$\\
$\alpha_{118}$ & $2 \ \ 3 \ \ 4 \ \ 6 \ \ 5 \ \ 4 \ \ 2 \ \ 1$ & $e_5+e_8$ & $27$\\
$\alpha_{119}$ & $2 \ \ 3 \ \ 4 \ \ 6 \ \ 5 \ \ 4 \ \ 3 \ \ 1$ & $e_6+e_8$ & $28$\\
$\alpha_{120}$ & $2 \ \ 3 \ \ 4 \ \ 6 \ \ 5 \ \ 4 \ \ 3 \ \ 2$ & $e_7+e_8$ & $29$\\
  \hline
 \end{longtable}
\end{center}


\begin{center}
 \begin{longtable}{|c|l|l|}
 \caption{Arms and subhooks for roots $\alpha \in \Phi_8^+$.} \label{tab:armse8} \\

 \hline
 \hspace{-0.2cm} Root \hspace{-0.2cm} & Arm $a(\alpha)$ & Subhooks $h'(\beta)$ \\ \hline
 \endfirsthead

 \multicolumn{3}{c}{{\tablename\ \thetable{} (cont.)}} \\
 \hline
 \hspace{-0.2cm} Root \hspace{-0.2cm} & Arm $a(\alpha)$ & Subhooks $h'(\beta)$ \\ \hline
 \endhead

 \hline 
 \endfoot

 \endlastfoot

$\alpha_{1}$ &&\\
$\alpha_{2}$ &&\\
$\alpha_{3}$ &&\\
$\alpha_{4}$ &&\\
$\alpha_{5}$ &&\\
$\alpha_{6}$ &&\\

$\alpha_{9}$ & $3$ &\\

$\alpha_{10}$ & $4$ &\\

$\alpha_{11}$ & $4$ &\\

$\alpha_{12}$ & $5$ &\\

$\alpha_{13}$ & $6$ &\\

$\alpha_{16}$ & $1, 4$ &\\

$\alpha_{17}$ & $2, 3$ &\\

$\alpha_{18}$ & $2, 5$ &\\

$\alpha_{19}$ & $3, 5$ &\\

$\alpha_{20}$ & $4, 6$ &\\

$\alpha_{23}$ & $1, 2, 10$ &\\

$\alpha_{24}$ & $1, 5, 12$ &\\

$\alpha_{25}$ & $2, 3, 5$ &\\

$\alpha_{26}$ & $2, 6, 13$ &\\

$\alpha_{27}$ & $3, 6, 13$ &\\

$\alpha_{30}$ & $1, 2, 5, 9$ &\\

$\alpha_{31}$ & $1, 6, 9, 13$ &\\

$\alpha_{32}$ & $4, 10, 11, 12$ &\\

$\alpha_{33}$ & $2, 3, 6, 13$ &\\

$\alpha_{37}$ & $1, 4, 10, 12, 18$ &\\

$\alpha_{38}$ & $1, 2, 6, 9, 13$ &\\

$\alpha_{40}$ & $4, 6, 10, 11, 17$ &\\

$\alpha_{44}$ & $3, 9, 11, 16, 17, 19$ &\\

$\alpha_{45}$ & $1, 4, 6, 10, 16, 20$ &
$\{2, 24, 30\}$, $\{2, 27, 33\}$\\

$\alpha_{48}$ & $5, 12, 13, 18, 19, 20$ &\\

$\alpha_{51}$ & $3, 6, 9, 11, 16, 17, 23$ &\\

$\alpha_{52}$ & $1, 5, 12, 13, 18, 20, 26$ & \\

$\alpha_{57}$ & $3, 5, 9, 13, 19, 24, 25, 30$ & $\{6, 32, 40\}$, $\{1, 40, 45\}$\\

$\alpha_{63}$ & $4, 11, 12, 16, 19, 20, 24, 27, 31$ &\\

$\alpha_{69}$ & $2, 10, 17, 18, 23, 25, 26, 30, 32, 33 $ &\\

\hline

$\alpha_{7}$ && \\
$\alpha_{14}$ & $6$ &\\
$\alpha_{21}$ & $5, 7$ &\\
$\alpha_{28}$ & $4, 7, 12$ &\\
$\alpha_{34}$ & $2, 7, 10, 14$ &\\
$\alpha_{35}$ & $3, 7, 11, 14$ &\\
$\alpha_{39}$ & $1, 7, 9, 14, 16$ &\\
$\alpha_{41}$ & $2, 3, 7, 14, 21$ &\\
$\alpha_{46}$ & $1, 2, 7, 9, 14, 21$ &\\
$\alpha_{49}$ & $4, 7, 10, 11, 14, 17$ &\\

$\alpha_{53}$ & $1, 4, 7, 10, 14, 16, 23$ & $\{ 11, 21, 35 \}$, $\{ 2, 35, 41 \}$\\

$\alpha_{55}$ & $5, 7, 12, 18, 19, 25, 32$ &\\

$\alpha_{58}$ & $3, 7, 9, 11, 14, 16, 17, 23$ &\\

$\alpha_{59}$ & $1, 5, 7, 12, 18, 21, 28, 34$ & $\{ 9, 20, 31 \}$, $\{ 2, 31, 38 \}$,\\
&&$\{ 4, 38, 45 \}$\\

$\alpha_{61}$ & $6, 13, 14, 20, 21, 26, 27, 28$ &\\

$\alpha_{64}$ & $3, 5, 7, 9, 19, 21, 24, 25, 30$ & $\{ 11, 34, 49 \}$, $\{ 11, 38, 51 \}$, \\
&& $\{ 1, 49, 53 \}$\\

$\alpha_{66}$ & $1, 6, 13, 14, 20, 21, 26, 28, 34$ &\\

$\alpha_{70}$ & $4, 7, 11, 12, 16, 19, 24, 32, 37, 44$ & $\{ 10, 21, 34 \}$, $\{ 3, 34, 41 \}$,\\
&&$\{ 1, 41, 46 \}$\\

$\alpha_{71}$ & $3, 6, 9, 13, 14, 21, 27, 31, 33, 38$ & $\{ 7, 40, 49 \}$, $\{ 1, 49, 53 \}$, \\
&&$\{ 5, 49, 55 \}$, $\{ 1, 55, 59 \}$\\

$\alpha_{75}$ & $2, 7, 10, 17, 18, 23, 25, 30, 32, 37, 44$ &\\

$\alpha_{76}$ & $4, 6, 11, 14, 16, 20, 27, 28, 31, 35, 39$ & $\{ 12, 33, 48 \}$, $\{ 1, 48, 52 \}$, \\
&& $\{ 7, 48, 55 \}$, $\{ 3, 52, 57 \}$, \\
&& $\{ 1, 55, 59 \}$, $\{ 3, 59, 64 \}$ \\

$\alpha_{80}$ & 
$2, 6, 10, 14, 17, 23, 26, 33, 38, 40, 45, 51$ & $\{ 7, 48, 55 \}$, $\{ 1, 55, 59 \}$,\\
&&$\{ 3, 59, 64 \}$, $\{ 4, 64, 70 \}$\\

$\alpha_{82}$ & $5, 12, 13, 19, 20, 21, 24, 27, 28, 31, 35, 39$ &\\

$\alpha_{85}$ & $2, 5, 13, 18, 21, 25, 26, 30, 33, 34, 38, 41, 46$ & $\{ 12, 51, 63 \}$, $\{ 7, 63, 70 \}$,\\
&&$\{ 6, 70, 76 \}$\\

$\alpha_{89}$ & 
$4, 10, 12, 18, 20, 26, 28, 32, 34, 37, 40, 45,$ & $\{ 5, 51, 57 \}$, $\{ 5, 58, 64 \}$,\\
& $49, 53$ &$\{ 6, 64, 71 \}$\\

$\alpha_{93}$ & $ 3, 11, 17, 19, 25, 27, 32, 33, 35, 40, 41, 48,$&\\
& $49, 55, 61$&\\

$\alpha_{97}$ & $1, 9, 16, 23, 24, 30, 31, 37, 38, 39, 44, 45, 46,$ &\\
& $51, 52, 53$ &\\

\hline

$\alpha_8$ &&\\
$\alpha_{15}$ & $7$ &\\
$\alpha_{22}$ & $6, 8$ &\\
$\alpha_{29}$ & $5, 8, 13$ &\\
$\alpha_{36}$ & $4, 8, 12, 15$ &\\
$\alpha_{42}$ & $2, 8, 10, 15, 18$ &\\
$\alpha_{43}$ & $3, 8, 11, 15, 19$ &\\
$\alpha_{47}$ & $1, 8, 9, 15, 16, 22$ &\\
$\alpha_{50}$ & $2, 3, 8, 15, 22, 29$ &\\
$\alpha_{54}$ & $1, 2, 8, 9, 15, 22, 29$ &\\
$\alpha_{56}$ & $4, 8, 10, 11, 15, 17, 22$ &\\
$\alpha_{60}$ & $1, 4, 8, 10, 15, 16, 22, 23$ & $\{11, 29, 43 \}$, $\{2, 43, 50 \}$\\
$\alpha_{62}$ & $5, 8, 12, 15, 18, 19, 25, 32$&\\
$\alpha_{65}$ & $3, 8, 9, 11, 15, 16, 17, 22, 23$&\\

$\alpha_{67}$ & $1, 5, 8, 12, 15, 18, 24, 30, 37$& $\{ 19, 22, 43\}$, $\{ 2, 43, 50\}$,\\
&& $\{ 4, 50, 56\}$\\

$\alpha_{68}$ & $6, 8, 13, 20, 26, 27, 33, 40, 48$&\\

$\alpha_{72}$ & $3, 5, 8, 9, 15, 19, 24, 25, 29, 30$& $\{11,38,51\}$, $\{11,42,56\}$, \\
&& $\{ 7, 51, 58 \}$, $\{ 1, 56, 60 \}$ \\

$\alpha_{73}$ & $1, 6, 8, 13, 20, 22, 26, 29, 36, 42$& $\{ 9, 28, 39 \}$, $\{ 2, 39, 46 \}$,\\
&& $\{ 4, 46, 53 \}$, $\{ 5, 53, 59 \}$\\

$\alpha_{74}$ & $7, 14, 15, 21, 22, 28, 29, 34, 35, 36$&\\

$\alpha_{77}$ & $4, 8, 11, 12, 15, 16, 19, 24, 32, 37, 44$& $\{ 10, 29, 42 \}$, $\{ 3, 42, 50 \}$,\\
&& $\{ 1, 50, 54 \}$\\

$\alpha_{78}$ & $3, 6, 8, 9, 13, 22, 27, 29, 31, 33, 38$& $\{ 11, 42, 56 \}, \{ 11, 46, 58 \}$,\\
&& $\{ 1, 56, 60 \}$, $\{ 5, 56, 62 \}$, \\
&& $\{ 5, 58, 64 \}$, $\{ 1, 62, 67 \}$ \\

$\alpha_{79}$ & $1, 7, 14, 15, 21, 22, 28, 29, 34, 36, 42$&\\
$\alpha_{81}$ & $2, 8, 10, 15, 17, 18, 23, 25, 30, 32, 37, 44$&\\

$\alpha_{83}$ & $4, 6, 8, 11, 16, 20, 22, 27, 31, 40, 45, 51$& $\{ 10, 29, 42 \}$, $\{ 3, 42, 50 \}$, \\
&& $\{ 1, 50, 54 \}$, $\{ 12, 50, 62 \}$, \\
&& $\{ 1, 62, 67 \}$, $\{ 12, 58, 70 \}$, \\
&& $\{ 3, 67, 72 \}$ \\

$\alpha_{84}$ & $3, 7, 9, 14, 15, 21, 22, 29, 35, 39, 41, 46$& $\{ 8, 49, 56 \}$, $\{ 1, 56, 60 \}$,\\
&&$\{ 5, 56, 62 \}$, $\{ 1, 62, 67 \}$,\\
&&$\{ 6, 62, 68 \}$, $\{ 1, 68, 73 \}$\\

$\alpha_{86}$ & $2, 6, 8, 10, 17, 22, 23, 26, 33, 38, 40, 45, 51$& $\{ 15, 48, 62 \}$, $\{ 1, 62, 67 \}$, \\
&& $\{ 3, 67, 72 \}$, $\{ 18, 58, 75 \}$, \\
&& $\{ 4, 72, 77 \}$ \\

$\alpha_{87}$ & $5, 8, 12, 13, 19, 20, 24, 27, 31, 48, 52, 57, 63$& $\{ 15, 26, 42 \}$, $\{ 3, 42, 50 \}$,\\
&&$\{ 1, 50, 54 \}$, $\{ 4, 50, 56 \}$,\\
&&$\{ 1, 56, 60 \}$, $\{ 3, 60, 65 \}$\\

$\alpha_{88}$ & $ 4, 7, 11, 14, 15, 16, 22, 28, 35, 39, 49, 53, 58 $& $\{ 8, 34, 42 \}$, $\{ 3, 42, 50 \}$,\\
&&$\{ 1, 50, 54 \}$, $\{ 8, 55, 62 \}$,\\
&&$\{ 1, 62, 67 \}$, $\{ 6, 62, 68 \}$,\\
&&$\{ 3, 67, 72 \}$, $\{ 1, 68, 73 \}$,\\
&&$\{ 3, 73, 78 \}$\\

$\alpha_{90}$ & $2, 5, 8, 13, 18, 25, 26, 29, 30, 33, 38, 48, 52, 57$& $\{ 10, 43, 56 \}$,$\{ 1, 56, 60 \}$, \\
&& $\{ 3, 60, 65 \}$, $\{ 10, 64, 75 \}$, \\
&& $\{ 12, 65, 77 \}$,$\{ 6, 75, 80 \}$, \\
&& $\{ 6, 77, 83 \}$ \\

$\alpha_{91}$ & $ 2, 7, 10, 14, 15, 17, 22, 23, 34, 41, 46, 49, 53, 58 $& $\{ 8, 55, 62 \}$, $\{ 1, 62, 67 \}$, \\
&&$\{ 6, 62, 68 \}$, $\{ 3, 67, 72 \}$,\\
&&$\{ 1, 68, 73 \}$, $\{ 4, 72, 77 \}$,\\
&&$\{ 3, 73, 78 \}$, $\{ 4, 78, 83 \}$\\

$\alpha_{92}$ & $5, 7, 12, 15, 19, 21, 24, 28, 29, 35, 36, 39, 43,$& $\{ 13, 49, 61 \}$,$\{ 1, 61, 66 \}$,\\
& $47$& $\{ 8, 61, 68 \}$, $\{ 3, 66, 71 \}$ \\
&& $\{ 1, 68, 73 \}$, $\{ 4, 71, 76 \}$, \\
&& $\{ 3, 73, 78 \}$, $\{ 4, 78, 83 \}$ \\

$\alpha_{94}$ & $ 4, 8, 10, 12, 18, 20, 26, 32, 37, 40, 45, 48, 52, 63,$& $\{ 11, 29, 43 \}$, $\{ 1, 43, 47 \}$,\\
& $69$ &$\{ 2, 43, 50 \}$, $\{ 1, 50, 54 \}$,\\
&&$\{ 11, 54, 65 \}$, $\{ 5, 65, 72 \}$,\\
&&$\{ 6, 72, 78 \}$\\

$\alpha_{95}$ & $2, 5, 7, 15, 18, 21, 25, 29, 30, 34, 41, 46, 55, 59,$&$\{ 8, 49, 56 \}$, $\{ 1, 56, 60 \}$, \\
& $64$ & $\{ 3, 60, 65 \}$, $\{ 8, 61, 68 \}$, \\
&& $\{ 1, 68, 73 \}$, $\{ 8, 70, 77 \}$, \\
&& $\{ 3, 73, 78 \}$, $\{ 10, 71, 80 \}$, \\
&& $\{ 8, 76, 83 \}$, $\{ 8, 80, 86 \}$, \\
&& $\{ 8, 82, 87 \}$, $\{ 14, 77, 88 \}$ \\

$\alpha_{96}$ & $6, 13, 14, 20, 21, 22, 27, 28, 29, 31, 35, 36, 39, 43,$&\\
& $47$&\\

$\alpha_{98}$ & $ 3, 8, 11, 17, 19, 25, 27, 32, 33, 40, 44, 48, 51, 57,$& $\{ 9, 36, 47 \}$, $\{ 2, 47, 54 \}$,\\
& $63, 69$ &$\{ 4, 54, 60 \}$, $\{ 5, 60, 67 \}$,\\
&&$\{ 6, 67, 73 \}$\\

$\alpha_{99}$ & $4, 7, 10, 12, 15, 18, 28, 32, 34, 36, 37, 42, 49, 53,$& $\{ 8, 58, 65 \}$, $\{ 8, 61, 68 \}$, \\
             & $55, 59$& $\{ 5, 65, 72 \}$, $\{ 1, 68, 73 \}$, \\
&& $\{ 11, 66, 76 \}$, $\{ 8, 71, 78 \}$, \\
&& $\{ 2, 76, 80 \}$, $\{ 5, 76, 82 \}$, \\
&& $\{ 8, 76, 83 \}$, $\{ 14, 72, 84 \}$, \\
&& $\{ 2, 82, 85 \}$, $\{ 2, 83, 86 \}$, \\
&& $\{ 5, 83, 87 \}$, $\{ 2, 87, 90 \}$ \\

$\alpha_{100}$ & $2, 6, 13, 14, 21, 22, 26, 29, 33, 34, 38, 41, 42, 46,$& $\{ 20, 58, 76 \}$, $\{ 5, 76, 82 \}$, \\
& $50, 54$& $\{ 8, 76, 83 \}$, $\{ 5, 83, 87 \}$, \\
&& $\{ 7, 83, 88 \}$, $\{ 5, 88, 92 \}$ \\

$\alpha_{101}$ & $ 1, 8, 9, 16, 23, 24, 30, 31, 37, 38, 44, 45, 51, 52,$&\\
& $57, 63, 69$&\\

$\alpha_{102}$ & $3, 7, 11, 15, 17, 19, 25, 32, 35, 41, 43, 44, 49, 50,$&$\{ 9, 61, 71 \}$, $\{ 4, 71, 76 \}$,\\
& $55, 56, 62$& $\{ 8, 71, 78 \}$, $\{ 14, 67, 79 \}$,\\
&& $\{ 2, 76, 80 \}$, $\{ 5, 76, 82 \}$,\\
&& $\{ 4, 78, 83 \}$, $\{ 2, 82, 85 \}$,\\
&& $\{ 2, 83, 86 \}$, $\{ 5, 83, 87 \}$,\\
&& $\{ 4, 85, 89 \}$, $\{ 2, 87, 90 \}$,\\
&& $\{ 4, 90, 94 \}$ \\

$\alpha_{103}$ & $4, 6, 10, 14, 20, 22, 26, 28, 34, 36, 40, 42, 45, 49,$&$\{ 13, 58, 71 \}$, $\{ 8, 71, 78 \}$,\\
& $53, 56, 60$& $\{ 12, 71, 82 \}$, $\{ 7, 78, 84 \}$,\\
&& $\{ 2, 82, 85 \}$, $\{ 8, 82, 87 \}$,\\
&& $\{ 2, 87, 90 \}$, $\{ 7, 87, 92 \}$,\\
&& $\{ 2, 92, 95 \}$ \\

$\alpha_{104}$ & 
$ 1, 7, 9, 15, 16, 23, 24, 30, 37, 39, 44, 46, 53, 58,$ & $\{ 8, 66, 73 \}$, $\{ 3, 73, 78 \}$,\\
& $59, 64, 70, 75$ &$\{ 4, 78, 83 \}$, $\{ 2, 83, 86 \}$,\\
&& $\{ 5, 83, 87 \}$, $\{ 2, 87, 90 \}$,\\
&& $\{ 4, 90, 94 \}$, $\{ 3, 94, 98 \}$\\

$\alpha_{105}$ & $ 3, 6, 11, 14, 17, 22, 27, 33, 35, 40, 41, 43, 49, 50,$ & $\{ 13, 53, 66 \}$, $\{ 8, 66, 73 \}$,\\
& $51, 56, 58, 65$ & $\{ 7, 73, 79 \}$, $\{ 13, 70, 82 \}$,\\
&& $\{ 2, 82, 85 \}$, $\{ 8, 82, 87 \}$,\\
&& $\{ 4, 85, 89 \}$, $\{ 2, 87, 90 \}$,\\
&& $\{ 7, 87, 92 \}$, $\{ 4, 90, 94 \}$,\\
&& $\{ 2, 92, 95 \}$, $\{ 4, 95, 99 \}$\\

$\alpha_{106}$ & $5, 12, 13, 18, 20, 21, 26, 28, 29, 34, 36, 42, 48,$ & $\{ 8, 64, 72 \}$, $\{ 4, 72, 77 \}$,\\
& $52, 55, 59, 61, 66$ &$\{ 6, 72, 78 \}$, $\{ 2, 77, 81 \}$,\\
&&$\{ 4, 78, 83 \}$, $\{ 7, 78, 84 \}$,\\
&&$\{ 2, 83, 86 \}$, $\{ 4, 84, 88 \}$,\\
&&$\{ 2, 88, 91 \}$\\

$\alpha_{107}$ & $1, 6, 9, 14, 16, 22, 23, 31, 38, 39, 45, 46, 47, 51,$&$\{ 13, 70, 82 \}$, $\{ 2, 82, 85 \}$,\\
& $53, 54, 58, 60, 65$& $\{ 8, 82, 87 \}$, $\{ 4, 85, 89 \}$,\\
&& $\{ 2, 87, 90 \}$, $\{ 7, 87, 92 \}$,\\
&& $\{ 3, 89, 93 \}$, $\{ 4, 90, 94 \}$,\\
&& $\{ 2, 92, 95 \}$, $\{ 3, 94, 98 \}$,\\
&& $\{ 4, 95, 99 \}$, $\{ 3, 99, 102 \}$\\

$\alpha_{108}$ & $3, 5, 13, 19, 21, 25, 27, 29, 33, 35, 41, 43, 48, 50,$ & $\{ 6, 59, 66 \}$, $\{ 6, 67, 73 \}$,\\
& $55, 57, 62, 64, 72$ & $\{ 6, 70, 76 \}$, $\{ 7, 73, 79 \}$,\\
&&$\{ 2, 76, 80 \}$, $\{ 6, 77, 83 \}$,\\
&&$\{ 2, 83, 86 \}$, $\{ 7, 83, 88 \}$,\\
&&$\{ 12, 80, 89 \}$, $\{ 2, 88, 91 \}$,\\
&&$\{ 8, 89, 94 \}$, $\{ 7, 94, 99 \}$,\\
&&$\{ 6, 99, 103 \}$\\

$\alpha_{109}$ & $1, 5, 9, 13, 21, 24, 29, 30, 31, 38, 39, 46, 47, 52,$ & $\{ 8, 70, 77 \}$, $\{ 2, 77, 81 \}$,\\
& $54, 57, 59, 64, 66, 71$ &$\{ 6, 77, 83 \}$, $\{ 2, 83, 86 \}$,\\
&&$\{ 7, 83, 88 \}$, $\{ 12, 80, 89 \}$,\\
&&$\{ 2, 88, 91 \}$, $\{ 3, 89, 93 \}$,\\
&&$\{ 8, 89, 94 \}$, $\{ 3, 94, 98 \}$,\\
&&$\{ 7, 94, 99 \}$, $\{ 3, 99, 102 \}$,\\
&&$\{ 6, 99, 103 \}$, $\{ 3, 103, 105 \}$\\

$\alpha_{110}$ & $4, 11, 12, 19, 20, 27, 28, 32, 35, 36, 40, 43, 48,$ & $\{ 10, 57, 69 \}$, $\{ 7, 69, 75 \}$,\\
& $49, 55, 56, 61, 62, 68, 74$ &$\{ 6, 75, 80 \}$, $\{ 8, 75, 81 \}$,\\
&&$\{ 5, 80, 85 \}$, $\{ 6, 81, 86 \}$,\\
&&$\{ 5, 86, 90 \}$, $\{ 7, 86, 91 \}$,\\
&&$\{ 5, 91, 95 \}$, $\{ 6, 95, 100 \}$\\

$\alpha_{111}$ & $1, 4, 12, 16, 20, 24, 28, 31, 36, 37, 39, 45, 47,$&$\{ 6, 72, 78 \}$, $\{ 6, 75, 80 \}$,\\
& $52, 53, 59, 60, 63, 66, 67, 70$& $\{ 8, 75, 81 \}$, $\{ 7, 78, 84 \}$,\\
&& $\{ 5, 80, 85 \}$, $\{ 6, 81, 86 \}$,\\
&& $\{ 5, 86, 90 \}$, $\{ 7, 86, 91 \}$,\\
&& $\{ 11, 85, 93 \}$, $\{ 5, 91, 95 \}$,\\
&& $\{ 8, 93, 98 \}$, $\{ 6, 95, 100 \}$,\\
&& $\{ 7, 98, 102 \}$, $\{ 6, 102, 105 \}$,\\
&& $\{ 5, 105, 108 \}$\\

$\alpha_{112}$ & $ 2, 10, 17, 18, 25, 26, 32, 33, 34, 40, 41, 42, 48,$&\\
& $49, 50, 55, 56, 61, 62, 68, 74$&\\

$\alpha_{113}$ & $1, 2, 10, 18, 23, 26, 30, 34, 37, 38, 42, 45, 46,$ & $\{ 17, 82, 93 \}$, $\{ 8, 93, 98 \}$,\\
& $52, 53, 54, 59, 60, 66, 67, 73, 79$ &$\{ 7, 98, 102 \}$, $\{ 6, 102, 105 \}$,\\
&&$\{ 5, 105, 108 \}$, $\{ 4, 108, 110 \}$\\

$\alpha_{114}$ & $3, 9, 11, 16, 19, 24, 27, 31, 35, 39, 43, 44, 47,$ & $\{ 8, 75, 81 \}$, $\{ 6, 81, 86 \}$,\\
& $51, 57, 58, 63, 64, 70, 71, 76, 82$ &$\{ 5, 86, 90 \}$, $\{ 7, 86, 91 \}$,\\
&&$\{ 4, 90, 94 \}$, $\{ 5, 91, 95 \}$,\\
&&$\{ 4, 95, 99 \}$, $\{ 6, 95, 100 \}$,\\
&&$\{ 4, 100, 103 \}$, $\{ 5, 103, 106 \}$\\

$\alpha_{115}$ & $2, 3, 9, 17, 23, 25, 30, 33, 38, 41, 44, 46, 50, 51,$ & $\{ 10, 82, 89 \}$, $\{ 8, 89, 94 \}$,\\
& $54, 57, 58, 64, 65, 71, 72, 78, 84$ & $\{ 7, 94, 99 \}$, $\{ 6, 99, 103 \}$,\\
&&$\{ 5, 103, 106 \}, \hspace{-0.05cm} \{ 11, 106, 110 \}$, \hspace{-0.3cm} \\
&&$\{ 1, 110, 111 \}$\\

$\alpha_{116}$ & $4, 10, 11, 16, 17, 23, 32, 37, 40, 44, 45, 49, 51,$ & $\{ 2, 82, 85 \}$, $\{ 2, 87, 90 \}$,\\
& $53, 56, 58, 60, 63, 65, 70, 76, 77, 83, 88$ & $\{ 2, 92, 95 \}, \hspace{-0.05cm} \{ 2, 96, 100 \}$,\\
&&$\{ 12, 100, 106 \}, \hspace{-0.05cm} \{ 3, 106, 108 \}$, \hspace{-0.3cm} \\
&&$\{ 1, 108, 109 \}$\\

$\alpha_{117}$ & $5, 12, 18, 19, 24, 25, 30, 32, 37, 44, 48, 52, 55,$&$\{ 13, 88, 96 \}$, $\{ 2, 96, 100 \}$,\\
& $57, 59, 62, 63, 64, 67, 69, 70, 72, 75, 77, 81$& $\{ 4, 100, 103 \}$, $\{ 3, 103, 105 \}$, \hspace{-0.3cm} \\
&& $\{ 1, 105, 107 \}$\\

$\alpha_{118}$ & $6, 13, 20, 26, 27, 31, 33, 38, 40, 45, 48, 51, 52,$ & $\{ 2, 88, 91 \}$, $\{ 2, 92, 95 \}$,\\
& $57, 61, 63, 66, 68, 69, 71, 73, 76, 78, 82, 83, 87$ &$\{ 4, 95, 99 \}$, $\{ 3, 99, 102 \}$,\\
&&$\{ 1, 102, 104 \}$\\

$\alpha_{119}$ & $ 7, 14, 21, 28, 34, 35, 39, 41, 46, 49, 53, 55, 58,$&\\
& $59, 61, 64, 66, 70, 71, 75, 76, 80, 82, 85, 89, 93,$&\\
& $97$ &\\

$\alpha_{120}$ & $8, 15, 22, 29, 36, 42, 43, 47, 50, 54, 56, 60, 62,$ &\\
& $65, 67, 68, 72, 73, 74, 77, 78, 79, 81, 83, 84, 86,$&\\
& $87, 88$ &\\
  \hline
 \end{longtable}
\end{center}


\begin{center}
 \begin{longtable}{|c|l|l|c|}
 \caption{Positive roots in the root system $\Phi_4$ of type $F_4$.} \label{tab:rootsf4} \\

 \hline
 Root & Linear combination & Linear combination of $e_1, e_2, e_3, e_4$ & Height\\
      & of simple roots    && \\
      & $\alpha_1 \alpha_2 \alpha_3 \alpha_4$ && \\ \hline
 \endfirsthead

 \multicolumn{4}{c}{{\tablename\ \thetable{} (cont.)}} \\
 \hline
 Root & $\alpha_1 \alpha_2 \alpha_3 \alpha_4$ & Linear combination of $e_1, e_2, e_3, e_4$ & Height\\ \hline
 \endhead

 \hline 
 \endfoot

 \endlastfoot

$\alpha_{1}$ & $1 \ \ 0 \ \ 0 \ \ 0$ & $e_2-e_3$ & $1$\\
$\alpha_{2}$ & $0 \ \ 1 \ \ 0 \ \ 0$ & $e_3-e_4$ & $1$\\
$\alpha_{3}$ & $0 \ \ 0 \ \ 1 \ \ 0$ & $e_4$ & $1$\\
$\alpha_{4}$ & $0 \ \ 0 \ \ 0 \ \ 1$ & $\frac12 (e_1-e_2-e_3-e_4)$ & $1$\\
$\alpha_{5}$ & $1 \ \ 1 \ \ 0 \ \ 0$ & $e_2-e_4$ & $2$\\
$\alpha_{6}$ & $0 \ \ 1 \ \ 1 \ \ 0$ & $e_3$ & $2$\\
$\alpha_{7}$ & $0 \ \ 0 \ \ 1 \ \ 1$ & $\frac12 (e_1-e_2-e_3+e_4)$ & $2$\\
$\alpha_{8}$ & $1 \ \ 1 \ \ 1 \ \ 0$ & $e_2$ & $3$\\
$\alpha_{9}$ & $0 \ \ 1 \ \ 2 \ \ 0$ & $e_3+e_4$ & $3$\\
$\alpha_{10}$ & $0 \ \ 1 \ \ 1 \ \ 1$ & $\frac12 (e_1-e_2+e_3-e_4)$ & $3$\\
$\alpha_{11}$ & $1 \ \ 1 \ \ 2 \ \ 0$ & $e_2+e_4$ & $4$\\
$\alpha_{12}$ & $1 \ \ 1 \ \ 1 \ \ 1$ & $\frac12 (e_1+e_2-e_3-e_4)$ & $4$\\
$\alpha_{13}$ & $0 \ \ 1 \ \ 2 \ \ 1$ & $\frac12 (e_1-e_2+e_3+e_4)$ & $4$\\
$\alpha_{14}$ & $1 \ \ 2 \ \ 2 \ \ 0$ & $e_2+e_3$ & $5$\\
$\alpha_{15}$ & $1 \ \ 1 \ \ 2 \ \ 1$ & $\frac12 (e_1+e_2-e_3+e_4)$ & $5$\\
$\alpha_{16}$ & $0 \ \ 1 \ \ 2 \ \ 2$ & $e_1-e_2$ & $5$\\
$\alpha_{17}$ & $1 \ \ 2 \ \ 2 \ \ 1$ & $\frac12 (e_1+e_2+e_3-e_4)$ & $6$\\
$\alpha_{18}$ & $1 \ \ 1 \ \ 2 \ \ 2$ & $e_1-e_3$ & $6$\\
$\alpha_{19}$ & $1 \ \ 2 \ \ 3 \ \ 1$ & $\frac12 (e_1+e_2+e_3+e_4)$ & $7$\\
$\alpha_{20}$ & $1 \ \ 2 \ \ 2 \ \ 2$ & $e_1-e_4$ & $7$\\
$\alpha_{21}$ & $1 \ \ 2 \ \ 3 \ \ 2$ & $e_1$ & $8$\\
$\alpha_{22}$ & $1 \ \ 2 \ \ 4 \ \ 2$ & $e_1+e_4$ & $9$\\
$\alpha_{23}$ & $1 \ \ 3 \ \ 4 \ \ 2$ & $e_1+e_3$ & $10$\\
$\alpha_{24}$ & $2 \ \ 3 \ \ 4 \ \ 2$ & $e_1+e_2$ & $11$\\
  \hline
 \end{longtable}
\end{center}


\begin{center}
 \begin{longtable}{|c|l|l|}
 \caption{Arms and subhooks for roots $\alpha \in \Phi_4^+$.} \label{tab:armsf4} \\

 \hline
 \hspace{-0.2cm} Root \hspace{-0.2cm} & Arm $a(\alpha)$ & Subhooks $h'(\beta)$ \\ \hline
 \endfirsthead

 \multicolumn{3}{c}{{\tablename\ \thetable{} (cont.)}} \\
 \hline
 \hspace{-0.2cm} Root \hspace{-0.2cm} & Arm $a(\alpha)$ & Subhooks $h'(\beta)$ \\ \hline
 \endhead

 \hline 
 \endfoot

 \endlastfoot

$\alpha_{1}$ &&\\
$\alpha_{2}$ &&\\
$\alpha_{3}$ &&\\
$\alpha_{4}$ &&\\

$\alpha_{5}$ & $2$ &\\

$\alpha_{6}$ & $3$ &\\

$\alpha_{7}$ & $4$ &\\

$\alpha_{8}$ & $1, 3$ &\\

$\alpha_{9}$ & $3$ &\\

$\alpha_{10}$ & $2, 4$ &\\

$\alpha_{11}$ & $1, 3$ &\\

$\alpha_{12}$ & $1, 4, 7$ &\\

$\alpha_{13}$ & $3, 4, 7$ &\\

$\alpha_{14}$ & $2, 5, 6$ &\\

$\alpha_{15}$ & $1, 3, 4, 7$ &\\

$\alpha_{16}$ & $4, 7$ &\\

$\alpha_{17}$ & $2, 4, 5, 6, 8$ &\\

$\alpha_{18}$ & $1, 4, 7$ &\\

$\alpha_{19}$ & $3, 6, 7, 8, 9, 10$ &\\

$\alpha_{20}$ & $2, 4, 5, 10$ & $\{1, 13, 15\}$\\

$\alpha_{21}$ & $3, 4, 6, 7, 8, 10, 13$ &\\

$\alpha_{22}$ & $3, 7, 9, 11, 13$ & 
$\{6, 12, 17\}$, $\{4, 17, 20\}$\\

$\alpha_{23}$ & $2, 6, 9, 10, 13, 16$ &\\

$\alpha_{24}$ & $1, 5, 8, 11, 12, 14, 15$ &\\
  \hline
 \end{longtable}
\end{center}


\end{document}